\documentclass[11pt]{amsart}
\usepackage{amsxtra}
\usepackage{amssymb}
\usepackage[mathscr]{eucal}

\input xy
\xyoption{matrix}\xyoption{frame}

\theoremstyle{plain}
\newcommand{\cleqn}{\setcounter{equation}{0}}
\newcommand{\clth}{\setcounter{theorem}{0}}
\newcommand {\sectionnew}[1]{\section{#1}\cleqn\clth}

\newtheorem{theorem}{Theorem}[section]
\newtheorem{lemma}[theorem]{Lemma}
\newtheorem{definition-theorem}[theorem]{Definition-Theorem}
\newtheorem{proposition}[theorem]{Proposition}
\newtheorem{corollary}[theorem]{Corollary}
\newtheorem{definition}[theorem]{Definition}
\newtheorem{example}[theorem]{Example}
\newtheorem{examples}[theorem]{Examples}
\newtheorem{remark}[theorem]{Remark}
\newtheorem{conjecture}[theorem]{Conjecture}
\newtheorem{notation}[theorem]{Notation}

\newtheorem*{probA*}{Problem A}
\newtheorem*{probB*}{Problem B}
\newtheorem*{maintheorem*}{Main Theorem}
\newcommand \bth[1] { \begin{theorem}\label{t#1} }
\newcommand \ble[1] { \begin{lemma}\label{l#1} }

\newcommand \bpr[1] { \begin{proposition}\label{p#1} }
\newcommand \bco[1] { \begin{corollary}\label{c#1} }
\newcommand \bde[1] { \begin{definition}\label{d#1}\rm }
\newcommand \bex[1] { \begin{example}\label{e#1}\rm }
\newcommand \bre[1] { \begin{remark}\label{r#1}\rm }
\newcommand \bcj[1] { \begin{conjecture}\label{j#1}\rm }

\newcommand \bnota[1] { \begin{notation}\label{n#1}\rm }
\renewcommand {\eth} { \end{theorem} }
\newcommand {\ele} { \end{lemma} }

\newcommand {\epr} { \end{proposition} }
\newcommand {\eco} { \end{corollary} }
\newcommand {\ede} { \end{definition} }
\newcommand {\eex} { \end{example} }
\newcommand {\ere} { \end{remark} }
\newcommand {\ecj} { \end{conjecture} }

\newcommand {\enota} { \end{notation} }
\newcommand \thref[1]{Theorem \ref{t#1}}
\newcommand \leref[1]{Lemma \ref{l#1}}
\newcommand \prref[1]{Proposition \ref{p#1}}
\newcommand \coref[1]{Corollary \ref{c#1}}

\newcommand \deref[1]{Definition \ref{d#1}}
\newcommand \exref[1]{Example \ref{e#1}}

\theoremstyle{definition}
\newtheorem*{remark*}{Remark}
\newtheorem*{definition*}{Definition}
\newtheorem*{step1*}{Step 1}
\newtheorem*{step2*}{Step 2}
\newtheorem*{step3*}{Step 3}
\newtheorem*{step4*}{Step 4}
\newtheorem*{step5*}{Step 5}
\newtheorem*{step6*}{Step 6}

\newcommand \circled[1]{\xymatrix{#1 \save**\frm{o} \restore}}

\def \Cset {{\mathbb C}}
\def \KK {{\mathbb K}}

\def \Zset {{\mathbb Z}}

\def \Qset {{\mathbb Q}}

\def \AA {{\mathcal{A}}} 

\def \FF {{\mathcal{F}}}

\def \TT {{\mathcal{T}}}

\def \OO {{\mathcal{O}}}

\def \UU {{\mathcal{U}}}

\def \SS {{\mathcal{S}}}

\def \TT {{\mathcal{T}}} 

\def \qb {{\bf{q}}}
\def \rb {{\bf{r}}}
\def \xb {{\bf{x}}}
\def \yb {{\bf{y}}}

\def \ex {{\bf{ex}}}

\def \inv {{\bf{inv}}}
\def \de {\delta}
\def \al {\alpha}
\def \be {\beta}

\def \la {\lambda}

\def \Om {\Omega}
\def \ga {\gamma}
\def \de {\delta}
\def \Ga {\Gamma}

\def \sig {\sigma}

\def \ep {\epsilon}
\def \De {\Delta}
\def \sig{\sigma}
\def \rb { {\bf{r}} }

\def \mt  {\mapsto}

\def \hra {\hookrightarrow}

\def \sy  {\ast}                         
\def \bu  {\bullet}                 
\def \ci  {\circ}

\def \ol {\overline}
\def \wt {\widetilde}

\def \id { {\mathrm{id}} }

\def \sign { {\mathrm{sign}} }

\DeclareMathOperator \rk { {\mathrm{rk}} }

\def \h  {\mathfrak{h}}

\def \m  {\mathfrak{m}}
\def \b  {\mathfrak{b}}

\def \Mmn {M_{m,n}}

\DeclareMathOperator \Gr { {\mathrm{Gr}} }

\DeclareMathOperator \lt  { {\mathrm{lt}} }
\DeclareMathOperator \lc  { {\mathrm{lc}} }
\DeclareMathOperator \supp { {\mathrm{supp}} }

\DeclareMathOperator \Fract { {\mathrm{Fract}} }

\newcommand\kx{\KK^*}

\newcommand\HH{{\mathcal{H}}}
\newcommand\xh{X(\HH)}
\DeclareMathOperator \Spec {Spec}

\DeclareMathOperator \chr {char}
\newcommand \Znn {\Zset_{\ge 0}}

\newcommand \Hmax {\HH_{\max}}

\def\lab{{\boldsymbol \lambda}}
\def \gab {{\boldsymbol \gamma}}

\DeclareMathOperator \Lie {Lie}
\newcommand \Phat {\widehat{P}}
\newcommand \Rhat {\widehat{R}}
\newcommand \OMmn {\OO(M_{m,n}(\KK))}
\DeclareMathOperator \rank {rank}
\newcommand \hmax {\h_{\max}}
\newcommand{\xbtil}{\tilde{\xb}}
\newcommand{\ybtil}{\tilde{\yb}}
\newcommand{\Ytil}{\wt{Y}}
\newcommand{\Btil}{\wt{B}}
\DeclareMathOperator \Pcore {P{.}core}
\DeclareMathOperator \Cas {Cas}
\newcommand \kbar {\overline{\KK}}

\begin{document}
\title[Cluster algebra structures on Poisson nilpotent algebras]
{Cluster algebra structures on \\ 
Poisson nilpotent algebras}
\author[K. R. Goodearl]{K. R. Goodearl}
\address{
Department of Mathematics \\
University of California\\
Santa Barbara, CA 93106 \\
U.S.A.
}
\email{goodearl@math.ucsb.edu}
\author[M. T. Yakimov]{M. T. Yakimov}
\thanks{The research of K.R.G. was partially supported by NSF grants DMS-0800948 and DMS-1601184, 
and that of M.T.Y. by NSF grants DMS--1303038 and DMS--1601862.}
\address{
Department of Mathematics \\
Louisiana State University \\
Baton Rouge, LA 70803 \\
U.S.A.
}
\email{yakimov@math.lsu.edu}
\date{}
\begin{abstract}
Various coordinate rings of varieties appearing in the theory of Poisson Lie groups and Poisson homogeneous spaces belong to the 
large, axiomatically defined class 
of symmetric Poisson nilpotent algebras, e.g. coordinate rings of Schubert cells for symmetrizable Kac--Moody groups, affine charts of Bott-Samelson varieties, coordinate rings of double Bruhat cells (in the last case after a localization). 
We prove that every Poisson  nilpotent algebra satisfying a mild condition on certain scalars is canonically isomorphic to a cluster algebra which coincides with the corresponding upper cluster algebra, 
without additional localizations 
by frozen variables. The constructed cluster structure is compatible with the Poisson structure in the sense of Gekhtman, Shapiro and 
Vainshtein. All Poisson nilpotent algebras are proved to be equivariant Poisson Unique Factorization Domains. Their seeds 
are constructed from sequences of Poisson-prime elements for chains of Poisson UFDs; mutation matrices are effectively determined 
from linear systems in terms of the underlying Poisson structure. Uniqueness, existence, mutation, and 
other properties  are established for these sequences of Poisson-prime elements. \end{abstract}
\maketitle
\tableofcontents

\sectionnew{Introduction}
\label{intro}
\subsection{Cluster algebras and coordinate rings}
\label{cluster-coord}
Cluster algebras comprise a large, axiomatically defined class of algebras, introduced by Fomin and Zelevinsky 
\cite{FZ1}. They play a key role in representation theory, combinatorics, mathematical physics, algebraic and Poisson geometry.
We refer the reader to the books \cite{GSVbook,M} and the papers \cite{F,Ke,L,W} for overviews of some of the applications of cluster algebras.

A {\em{cluster algebra}} (of {\em{geometric type}}) is defined by iteratively mutating generating sets of a rational function field $\FF$ 
in $N$ variables over a base field $\KK$. The generating sets $(y_1, \ldots, y_N)$ of $\FF$ are called \emph{clusters}, and the process depends 
on a mutation matrix $\wt{B}$ associated to an initial cluster. 
The cluster algebra $\AA(\wt{B})_\KK$ is the $\KK$-subalgebra of $\FF$ generated by all cluster variables (both mutable and frozen), none of which are automatically inverted. We will write $\AA(\wt{B},\inv)_\KK$ for the localization of $\AA(\wt{B})_\KK$ obtained by inverting a subset $\inv$ of the frozen variables.
By the Fomin-Zelevinsky Laurent Phenomenon \cite{FZl}, the algebra $\AA(\wt{B})_\KK$ is a 
subalgebra of the {\em{upper cluster algebra}} $\overline{\AA}(\wt{B})_\KK$, obtained by intersecting in $\FF$ the mixed polynomial/Laurent polynomial rings of all clusters, in which the cluster variables are inverted but  the frozen variables are not. The analogous intersection of polynomial/Laurent polynomial rings in which a set $\inv$ of frozen variables is inverted will be denoted $\overline{\AA}(\wt{B},\inv)_\KK$. One of the main goals of the program of Fomin and Zelevinsky \cite{FZ1} is 

\begin{probA*}
Prove that the coordinate rings $\KK[V]$ of algebraic varieties $V$ that appear in Lie theory 
have canonical cluster algebra structures, i.e., $\KK[V] =  \AA(\wt{B})_\KK$ or, more generally, $\KK[V] = \AA(\wt{B},\inv)_\KK$.
\end{probA*}
The more general use of $\AA(\wt{B},\inv)_\KK$ is necessary due to the concrete properties of the varieties $V$ in question.
The consequences of this are that one can apply the cluster algebra structure of $\KK[V]$ to the study of the totally nonnegative part of $V$ 
and the canonical basis of $\KK[V]$, in the cases when they are defined. After the work of Berenstein, Fomin and Zelevinsky \cite{BFZ}, it became clear that 
another desirable property of cluster algebras is their equality with upper cluster algebras (when it holds), thus prompting

\begin{probB*}
Classify the coordinate rings $\KK[V]$ from Problem A for which the corresponding cluster algebras have the property 
$\AA(\wt{B})_\KK = \overline{\AA}(\wt{B})_\KK$ or, more generally,  
$\AA(\wt{B},\inv)_\KK = \overline{\AA}(\wt{B},\inv)_\KK$ for the same set $\inv$ that is used in the solution of Problem A
{\em{(}}equality of cluster and upper cluster algebras{\em{)}}.
\end{probB*}

Here is an illustration of Problems A and B: Berenstein, Fomin and Zelevinsky \cite{BFZ} proved that the coordinate rings of all double Bruhat cells $G^{w,u}$ 
in complex simple Lie groups $G$ possess upper cluster algebra structures: $\Cset[G^{w,u}] = \overline{\AA}(\wt{B},\inv)_\KK$. 
The problem of whether $\Cset[G^{w,u}] = \AA(\wt{B},\inv)_\KK$, or equivalently 
of whether Problem B has a positive solution in this case, was left open.

Gekhtman, Shapiro and Vainshtein \cite{GSV,GSVbook} proved that under mild assumptions a cluster algebra $\AA(\wt{B})_\KK$ has a Poisson algebra structure 
with respect to which all clusters are  \emph{log canonical}, which means that 
\[
\{y_l, y_j \} = \la_{lj} y_l y_j, \quad \la_{lj} \in \KK^*
\]
for the cluster variables $y_1, \ldots, y_N$ in each cluster of $\AA(\wt{B})_\KK$. Gekhtman, Shapiro and Vainshtein set up a program of proving that 
the coordinate rings of Poisson varieties $V$ in Lie theory possess natural cluster structures by constructing sufficiently many log-canonical 
clusters on $V$ (connected by mutations) and using their coordinates as cluster variables \cite{GSV1,GSVbook}. This has had remarkable success, and led to the construction 
of compatible upper cluster algebra structures for the Belavin-Drinfeld Poisson structures on $GL_n$ and their doubles, 
\cite{GSV1,GSV2}. However, similarly to \cite{BFZ}, the methods of \cite{GSV,GSVbook,GSV1,GSV2} rely on codimension $2$ arguments and only show that certain $\KK[V]$'s 
are isomorphic to upper cluster algebras rather than to cluster algebras.

There are three instances of positive solutions to Problem B. For acyclic cluster algebras this was shown in \cite{BFZ}. Muller \cite{Mu} proved that 
Problem B has a positive solution for the larger class of locally acyclic cluster algebras, but only in the form $\AA(\wt{B},\inv)_\KK = \overline{\AA}(\wt{B},\inv)_\KK$ where $\inv$ includes all frozen variables. (This extra localization does not allow to deduce that the coordinate rings of double Bruhat cells of type $A_n$ 
are cluster algebras.)
 Geiss, Leclerc and Schr\"oer \cite{GLSh1} proved that the coordinate rings of Schubert cells in symmetric Kac-Moody groups 
possess cluster algebra structures such that $\AA(\wt{B})_\KK = \overline{\AA}(\wt{B})_\KK$, using the representation theory of the related preprojective algebra.  
\subsection{Overview of results}
\label{overview}
The results obtained in this paper and the sequel \cite{GY-partII} are as follows: 

(I) In the present paper we show that, modulo a mild condition on scalars, every Poisson algebra $R$ in the very large axiomatically defined class of symmetric Poisson nilpotent algebras 
possesses a natural cluster algebra structure and for it $R= \AA(\wt{B})_\KK = \overline{\AA}(\wt{B})_\KK$.  For each such cluster algebra 
we construct a large family of clusters and explicit mutations between them.

(II) In the sequel \cite{GY-partII}, we prove that the coordinate rings of all Schubert cells and double Bruhat 
cells for a symmetrizable Kac-Moody group, equipped with their standard Poisson structures, are symmetric Poisson nilpotent algebras or localizations thereof. 
We then establish that, for the double Bruhat cells of every complex simple Lie group,
the Berenstein-Fomin-Zelevinsky cluster algebras equal the corresponding upper cluster algebras. 
We also prove such a result for the upper cluster algebra structures on the
double Bruhat cells of all symmetrizable Kac-Moody groups, constructed by Williams \cite{Wi}. 

Compared to the approaches in \cite{GLSh1,GSV,GSVbook,GSV1} to Problems A and B, our general results do not rely on 
concrete initial combinatorial data and apply to algebras coming from symmetrizable Kac-Moody groups, not only symmetric Kac-Moody groups.
Compared to \cite{Mu}, the above results prove equalities of cluster algebras and upper cluster algebras without an additional localization by frozen variables.
\subsection{Poisson nilpotent algebras and Main Theorem}
\label{P-CGL-intro}
Fix a base field $\KK$ throughout, with $\chr \KK = 0$, and consider everything to be done over $\KK$. In particular, the term ``algebra" will always mean a $\KK$-algebra, all automorphisms are assumed to be $\KK$-algebra automorphisms, and all derivations are assumed to be $\KK$-linear.

The Poisson nilpotent algebras (or Poisson-CGL extensions) with which we work are iterated \emph{Poisson-Ore extensions} $B[x; \sig, \de]_p$, where $B$ is a Poisson algebra, $B[x; \sig, \de]_p = B[x]$ is a polynomial ring, $\sig$ and $\de$ are suitable Poisson derivations on $B$, and
\[
\{x, b\} = \sig(b) x + \de(b), \; \; \forall b \in B
\]
(see \S\ref{p-ore}). For an iterated  Poisson-Ore extension
$$
R := \KK[x_1]_p [x_2; \sig_2, \delta_2]_p \cdots [x_N; \sig_N, \delta_N]_p
$$
and $k \in [0,N]$, denote
$$
R_k := \KK [ x_1, \dots, x_k ] = \KK[x_1]_p [x_2; \sig_2, \delta_2]_p \cdots [x_k; \sig_k, \delta_k]_p.
$$

\begin{definition*}
An iterated Poisson-Ore extension $R$ as above will be called a \emph{Poisson-CGL} \emph{extension} 
 if it is equipped with a rational Poisson action of a torus $\HH$ such that
\begin{enumerate}
\item[(i)] The elements $x_1, \ldots, x_N$ are $\HH$-eigenvectors.
\item[(ii)] For every $k \in [2,N]$, the map $\de_k$ on $R_{k-1}$ is locally nilpotent. 
\item[(iii)] For every $k \in [1,N]$, there exists $h_k \in \Lie \HH$ such that 
$\sig_k = (h_k \cdot)|_{R_{k-1}}$ and the $h_k$-eigen\-value of $x_k$, to be denoted by $\la_k$, is nonzero.
\end{enumerate}
There is a unique choice of $\HH$, maximal among tori acting faithfully and rationally on $R$ by Poisson automorphisms such that the above conditions hold, as we prove in \S\ref{maxtori}.

We say that the Poisson-CGL extension $R$ is \emph{symmetric} if it has the above properties (with respect to the same torus action) 
when the generators $x_1, \ldots, x_N$ are adjoined in reverse order. (This condition is of a different nature than the symmetry condition
for Kac-Moody groups; all symmetrizable Kac-Moody groups give rise to symmetric Poisson CGL extensions in different ways.)
\end{definition*}

Poisson-CGL extensions are semiclassical analogs of the {\em{Cauchon-Goodearl-Letzter extensions}}, appearing in the theory of quantum groups.
They have been extensively studied in recent years beginning with the papers \cite{Cau,GL,LLR}. CGL extensions provide a natural axiomatic 
framework for {\em{quantum nilpotent algebras}}: deformations of universal enveloping algebras of nilpotent Lie algebras at nonroots of unity \cite[p. 9697]{GYann}.
Because of this, we view Poisson-CGL extensions as an axiomatic framework for {\em{Poisson nilpotent algebras}}, i.e., {\em{Poisson algebra structures on 
polynomial rings which are semiclassical limits of quantum nilpotent algebras}}. We note that a Poisson CGL extension is not a nilpotent Lie algebra itself
under the Poisson bracket (Poisson algebras with such a property are very rare in the theory of Poisson Lie groups and Poisson homogeneous spaces).

Varieties that appear in the theory of Poisson Lie groups 
and Poisson homogeneous spaces provide numerous examples of Poisson-CGL extensions, see for example 
\cite{EL} for coordinate rings of Schubert cells and affine charts of Bott-Samelson varieties.
The coordinate rings of such varieties are semiclassical limits of quantized algebras of functions and the above axioms are Poisson
incarnations of the Levendorskii-Soibelman straightening law for quantum groups \cite[Proposition I.6.10]{BrGo}. 
A simple example of such a Poisson algebra, which is used as a running example through the paper, 
is the coordinate ring of the matrix affine Poisson space $M_{m,n}(\KK)$ of $m\times n$ matrices over $\KK$. It is the polynomial algebra 
$$
\OMmn = \KK[ t_{ij} \mid 1 \le i \le m, \; 1 \le j \le n]
$$
with the Poisson bracket
$$
\{ t_{ij}, t_{kl} \} = (\sign(k-i) + \sign(l-j)) t_{il} t_{kj},
$$
see \exref{OMmn} for details. This Poisson structure arises when the space $M_{m,n}(\KK)$ is identified with the open Schubert cell of the Grassmannian 
$\Gr(m, m+n)$ equipped with the standard Poisson structure, \cite{BGY,GSV}. 

The next theorem summarizes the main result of the paper:

\begin{maintheorem*}
Every symmetric Poisson-CGL extension $R$ such that $\la_l/\la_j \in \Qset_{>0}$ for all $l,j$ has a canonical structure of cluster algebra 
which coincides with its upper cluster algebra.
\end{maintheorem*}

A detailed formulation of this result is given in \thref{cluster}. Additional features of the full formulation of the theorem are: 
\begin{enumerate}
\item For each such Poisson-CGL extension $R$, the theorem constructs a large family of explicit seeds and mutations between them.
\item The Poisson structure on $R$ is compatible with the cluster structure 
in the sense of Gekhtman, Shapiro and Vainshtein.
\item The cluster variables in the constructed seeds are the unique homogeneous Poisson-prime elements 
of Poisson-CGL (sub)extensions that do not belong to smaller subextensions.
\item The mutation matrices of the constructed seeds can be effectively computed using linear systems of equations, 
coming from the Poisson structure.
\item For each generator $x_k$, an appropriate rescaling of it is a cluster variable in one of our seeds. In particular, the cluster 
variables in our finite collection of seeds generate the cluster algebra in question.
\item The above features also hold for the localization of $R$ with respect to any set $\inv$ of frozen variables; in particular, $R = \AA(\Btil,\inv)_\KK = \ol{\AA}(\Btil,\inv)_\KK$. 
\end{enumerate}
\subsection{Poisson-prime elements, Poisson-UFDs and proof of the Main Theorem}
\label{strategy-intro} Our proofs rely on arguments with Unique Factorization Domains in the Poisson algebra setting. 
Let $R$ be a noetherian Poisson domain (i.e., a Poisson algebra which is a noetherian integral domain). An element $p \in R$ is called a 
{\em{Poisson-prime element}} if any of the following  equivalent conditions is satisfied:
\begin{enumerate}
\item The ideal $(p)$ is a prime ideal and a Poisson ideal.
\item $p$ is a prime element of $R$ (in the commutative sense) such that $p \mid \{ p, - \}$. 
\item {[Assuming $\KK=\Cset$.]} $p$ is a prime element of $R$ and the zero locus $V(p)$ is a union of symplectic leaves of the maximal spectrum of $R$,
\end{enumerate}
see \S \ref{p-ufd} for details. The importance of the third formulation is that one can use Poisson geometry to classify Poisson-prime 
elements. 

Recall that by Nagata's lemma, a noetherian integral domain $R$ is a UFD if and only if every nonzero prime
ideal contains a prime element. We call a Poisson domain $R$ a {\em{Poisson-UFD}} if every nonzero Poisson-prime ideal contains a
Poisson-prime element. If $R$ is equipped with a Poisson action of a group $\HH$, we say that $R$ is an {\em{$\HH$-Poisson-UFD}}
if every nonzero $\HH$-Poisson-prime ideal of $R$ contains a Poisson-prime $\HH$-eigenvector. In the case of a $\KK$-torus $\HH$ 
acting rationally on $R$, the $\HH$-eigenvectors are precisely the nonzero homogeneous elements of $R$ with respect to its grading by the 
character group of $\HH$, and we will use the second terminology. 

Our proof of the Main Theorem is based on the following steps. The results in the individual steps are of independent interest and often 
admit wider generality than that of the Main Theorem.

\begin{step1*}
 Let $B$ be a noetherian $\HH$-Poisson-UFD for a $\KK$-torus $\HH$, and $B[x; \sig, \de]_p$  a Poisson-Ore extension with 
$\sig$ and $\de$ satisfying the properties (i)-(iii) in the definition of a Poisson-CGL extension. We prove that $B[x; \sig, \de]_p$
is also a noetherian $\HH$-Poisson-UFD and give an explicit classification of its homogeneous Poisson-prime elements in terms of  
those of $B$. These results appear in Theorems \ref{tPoissonLLR} and \ref{tBlistR}. Section \ref{PpPO} contains other details 
on the Poisson-UFD properties of the extension $B \subset B[x; \sig, \de]_p$.
\end{step1*}

\begin{step2*}
For each Poisson-CGL extension $R$ and $k \in [1,N]$, we prove that the Poisson algebra $R_k$ has a unique (up to rescaling)
homogeneous Poisson-prime element $y_k$ that does not belong to $R_{k-1}$. For the sequence 
\begin{multline}
\label{y-int}
\mbox{$y_1, \ldots, y_N$ of homogeneous 
Poisson-prime elements of} \\  \mbox{the algebras in the chain $R_1 \subset \ldots \subset R_N$}, 
\end{multline}
we prove that each $y_k$ is linear 
in $x_k$ with leading term $y_j x_k$ for some $j <k$ if $\de_k \neq 0$, and $y_k = x_k$ if $\de_k =0$. These facts are proved in 
\thref{mainPprime} and Section \ref{iterp-o} contains additional facts for the sequences of Poisson-primes \eqref{y-int}.
\end{step2*}

\begin{step3*}
Let $R$ be a symmetric Poisson-CGL extension. Each element of the set
\begin{equation}
\label{XiN-int}
\Xi_N := \{ \tau \in S_N \mid 
\tau([1,k]) \; \; \mbox{is an interval for all} 
\; \; k \in [2,N] \}
\end{equation}
gives rise to a presentation of $R$ as a Poisson-CGL extension with the generators of $R$  adjoined in the order
$x_{\tau(1)}, \ldots, x_{\tau(N)}$. The associated sequence of Poisson-prime elements from Step 2 will be denoted by
\begin{equation}
\label{tau-int}
y_{\tau, 1}, \ldots, y_{\tau, N}.
\end{equation}
For a symmetric Poisson-CGL extension $R$ and $j, k \in [1,N]$, the ``interval subalgebra'' 
$$
R_{[j,k]} := \KK[x_i \mid j \leq i \leq k]
$$ 
is an $\HH$-stable 
Poisson subalgebra of $R$ which is a Poisson-CGL extension in its own right with respect to the same torus action of $\HH$. 

In \thref{sym-Pprime} we express the elements of each sequence \eqref{tau-int} in terms of certain ``final Poisson-prime elements''
of the interval subalgebras $R_{[j,k]}$. Sections \ref{maxtorisymm} and \ref{Pprime-sym} contain further facts on symmetric Poisson 
CGL extensions and their interval subalgebras.
\end{step3*}

\begin{step4*}
 Next we link the clusters from the previous step by mutations. In Theorems \ref{tswapkk+1} and \ref{tswap2} we prove that if a Poisson-CGL extension has two presentations 
by adjoining the generators in the orders $x_1, \ldots, x_N$ and $x_1, \ldots, x_{k-1}, x_{k+1}, x_k, x_{k+2}, \ldots, x_N$, then the corresponding sequences 
of Poisson-prime elements \eqref{y-int} are either a permutation of each other or a one-step ``almost mutation'' of each other, 
meaning that the new elements $y$ are given by a mutation formula in terms of the old ones where the coefficient of one of the two monomials 
is not necessarily $1$ but a nonzero scalar. 

In \thref{PCGLmuta} we prove that the generators $x_1, \ldots, x_N$ of any symmetric Poisson-CGL extension can be rescaled so that 
the sequences of Poisson-prime elements \eqref{tau-int} corresponding to any pair $\tau, \tau' \in \Xi_N$ such that $\tau' = \tau (i,\, i+1)$ 
are either a permutation or a one-step mutation of each other. One of the upshots of Steps 2 and 4 is that the rescaling of each of the 
generators $x_k$ equals the cluster variable $y_{\tau, 1}$, for all $\tau \in \Xi_N$ with $\tau(1) = k$.    
Sections \ref{mCGL} and \ref{mut-sym} contain further details on the rescaling and mutations.
\end{step4*}

\begin{step5*}
 In Section \ref{Integr} we describe a method to control the size of the involved upper cluster algebras from above. In Theorems \ref{tdivision} and \ref{tPCGLalmostcluster} 
we prove that the intersection inside the fraction field of $R$ of all mixed polynomial/Laurent polynomial rings associated to the clusters indexed by $\Xi_N$ 
(without the frozen variables inverted) equals $R$.
\end{step5*}

\begin{step6*}
 Using the mutations from Step 4, we upgrade the clusters \eqref{tau-int} to seeds of cluster algebras and prove that 
the corresponding mutation matrices can be effectively computed using linear systems of equations, coming from the Poisson structure
(\thref{cluster}(a)). Using Steps 3--5 and the Laurent Phenomenon \cite{FZl} we get the chain of inclusions for the constructed cluster algebras $\AA(\wt{B})_\KK$
$$
R \subseteq \AA(\wt{B})_\KK \subseteq \overline{\AA}(\wt{B})_\KK \subseteq R, 
$$ 
which forces the equalities $R = \AA(\wt{B})_\KK = \overline{\AA}(\wt{B})_\KK$; this is carried out in Section \ref{main}. 
\end{step6*}

One of the consequences of the Main Theorem is that every symmetric Poisson-CGL extension with $\la_l/\la_j \in \Qset_{>0}$ possesses a canonical quantization given by the
quantum cluster algebra \cite{BZ} associated to the constructed Poisson cluster algebra. 
Recently, Y. Mi \cite{Mi} independently constructed a quantization of
each integral symmetric Poisson CGL extension that is a symmetric CGL
extension, and proved that in a certain sense this quantization is
unique. (Here integrality means that the elements $h_k$
in the definition in \S\ref{P-CGL-intro} belong to the cocharacter lattice of $\HH$.) The two quantizations can then be linked by 
the quantum cluster algebra structure on the latter constructed in \cite[Theorem 8.2]{GYbig}. 
\subsection{Notation and conventions}
\label{notaconv}
We write $\xh$ for the (rational) character group of a $\KK$-torus $\HH$, and we view $\xh$ as an additive group. Rational actions of $\HH$ on an algebra $R$ are equivalent to $\xh$-gradings \cite[Lemma II.2.11]{BrGo}, and we will use the terminology of gradings whenever convenient. In particular, when $R$ is equipped with a rational $\HH$-action, its homogeneous elements (with respect to the $\xh$-grading) are the $\HH$-eigenvectors together with $0$. The degree of a nonzero homogeneous element $r \in R$ will be denoted $\chi_r$; this is just its $\HH$-eigenvalue.

For any positive integer $N$, we view the elements of $\Zset^N$ as column vectors, and we write $\{ e_1, \dots, e_N \}$ for the standard basis of $\Zset^N$. The row vector corresponding to any $f \in \Zset^N$ is the transpose of $f$, denoted $f^T$. If $f = (m_1,\dots,m_N)^T$, denote
\begin{equation}
\label{f+-}
[f]_+ := \sum_{j=1}^N \max(m_j,0) e_j \qquad \text{and} \qquad [f]_-  := \sum_{j=1}^N \min(m_j,0) e_j .
\end{equation}
 We identify permutations in $S_N$ with their corresponding permutation matrices in $GL_N(\Zset)$, so that each $\tau \in S_N$ acts on $\Zset^N$ by $\tau(e_i) = e_{\tau(i)}$, for all $i \in [1,N]$.  Corresponding to any skew-symmetric matrix $\qb = (q_{kj}) \in M_N(\KK)$ is a skew-symmetric bicharacter
\begin{equation}
\label{Om}
\Om_\qb : \Zset^N \times \Zset^N \rightarrow \KK, \quad \text{given by} \quad \Om_\qb(e_k, e_j) = q_{kj}, \; \; \forall j,k \in [1,N].
\end{equation}
Given $f = (m_1,\dots,m_N)^T \in \Zset^N$ and an $N$-tuple $\xb = (x_1,\dots,x_N)$ of elements from a commutative ring, we denote
\begin{equation}
\label{x^f}
\xb^f := \prod_{j=1}^N x_j^{m_j}.
\end{equation}

We recall definitions and establish some auxiliary results for Poisson algebras and cluster algebras in Sections \ref{Palgs} and \ref{appendix}, respectively. Here we just mention that since our aim is to produce cluster algebra structures on $\KK$-algebras, we build cluster algebras $\AA(\xbtil,\Btil)_\KK$ over the field $\KK$ directly, rather than first building cluster algebras over an integral (semi)group ring and then extending scalars to $\KK$. Moreover, we work exclusively with cluster algebras of geometric type.

\sectionnew{Poisson algebras}
\label{Palgs}

\subsection{Poisson algebras and Poisson ideals}
\label{p-alg}

We recall that a \emph{Poisson algebra} is a commutative algebra $R$ equipped with a Poisson bracket $\{-,-\}$, that is, a Lie bracket which is also a derivation in each variable (for the associative product). For $a\in R$, the derivation $\{a,-\}$ is called the \emph{Hamiltonian associated to $a$}. The Poisson bracket on $R$ induces unique Poisson brackets on any quotient of $R$ modulo  a \emph{Poisson ideal}, meaning an ideal $I$ such that $\{R,I\} \subseteq I$, and on any localization of $R$ with respect to a multiplicative set (e.g.,
\cite[Proposition 1.7]{Loo}). A \emph{Poisson automorphism} of $R$ is any algebra automorphism which preserves the Poisson bracket. We use the term \emph{Poisson action} to refer to an action of a group on $R$ by Poisson automorphisms. 

\begin{examples}  \label{PaffPtorus}  {\rm
A \emph{Poisson-Weyl algebra} is  a polynomial algebra
$\KK[x_1,\dots,x_{2n}]$ equipped with the Poisson bracket such that
$$\{x_i,x_j\} = \{x_{n+i}, x_{n+j}\} =0 \;\; \text{and} \;\; \{x_i,x_{n+j}\} = \delta_{ij}\,, \;\; \forall\, i,j \in [1,n].$$ 

Given any skew-symmetric matrix $(q_{kj}) \in M_N(\KK)$, there are compatible Poisson brackets on the polynomial algebra $\KK[x_1,\dots,x_N]$ and the Laurent polynomial algebra $\KK[x_1^{\pm1},\dots,x_N^{\pm1}]$ such that $\{x_k, x_j\} = q_{kj} x_k x_j$ for all $k$, $j$. The Poisson algebras $\KK[x_1,\dots,x_N]$ and $\KK[x_1^{\pm1},\dots,x_N^{\pm1}]$ are known as a \emph{Poisson affine space algebra} and a \emph{Poisson torus}, respectively.
}\end{examples}

The \emph{Poisson center} of a Poisson algebra $R$ is the subalgebra
$$Z_p(R) := \{ z \in R \mid \{z,-\} = 0 \}.$$
Its elements are sometimes called \emph{Casimirs}, in which case $Z_p(R)$ is denoted $\Cas(R)$.
An element $c \in R$ is said to be \emph{Poisson-normal} if the principal ideal $Rc$ is a Poisson ideal, that is, if $\{c,R\} \subseteq Rc$. If also $c$ is a non-zero-divisor, then $\{c,-\}$ determines a derivation $\partial_c$ on $R$ such that
\begin{equation}
\label{defdelc}
\{c,a\} = \partial_c(a) c, \; \; \forall a \in R.
\end{equation}
Moreover, it follows from the Jacobi identity for $\{ -,- \}$ that
$$
\partial_c( \{a,b\} ) = \{ \partial_c(a) ,b\} + \{a, \partial_c(b) \}, \; \; \forall\, a,b \in R,
$$
so that $\partial_c$ is also a derivation with respect to the Poisson bracket.

Given an arbitrary ideal $J$ in $R$, there is a largest Poisson
ideal contained in $J$, which, following \cite{BrGr}, we call the \emph{Poisson
core of $J$}; we shall denote it $\Pcore(J)$. The \emph{Poisson primitive
ideals} of
$R$ are the Poisson cores of the maximal
ideals of $R$ (in
\cite{Oh.symp}, these are called {\it symplectic ideals\/}). 
A \emph{Poisson-prime} ideal is any proper Poisson ideal $P$ of $R$ such that $(IJ \subseteq P \implies I \subseteq P \; \text{or} \; J \subseteq P)$ for all Poisson ideals $I$ and $J$ of $R$.
The above concepts are related by the following Poisson version of  \cite[Lemma 1.1]{GpDixMo}; we repeat the short arguments for convenience.

\ble{Pprimeinfo}
Let $R$ be a Poisson algebra.

{\rm(a)} $\Pcore(P)$ is prime for all prime ideals $P$ of $R$.

{\rm(b)} Every Poisson primitive ideal of $R$ is prime.

{\rm(c)} Every prime ideal minimal over a Poisson ideal is a Poisson ideal.

{\rm(d)} If $R$ is noetherian, every Poisson-prime ideal of $R$ is prime.

{\rm(e)} If $R$ is affine over $\KK$, every Poisson-prime ideal of $R$ is an intersection of Poisson primitive ideals.
\ele

\begin{remark*}
If $R$ is a noetherian Poisson algebra, \leref{Pprimeinfo}(d) implies that the Poisson-prime ideals
 in $R$ are precisely the ideals which are both
Poisson ideals and prime ideals; in that case, the hyphen in the term 
``Poisson-prime'' becomes unnecessary.
\end{remark*}

\begin{proof}
(a) Since $\Pcore(P)$ is the largest ideal contained in $P$ and stable under all the derivations $\{a,-\}$, this follows from \cite[Lemma 3.3.2]{Dixbook}.  

(b) and (c) are immediate from (a).

(d) Let $Q$ be a Poisson-prime ideal of $R$. There exist prime ideals $Q_1,\dots,Q_n$
minimal over $Q$ such that $Q_1Q_2\cdots Q_n\subseteq Q$. The $Q_i$ are Poisson ideals by (c), so the Poisson-primeness of $Q$ implies that some $Q_i = Q$.

(e) Let $Q$ be a Poisson-prime ideal of $R$. Since $Q$ is prime (by (d)), the Nullstellensatz implies that $Q$ is an intersection of maximal ideals $M_i$. Consequently, $Q = \bigcap_i \Pcore(M_i)$.
\end{proof}

If $\KK = \kbar$ and the Poisson algebra $R$ is the coordinate ring of an affine algebraic variety $V$, then $V$ may be partitioned into \emph{symplectic cores} \cite[\S3.3]{BrGr} parametrized by the Poisson primitive ideals of $R$, where the symplectic core corresponding to a Poisson primitive ideal $P$ is the set
$$\{ x \in V \mid \Pcore(\m_x) = P \}.$$
(Here $\m_x$ denotes the maximal ideal of $R$ corresponding to $x$.) In case $\KK = \Cset$, the variety $V$ is a union of smooth complex Poisson manifolds and, as such, is partitioned into symplectic leaves (see \cite[\S3.5]{BrGr}). In this case, it follows from \cite[Lemma 3.5]{BrGr} that the Zariski closure of any symplectic leaf $L$ of $V$ is a union of symplectic cores. In fact, the lemma shows that the defining ideal $I(\ol{L})$ is a Poisson primitive ideal $P$, equal to $\Pcore(\m_x)$ for any $x \in L$. Consequently, $\ol{L}$ is the union of the symplectic cores corresponding to the Poisson primitive ideals containing $P$.

These observations lead to geometric ways to verify that certain ideals are Poisson ideals, as follows. 

\bpr{Pvleaves}
Let $\KK = \Cset$ and $R$ the coordinate ring of an affine algebraic variety $V$. Suppose $R$ is a Poisson algebra and $I$ is a radical ideal of $R$. Then $I$ is a Poisson ideal if and only if the subvariety $V(I)$ is a union of symplectic leaves, if and only if $V(I)$is a union of Zariski closures of symplectic leaves.
\epr

\begin{proof} Suppose first that $V(I)$ is a union of closures of symplectic leaves $L_j$. By \cite[Lemma 3.5]{BrGr}, $I(\ol{L_j})$ is a Poisson ideal of $R$, and therefore $I = \bigcap_j I(\ol{L_j})$ is a Poisson ideal.

Now assume that $I$ is a Poisson ideal. By \leref{Pprimeinfo}(e), $I$ equals the intersection of the Poisson primitive ideals $P_j$ that contain $I$. If $x \in V(I)$, then from $\m_x \supseteq I$ we get $\Pcore(\m_x) \supseteq I$, whence $\Pcore(\m_x) = P_j$ for some $j$. By \cite[Lemma 3.5]{BrGr}, the closure of the symplectic leaf $L$ containing $x$ satisfies $I(\ol{L}) = P_j$, whence $L \subseteq V(P_j) \subseteq V(I)$. Therefore $V(I)$ is a union of symplectic leaves.

Finally, if $V(I)$ is a union of symplectic leaves, then, being closed in $V$, it is also a union of closures of symplectic leaves.
\end{proof}

A version of \prref{Pvleaves} over more general fields, with symplectic leaves replaced by symplectic cores, is available under additional hypotheses, such as uncountability of the base field, as we now show. 

\ble{BLSMarg}
Let $T$ be an affine Poisson algebra. Assume that $\KK$ is uncountable, $T$ is a domain, and $Z_p(\Fract T) = \KK$. Then there is a countable set $\SS$ of nonzero Poisson prime ideals of $T$ such that every nonzero Poisson prime ideal of $T$ contains a member of $\SS$.
\ele

\begin{proof} Let $\SS$ be the set of those nonzero Poisson prime ideals of $T$ which do not properly contain any nonzero Poisson prime ideal. Due to the descending chain condition on prime ideals in $T$, every nonzero Poisson prime ideal of $T$ contains a member of $\SS$. The argument in the proof of \cite[Theorem 3.2]{BLSM} shows that $\SS$ is countable. (The assumption $\KK= \Cset$ in that theorem is not used in this part of the proof.)
\end{proof}

The key argument of the following result was communicated to us by Jason Bell; we thank him for permission to include it here.

\ble{P.Jacobson.cond}
If $R$ is an affine Poisson algebra and $\KK$ is uncountable, then
\begin{equation}  \label{Pprimcap}
P = \bigcap\, \{ \m \in \max R \mid \Pcore(\m) = P \}, \;\; \forall \, \text{Poisson primitive ideals} \;\, P \;\, \text{of} \;\, R.
\end{equation}
\ele

\begin{proof}  It suffices to consider the case where the Poisson primitive ideal $P$ is zero. Hence, there is a maximal ideal $\m_0$ of $R$ such that $\Pcore(\m_0) = 0$.

Set $\ol{T} := R \otimes_\KK \kbar$, which is an affine Poisson $\kbar$-algebra, identify $R$ with its image in $\ol{T}$, and observe that $\ol{T}$ is integral over $R$. We will use several standard results about prime ideals in integral extensions, such as \cite[Theorems 44, 46, 48]{Kap}. In particular, since $R$ and $\ol{T}$ have the same finite Krull dimension, it follows that any non-minimal prime ideal of $\ol{T}$ has nonzero contraction to $R$.

There is a prime ideal $Q_0$ of $\ol{T}$ such that $Q_0 \cap R = 0$. After possibly shrinking $Q_0$, we may assume that it is a minimal prime, and thus a Poisson prime (\leref{Pprimeinfo}(c)). Hence, the domain $T := \ol{T}/Q_0$ is an affine Poisson $\kbar$-algebra. We identify $R$ with its image in $T$.

There is a maximal ideal $\m_0^* \in \max T$ such that $\m_0^* \cap R = \m_0$. Since $\Pcore(\m_0^*)$ is a Poisson prime ideal of $T$, the contraction $\Pcore(\m_0^*) \cap R$ is a Poisson prime ideal of $R$ contained in $\m_0$, whence $\Pcore(\m_0^*) \cap R = 0$. Consequently, $\Pcore(\m_0^*) =0$, and so $0$ is a Poisson primitive ideal of $T$. 

By \cite[Proposition 1.10]{Oh.symp}, $0$ is a Poisson-rational ideal of $T$, meaning that the Poisson center of $\Fract T$ is (algebraic over) $\kbar$. Hence, \leref{BLSMarg} provides us with a countable set $\{Q_1,Q_2,\dots \}$ of nonzero Poisson prime ideals of $T$ such that every nonzero Poisson prime ideal of $T$ contains some $Q_j$. Moreover, the ideals $P_j := Q_j \cap R$ are nonzero Poisson prime ideals of $R$.

We claim that any nonzero Poisson prime ideal $P$ of $R$ contains some $P_j$. There is a prime ideal $Q$ of $T$ lying over $P$, and after possibly shrinking $Q$, we may assume that $Q$ is minimal over $PT$. Since $PT$ is a Poisson ideal, so is $Q$. Then $Q \supseteq Q_j$ for some $j$, whence $P \supseteq P_j$, validating the claim. 

Choose a nonzero element $f_j \in P_j$ for each $j$, let $X$ be the multiplicative set generated by the $f_j$, and set $R' := R[X^{-1}]$. Since $R'$ is a countably generated $\KK$-algebra and $\KK$ is uncountable, $R'$ satisfies the Nullstellensatz over $\KK$ (e.g., \cite[Proposition II.7.16]{BrGo}). This means that the Jacobson radical of $R'$ is zero and $R'/\m'$ is 
algebraic over $\KK$ for all $\m' \in \max R'$. Consequently, $\m' \cap R \in \max R$ for any $\m' \in \max R'$. Moreover, $\m' \cap R$ cannot contain any $P_j$, and so the Poisson prime ideal $\Pcore(\m' \cap R)$ must be zero. Since $\bigcap \{ \m'\cap R \mid \m' \in \max R' \} = 0$, \eqref{Pprimcap} is proved.
\end{proof}

The condition \eqref{Pprimcap} also holds if $R$ is affine and all Poisson primitive ideals of $R$ are locally closed points of the Poisson-prime spectrum of $R$. This follows from the proof of \cite[Theorem 1.5]{GpDixMo}.

\bpr{Pvcores}
Let $\KK = \kbar$ and $R$ the coordinate ring of an affine algebraic variety $V$. Suppose $R$ is a Poisson algebra satisfying \eqref{Pprimcap}, and let $I$ be a radical ideal of $R$. Then $I$ is a Poisson ideal if and only if  $V(I)$ is a union of symplectic cores.
\epr

\begin{proof}
Suppose first that $I$ is a Poisson ideal. By \leref{Pprimeinfo}(e), $I$ equals the intersection of the Poisson primitive ideals $P_j$ that contain $I$. Obviously any point in the symplectic core $C_j$ corresponding to $P_j$ lies in $V(I)$. Conversely, if $x \in V(I)$, then from $\m_x \supseteq I$ we get $\Pcore(\m_x) \supseteq I$, whence $\Pcore(\m_x) = P_j$ and $x\in C_j$, for some $j$. Thus, $V(I) = \bigcup_j C_j$.

Now asume that $V(I)$ is a union of symplectic cores $C_j$ corresponding to Poisson primitive ideals $P_j$. Then
$$I = \bigcap_j \biggl( \bigcap \{ \m_x \mid x\in C_j \} \biggr).$$
Because of \eqref{Pprimcap}, $\bigcap \{ \m_x \mid x\in C_j \} = P_j$ for all $j$. Therefore $I = \bigcap_j P_j$, whence $I$ is a Poisson ideal.
\end{proof}

\subsection{Poisson polynomial rings}
\label{p-ore}

\bde{p-deriv}
A \emph{Poisson derivation} on a Poisson algebra $B$ is any $\KK$-linear map $\sig$ on $B$ which is a derivation with respect to both the associative multiplication and the Poisson bracket, that is,
\begin{equation}
\label{pderiv}
\sig(ab) = \sig(a) b + a \sig(b) \quad \text{and} \quad \sig(\{a,b\}) = \{\sig(a), b\} + \{a, \sig(b)\}, \; \; \forall a,b \in B.
\end{equation}
For instance, the maps $\partial_c$ defined in \eqref{defdelc} are Poisson derivations.

Assume that $\sig$ is a Poisson derivation on $B$. 
Following the terminology of \cite[\S 1.1.2]{Dumas}, a \emph{Poisson $\sig$-derivation} on $B$ is any derivation $\de$ (with respect to the associative multiplication) such that
\begin{equation}
\label{p-sig-deriv}
\de(\{a,b\}) = \{\de(a), b\} + \{a, \de(b)\} + \sig(a) \de(b) - \de(a) \sig(b), \; \; \forall a,b \in B.
\end{equation}
For any $c \in B$, the map
$$
a \mapsto \{c,a\} - \sig(a) c, \; \; \forall a \in B,
$$
is a Poisson $\sig$-derivation, and we refer to the ones of this form as \emph{inner Poisson $\sig$-der\-i\-vations} of $B$.
\ede

\ble{p-ore-construct}
{\rm \cite[Theorem 1.1]{Oh}}
Let $B$ be a Poisson algebra, $\sig$ a Poisson derivation on $B$, and $\de$ a Poisson $\sig$-derivation on $B$. The Poisson bracket on $B$ extends uniquely to a Poisson bracket on the polynomial algebra $B[x]$ such that
\begin{equation}
\label{p-ore-bracket}
\{x, b\} = \sig(b) x + \de(b), \; \; \forall b \in B.
\end{equation}
\ele

\bde{p-ore-def}
Suppose that $B$ is a Poisson algebra, $\sig$ a Poisson derivation on $B$, and $\de$ a Poisson $\sig$-derivation on $B$. The polynomial algebra $B[x]$, equipped with the Poisson bracket of \leref{p-ore-construct}, is just called a \emph{Poisson polynomial ring} in \cite{Oh}. Here we will say that $B[x]$ is a \emph{Poisson-Ore extension} of $B$, as in \cite{Dumas}. The Poisson bracket on $B[x]$ extends uniquely to one on the Laurent polynomial algebra $B[x^{\pm1}]$, and we will call the latter Poisson algebra a \emph{Poisson-Laurent extension} of $B$. To express that $B[x]$ and $B[x^{\pm1}]$ are Poisson algebras of the types just defined, we use notation analogous to that for Ore and skew-Laurent extensions, namely
$$B[x; \sig, \de]_p \quad \text{and} \quad B[x^{\pm1}; \sig, \de]_p \,.$$
In case $\de$ is identically zero, we omit it from the notation, writing just $B[x; \sig]_p$ and $B[x^{\pm1}; \sig]_p$.
\ede

The converse of \leref{p-ore-construct} is worth noting; it may be phrased as follows: If  a
polynomial ring $R=B[x]$ supports a Poisson bracket such that $B$ is a
Poisson subalgebra and $\{x,B\} \subseteq Bx+B$, then $R=
B[x;\sig,\de]_p$ for suitable $\sig$ and $\de$.

Basic examples of iterated Poisson-Ore extensions of the form
$$\KK[x_1] [x_2; \sig_2, \delta_2]_p \cdots [x_N; \sig_N, \delta_N]_p$$
are Poisson-Weyl algebras and Poisson affine space algebras (Examples \ref{PaffPtorus}).

\subsection{Differentials of torus actions}
\label{dHaction}
Suppose that $R$ is an algebra equipped with a rational action of a torus $\HH$ (by algebra automorphisms). The differential of this action provides an action of $\h := \Lie \HH$ on $R$ by derivations, as discussed in \cite[\S1.2]{GLaun}, where it is also noted that the $\HH$-stable $\KK$-subspaces of $A$ coincide with the $\h$-stable subspaces. It follows from this discussion and \cite[Lemma 1.3]{GLaun} that 
\begin{gather}
\text{The\;} \HH\text{-eigenspaces of\;} R\; \text{coincide with the\;} \h\text{-eigenspaces;}
\label{hHeigenspaces}  \\
\text{The\;} \HH\text{-eigenvectors in\;} R\; \text{coincide with the\;} \h\text{-eigenvectors;}
\label{hHeigenvectors}  \\
\text{The\;} \h\text{-action on} \; R \; \text{commutes with the} \; \HH\text{-action};
\label{hHcomm}  \\
\h.a=0, \; \; \forall \text{\;homogeneous\;} a\in A \text{\;of degree\;} 0.
\label{h.0}
\end{gather}

\ble{rat-p-action}
{\rm\cite[Lemma 1.4]{GLaun}}
Suppose a Poisson algebra $R$ is equipped with a rational Poisson action of a torus $\HH$. Then $\Lie \HH$ acts on $R$ by Poisson derivations.
\ele

\sectionnew{Cluster algebras and Poisson cluster algebras}
\label{appendix}

In this section, we collect basic definitions and notation concerning cluster algebras and compatible Poisson structures on them, together with some auxiliary results. We consider exclusively cluster algebras of geometric type, recall \cite[Definition 5.7]{FZ1}, and we do not invert frozen variables unless specifically indicated.

\subsection{Cluster algebras}
\label{clusternota}

Fix positive integers $n \le N$ and a subset $\ex \subset [1,N]$ of cardinality $n$. The indices in $\ex$ are called \emph{exchangeable}, and those in $[1,N] \setminus \ex$ \emph{frozen}. Matrices indexed by $[1,N] \times \ex$ will be called $N \times \ex$ matrices. Fix a purely transcendental field extension $\FF \supset \KK$ of transcendence degree $N$. We define cluster algebras over $\KK$ directly, rather than first defining them over $\Zset$ and then extending scalars. Moreover, we follow \cite{FZ1} rather than \cite{BFZ} in that frozen variables in cluster algebras are not automatically inverted. In our notation, we combine the mutation matrices for cluster and frozen variables as in \cite[p.~515]{FZ1} and \cite[p.~6]{BFZ}.

\bde{seeds}
A \emph{seed} (of geometric type) in $\FF$ is a pair $(\xbtil,\Btil)$ where
\begin{enumerate}
\item[(i)] $\xbtil = (x_1,\dots,x_N) \in \FF^N$ such that $\{x_1,\dots,x_N\}$ is algebraically independent over $\KK$ and generates $\FF$ as a field extension of $\KK$.
\item[(ii)] $\Btil = (b_{ij})$ is an $N \times \ex$ integer matrix of full rank $n$.
\item[(iii)] The $\ex \times \ex$ submatrix $B$ of $\Btil$ is \emph{skew-symmetrizable}: there exist $d_i \in \Zset_{>0}$ such that $d_i b_{ij} = - d_j b_{ji}$ for all $i,j \in \ex$.
\end{enumerate}
Note that (iii) implies $b_{kk} = 0$ for all $k \in \ex$.

The set $\xb := \{ x_j \mid j \in \ex\} \subseteq \xbtil$ is a \emph{cluster}, and its elements are \emph{cluster variables}. Elements of $\xbtil \setminus \xb$ are \emph{frozen variables}. The submatrix $B$ of $\Btil$ is the \emph{principal part} of $\Btil$, and is called the \emph{exchange matrix} of the seed $(\xbtil,\Btil)$.\ede

\bde{mukxtil}
Let $(\xbtil,\Btil)$ be a seed, with $\xbtil = (x_1,\dots,x_N)$. For $k \in \ex$, define
\begin{equation}
\label{xk'}
x'_k := \dfrac{\prod_{b_{ik}>0} x_i^{b_{ik}} + \prod_{b_{ik}<0} x_i^{-b_{ik}}} {x_k} \in \FF
\end{equation}
and $\mu_k(\xbtil) := (x_1,\dots, x_{k-1}, x'_k, x_{k+1}, \dots,x_N)$. (Then $\mu_k(\xbtil)$ is another sequence of algebraically independent elements generating $\FF$ over $\KK$.)
\ede

\bde{mutation}
Let $(\xbtil,\Btil)$ be a seed, $\Btil = (b_{ij})$, and $k \in \ex$. Define $\Btil' = \mu_k(\Btil) = (b'_{ij})$ by
\begin{equation}
\label{bij'}
b'_{ij} := \begin{cases}  -b_{ij}  &\text{if} \;\; i=k \;\; \text{or} \;\; j=k  \\
b_{ij} + \dfrac{ |b_{ik}| b_{kj} + b_{ik} |b_{kj}| }{2}  &\text{otherwise.} \end{cases}
\end{equation}
The matrix $\mu_k(\Btil)$ is the \emph{mutation of $\Btil$ in direction $k$}.

If the principal part of $\Btil$ is $B$, then the principal part of $\mu_k(\Btil)$ is $\mu_k(B)$. Moreover, $\mu_k(\Btil)$ preserves the conditions from \deref{seeds}: it has full rank $n$ \cite[Lemma 3.2]{BFZ}, and $\mu_k(B)$ is skew-symmetrizable \cite[Proposition 4.5]{FZ1}.

The \emph{mutation of $(\xbtil,\Btil)$ in direction $k$} is the seed $\mu_k(\xbtil,\Btil) := (\mu_k(\xbtil), \mu_k(\Btil))$. Mutation is involutive: $\mu_k( \mu_k(\xbtil,\Btil)) = (\xbtil,\Btil)$.

In general, seeds $(\xbtil,\Btil)$ and $(\xbtil',\Btil')$ are \emph{mutation-equivalent}, written $(\xbtil,\Btil) \sim (\xbtil',\Btil')$, if $(\xbtil',\Btil')$ can be obtained from $(\xbtil,\Btil)$ by a sequence of seed mutations. Mutation-equivalent seeds share the same set of frozen variables.
\ede

\bde{uppercluster}
Upper bounds and upper cluster algebras associated to the seed $(\xbtil,\Btil)$ are intersections of mixed polynomial-Laurent polynomial algebras generated by the variables from some cluster along with inverses of some of them. For any set $\inv \subseteq [1,N] \setminus \ex$ of frozen indices, let us define
\begin{equation}
\label{TAexinv}
\TT\AA (\xbtil,\Btil, \inv) := \KK [x_i \mid i\in [1,N] \setminus (\ex \sqcup \inv)] [x_j^{\pm1} \mid j \in \ex \sqcup \inv].
\end{equation}
The \emph{upper bound} for $(\xbtil,\Btil, \inv)$ is
\begin{equation}
\label{UB}
\UU (\xbtil,\Btil, \inv)_\KK := \TT\AA (\xbtil,\Btil, \inv) \cap \bigcap_{k \in \ex} \TT\AA (\mu_k(\xbtil),\mu_k(\Btil), \inv),
\end{equation}
and the corresponding \emph{upper cluster algebra} is
\begin{equation}
\label{uppercluster}
\overline{\AA} (\xbtil,\Btil, \inv)_\KK := \bigcap_{ (\xbtil', \Btil') \, \sim \, (\xbtil, \Btil)} \TT\AA (\xbtil', \Btil', \inv).
\end{equation}
\ede

Because we are restricting attention to seeds of geometric type, upper bounds coincide with upper cluster algebras:

\bth{U=Abar}
$\UU (\xbtil,\Btil, \inv)_\KK = \overline{\AA} (\xbtil,\Btil, \inv)_\KK$ for any seed $(\xbtil, \Btil)$ in $\FF$ and any $\inv \subseteq [1,N] \setminus \ex$.
\eth

\begin{proof}
The case $\inv = \varnothing$ follows explicitly from \cite[Corollary 1.9]{BFZ}, as a consequence of the relation
\begin{equation}
\label{UxB=Ux'B'}
\UU (\xbtil,\Btil, \varnothing)_\KK = \UU (\xbtil', \Btil', \varnothing)_\KK , \; \; \forall (\xbtil', \Btil') \sim (\xbtil, \Btil).
\end{equation}

For any seed $(\xbtil'', \Btil'')$, we have 
$$
\TT\AA (\xbtil'', \Btil'', \inv) = \TT\AA (\xbtil'', \Btil'', \varnothing) [x_j^{-1} \mid j \in \inv].
$$
Consequently, any element $u \in \UU (\xbtil,\Btil, \inv)_\KK$ can be expressed in the form $u = u_0 z^{-1}$ where $u_0 \in \UU (\xbtil, \Btil, \varnothing)$ and $z$ is a monomial in the $x_j$ for $j \in \inv$. In view of \eqref{UxB=Ux'B'}, $u$ must belong to $\TT\AA (\xbtil', \Btil', \inv)$ for all seeds $(\xbtil', \Btil')$ mutation-equivalent to $(\xbtil, \Btil)$. It follows that $\UU (\xbtil,\Btil, \inv)_\KK \subseteq \overline{\AA} (\xbtil,\Btil, \inv)_\KK$. Since the reverse inclusion holds a priori, the theorem is proved.
\end{proof}

Following this result, we refer to $\UU (\xbtil,\Btil, \inv)_\KK$ as the upper cluster algebra corresponding to $(\xbtil,\Btil, \inv)$.

\bde{clusteralg}
The \emph{cluster algebra} associated to a seed $(\xbtil,\Btil)$ in $\FF$ is the $\KK$-subalgebra $\AA(\xbtil,\Btil)_\KK$ of $\FF$ generated by the cluster and frozen variables from all the seeds mutation-equivalent to $(\xbtil,\Btil)$. More generally, for any set $\inv \subseteq [1,N] \setminus \ex$, we define 
$$\AA(\xbtil,\Btil,\inv)_\KK := \AA(\xbtil,\Btil)_\KK[x_l^{-1} \mid l \in \inv].
$$
These algebras are also called \emph{cluster algebras}.
\ede

The \emph{Laurent Phenomenon} for cluster algebras may be expressed as follows; see \cite[Corollary 1.12]{BFZ} for the version in terms of upper cluster algebras.

\bth{phenom}
{\rm \cite[Theorem 3.1]{FZ1}}
Every cluster algebra is contained in the corresponding upper cluster algebra: $\AA(\xbtil,\Btil, \inv)_\KK \subseteq \UU(\xbtil,\Btil, \inv)_\KK$ for all seeds $(\xbtil,\Btil)$ in $\FF$ and all $\inv \subseteq [1,N] \setminus \ex$. Equivalently,
$$
\AA(\xbtil,\Btil, \inv)_\KK \subseteq \TT\AA (\xbtil, \Btil, \inv)
$$
for all $(\xbtil,\Btil)$ and $\inv$.
\eth

If $G$ is a group acting on $\FF$ by automorphisms, we denote by $\FF^G$ the fixed field (or Galois subfield) of this action. For a $G$-eigenvector $u$, we denote its $G$-eigenvalue by $\chi_u$, or by $\chi[u]$ in case $u$ is a lengthy expression. The character lattice $X(G)$ is viewed as an additive group.

The following equivariance of mutations of seeds is probably well known, but we could not locate a reference.

\ble{equivar-mut}
Let $G$ be a group acting on $\FF$ by $\KK$-algebra automorphisms, and let $(\xbtil, \Btil)$ be a seed in $\FF$. Assume that all entries of $\xbtil$ are $G$-eigenvectors, and that $\xbtil^{b^j} \in \FF^G$ for all columns $b^j$ of $\Btil$. Then all seeds $(\xbtil', \Btil')$ mutation-equivalent to $(\xbtil, \Btil)$ have the same properties.
\ele

\begin{proof} It suffices to prove the lemma when $(\xbtil', \Btil')= \mu_k(\xbtil, \Btil)$ for some $k \in \ex$. Denote the entries of $\xbtil$, $\xbtil'$, $\Btil$, and $\Btil'$ by $x_j$, $x'_j$, $b_{lj}$, and $b'_{lj}$, respectively, and denote the columns of $\Btil'$ by $(b')^j$.

We may write $\xbtil^{b^k}$ in the form
$$
\xbtil^{b^k} = \xbtil^{[b^k]_+} \, \xbtil^{[b^k]_-}
$$
(recall \eqref{f+-}).
As $\xbtil^{[b^k]_\pm}$ are $G$-eigenvectors and $\chi[\xbtil^{b^k}] = 0$ by hypothesis, we see that
$$
\chi[ \xbtil^{[b^k]_-} ] = - \chi[ \xbtil^{[b^k]_+} ].
$$
Rewriting \eqref{xk'} in the form
$$
x'_k = x_k^{-1} \bigl( \xbtil^{[b^k]_+} + \xbtil^{- [b^k]_-} \bigr),
$$
we find that $x'_k$ is a $G$-eigenvector with
\begin{equation}
\label{chixk'}
\chi[x'_k] = - \chi[x_k] + \chi[ \xbtil^{[b^k]_+} ] = - \chi[x_k] - \chi[ \xbtil^{[b^k]_-} ].
\end{equation}
Of course, $x'_j = x_j$ is a $G$-eigenvector for all $j \ne k$, by assumption.

Turning to the second condition, we have $(\xbtil')^{(b')^k} = \xbtil^{-b^k}$ because $b'_{kk} = 0$ and $b'_{ik} = -b_{ik}$ for $i \ne k$, whence $(\xbtil')^{(b')^k} \in \FF^G$. Now let $j \ne k$. If $b_{kj} = 0$, then $(b')^j = b^j$ and $(\xbtil')^{(b')^j} = \xbtil^{b^j}$, yielding $(\xbtil')^{(b')^j} \in \FF^G$. Finally, if $b_{kj} \ne 0$, set $\ep = \sign(b_{kj})$, and observe from \eqref{bij'} that
$$
(b')^j = b^j + \ep b_{kj} [b^k]_\ep - 2b_{kj} e_k .
$$
Since the $k$-th entry of $b^j + \ep b_{kj} [b^k]_\ep - b_{kj} e_k$ is zero,
$$
(\xbtil')^{(b')^j} = \xbtil^{ b^j + \ep b_{kj} [b^k]_\ep - b_{kj} e_k } \, (x'_k)^{-b_{kj}}, 
$$
and using \eqref{chixk'} we obtain
$$
\chi[ (\xbtil')^{(b')^j} ] = \chi[ \xbtil^{b^j} ] + \ep b_{kj} \chi[ \xbtil^{\b^k]_\ep} ] - b_{kj} \chi[x_k] - b_{kj} \bigl( - \chi[x_k] + \ep \chi[ \xbtil^{\b^k]_\ep} ] \bigr) = 0, 
$$
completing the proof.
\end{proof}

\subsection{Compatible pairs}
\label{compatible}

Fix $n$, $N$, $\ex$, and $\FF$ as in \S\ref{clusternota}. We use a slight generalization of the notion of compatible pairs from \cite[Definition 3.1]{BZ}.

\bde{compat}
Let $\Btil = (b_{kj}) \in M_{N\times \ex}(\Zset)$, and let $\rb = (r_{ij}) \in M_N(\KK)$ be a skew-symmetric scalar matrix. We say that the pair $(\rb, \Btil)$ is \emph{compatible} if the following two conditions are satisfied:
\begin{align}
&
\sum_{i=1}^N b_{ik} r_{ij} = 0, \; \; \forall k \in \ex, \; j \in [1,N], \; k \neq j
\quad \mbox{and}
\label{comp1}
\\
&\sum_{i=1}^N b_{ik} r_{ik} \ne 0, \; \; \forall k \in \ex.
\label{comp2}
\end{align}
Note that the scalars appearing in \eqref{comp1} and \eqref{comp2} are just the entries of the matrix $\Btil^T \rb$. Moreover, due to the skew-symmetry of $\rb$ and the bicharacter $\Om_\rb$ associated to $\rb$ as in \eqref{Om}, we have
\begin{equation}  \label{compscalars}
\sum_{i=1}^N b_{ik} r_{ij} = \Om_\rb(b^k, e_j) \;\; \forall k \in \ex, \; j \in [1,N],
\end{equation}
where $b^k$ denotes the $k$-th column of $\Btil$.
\ede

\bpr{fullrank} If the pair $(\rb, \wt{B})$ is compatible,
then $\Btil$ has full rank, and the nonzero scalars $\beta_k := (\Btil^T \rb)_{kk}$ satisfy
\begin{equation}
\label{scalarsymm}
\beta_k b_{kj} = - \beta_j b_{jk}, \; \; \forall k,j \in \ex.
\end{equation}
\epr

\begin{proof} This is proved just as \cite[Proposition 3.3]{BZ}.
\end{proof}

Unlike the case in \cite[Proposition 3.3]{BZ}, compatibility of $(\rb, \wt{B})$ does not in general imply that the principal part of $\Btil$ is skew-symmetrizable. The next lemma describes an instance 
when this condition appears naturally.

\ble{alternative} Assume that $\Btil = (b_{kj}) \in M_{N\times \ex}(\Zset)$ and $\rb \in M_N(\KK)$ form a compatible pair, and set $\beta_k := (\Btil^T \rb)_{kk}$ for $k \in \ex$. If there exist positive integers 
$d_k$, for $k \in \ex$, such that   
\begin{equation}
\label{dbeta}
d_j \beta_k = d_k \beta_j, \; \;
\forall j,k \in \ex,
\end{equation}
then the principal part of $\wt{B}$ is skew-symmetrizable via these $d_k$, that is, $d_k b_{kj} = - d_j b_{jk}$ for all $k,j \in \ex$.
\ele

\begin{proof} By \prref{fullrank},
$\beta_k b_{kj} = - \beta_j b_{jk}$ for all $j, k \in \ex$.
Combining this with \eqref{dbeta} leads to 
$$
d_k \beta_j b_{kj} = d_j \beta_k b_{kj} = - d_j \beta_j b_{jk} .
$$ 
Thus $d_k b_{kj} = - d_j b_{jk}$ because of \eqref{comp2}.
\end{proof}

\bre{powers} The condition \eqref{dbeta} is satisfied if and only if there exists $q \in \kx$ such that each $\beta_k = m_k q$ for some 
$m_k \in \Zset_{> 0}$. Then one can set $d_k := m_k$. 
\ere

We define mutations of compatible pairs as in \cite[Definition 3.5]{BZ}.

\bde{pairmuta}
Let $(\rb, \Btil)$ be a compatible pair, with $\Btil = (b_{lj}) \in M_{N\times \ex}(\Zset)$, and let $k \in \ex$.
By \cite[Eq. (3.2)]{BFZ}, the mutated matrix $\mu_k(\wt{B})$ 
can be expressed as 
$$
\mu_k ( \wt{B}) = E_\ep \wt{B} F_\ep 
$$ 
for both choices of sign $\ep = \pm$, where $E_\ep = E_\ep^{\Btil, k}$ and $F_\ep = F_\ep^{\Btil, k}$ are 
the $N \times N$ and $\ex \times \ex$ matrices with entries given by
\begin{equation}
\label{Eep}
\begin{aligned}
(E_\ep)_{ij} &:= 
\begin{cases}
\de_{ij}, & \mbox{if} \; j \neq k \\
-1, & \mbox{if} \; i=j=k \\
\max(0, - \ep b_{ik}), & \mbox{if} \; i \neq j = k
\end{cases}  \\
(F_\ep)_{ij} &:= 
\begin{cases}
\de_{ij}, & \mbox{if} \; i \neq k \\
-1, & \mbox{if} \; i=j=k \\
\max(0, \ep b_{kj}), & \mbox{if} \; i =k \neq j.
\end{cases}
\end{aligned}
\end{equation}
We define the mutation in direction $k$ 
of the matrix $\rb$ by
\begin{equation}
\label{r-mut-k}
\mu_k(\rb) := E_\ep^T \rb E_\ep.
\end{equation}
\ede

\bpr{pair-mut} Let $(\rb, \wt{B})$ be a compatible pair and $k \in \ex$.

{\rm(a)} The matrix $\mu_k(\rb)$, defined in \eqref{r-mut-k}, 
does not depend on the choice of sign $\ep= \pm$. It is skew-symmetric, and the pair $(\mu_k(\rb),\mu_k(\wt{B}))$ is compatible.

{\rm(b)} Assume also that the principal part of $\wt{B}$ is skew-symmetrizable. Then the principal part of $\mu_k(\wt{B})$ is skew-symmetrizable, and 
\begin{equation}
\label{mutaBTr}
\mu_k(\Btil)^T \mu_k(\rb) = \Btil^T \rb .
\end{equation}
\epr

We define the \emph{mutation in direction $k \in \ex$} of the compatible 
pair $(\rb, \wt{B})$ to be the compatible pair $(\mu_k(\rb),\mu_k(\wt{B}))$.
\medskip

\begin{proof}[Proof of Proposition {\rm\ref{ppair-mut}}]
Part (a) is proved as \cite[Proposition 3.4]{BZ}.

(b) The principal part of $\mu_k(\wt{B})$, namely $\mu_k(B)$, is skew-symmetrizable (for the same choice of positive integers $d_j$, $j \in \ex$ as for $B$), by the observations in \cite[Proposition 4.5]{FZ}.

We have
$$
\mu_k(\wt{B})^T \mu_k(\rb)) = F_\ep^T \wt{B}^T E_\ep^T E_\ep^T \rb E_\ep  = 
 F_\ep^T \wt{B}^T \rb E_\ep.
$$
The second statement in \prref{fullrank} and the fact that for all $i, j \in \ex$,
$b_{ij}$ and $-b_{ji}$ have the same signs (which follows from the skew-symmetrizability assumption) imply 
$$
\wt{B}^T \rb E_\ep = F_\ep^T\wt{B}^T \rb. 
$$
Therefore $
\mu_k(\wt{B})^T \mu_k(\rb) =  F_\ep^T \wt{B}^T \rb E_\ep = \wt{B}^T \rb.$ 
\end{proof}

\subsection{Poisson cluster algebras}
\label{Pcluster}

Continue with $n$, $N$, $\ex$, and $\FF$ as in \S\ref{clusternota}. We first slightly extend two concepts from \cite[Introduction and \S1.3]{GSV}, \cite[\S4.1.1]{GSVbook}.

\bde{P-compat}
Suppose we have a Poisson $\KK$-algebra structure on $\FF$. An $N$-tuple $\xbtil = (x_1, \dots, x_N) \in \FF^N$ is \emph{log-canonical} (with respect to the given Poisson structure) if $\{x_l,x_j\} \in \KK x_l x_j$ for all $l,j \in [1,N]$. We say that a cluster algebra $\AA(\xbtil, \Btil)_\KK \subset \FF$ is \emph{Poisson-compatible} in case $\xbtil'$ is log-canonical for all seeds $(\xbtil', \Btil')$ mutation-equivalent to $(\xbtil, \Btil)$. The same terminology is used with the localized cluster algebras $\AA(\xbtil, \Btil, \inv)_\KK$. In the Poisson-compatible case, the upper cluster algebra $\UU(\xbtil, \Btil,\inv)_\KK$ is a Poisson subalgebra of $\FF$.
\ede

\bpr{brack-mut}
Fix a $\KK$-algebra Poisson bracket on $\FF$. Let $(\xbtil, \Btil)$ be a seed in $\FF$ and $k \in \ex$, and write $\xbtil = (x_1,\dots,x_N)$ and $\mu_k(\xbtil) = (x'_1,\dots,x'_N)$. Suppose that $\rb \in M_N(\KK)$ is a matrix such that $(\rb, \Btil)$ is a compatible pair. If
\begin{equation}
\label{brackxx}
\{x_l, x_j\} = \Om_\rb(e_l,e_j) x_l x_j , \; \; \forall l,j \in [1,N],
\end{equation}
then
\begin{equation}
\label{brackx'x'}
\{x'_l, x'_j\} = \Om_{\mu_k(\rb)}(e_l,e_j) x'_l x'_j , \; \; \forall l,j \in [1,N].
\end{equation}
\epr

\begin{proof} Let $\rb' := \mu_k(\rb) = E_\ep^T \rb E_\ep$, where $E_\ep = E_\ep^{\Btil, k}$. Then
\begin{equation}
\label{Ommur}
\Om_{\rb'}(f,g) = f^T E_\ep^T \rb E_\ep g = \Om_{\rb}( E_\ep f, E_\ep g), \; \; \forall f,g \in \Zset^N, \; \ep = \pm.
\end{equation}
Write $b^k$ for the $k$-th column of $\Btil$. From the definitions \eqref{Eep} and \eqref{xk'}, we see that
\begin{equation}
\label{Eepel}
\begin{aligned}
E_\ep e_l &= e_l , \; \; \forall \ep = \pm, \; l \in [1,N], \; l \ne k  \\
E_+ e_k &= - e_k - [b^k]_-  \\
E_- e_k &=  - e_k + [b^k]_+ ,
\end{aligned}
\end{equation}
and so
\begin{equation}
\label{xk'alt}
x'_k = \xbtil^{E_- e_k} + \xbtil^{E_+ e_k} .
\end{equation}

If $l,j \in [1,N]$ with $l,j \ne k$, then \eqref{Ommur} and \eqref{Eepel} imply that $\Om_{\rb'}(e_l, e_j) = \Om_\rb(e_l, e_j)$, whence
$$
\{ x'_l, x'_j \} = \{ x_l, x_j \} = \Om_\rb(e_l,e_j) x_l x_j = \Om_{\rb'}(e_l,e_j) x'_l x'_j  .
$$
When $l \ne k$, we have $\Om_{\rb'}( e_l, e_k) = \Om_\rb( e_l, E_\ep e_k)$ for $\ep = \pm$, and so, taking account of \eqref{xk'alt}, we have
$$
\{ x'_l, x'_k \} = \Om_\rb( e_l, E_- e_k) x_l \xbtil^{E_- e_k} + \Om_\rb( e_l, E_+ e_k) x_l \xbtil^{E_+ e_k} = \Om_{\rb'}( e_l, e_k) x'_l x'_k .
$$
The remaining cases of \eqref{brackx'x'} follow via the skew-symmetry of $\{-,-\}$ and $\Om_{\rb'}$.
\end{proof}

\bth{P-cluster-main}
Suppose that $\FF$ is a Poisson $\KK$-algebra. Let $(\xbtil, \Btil)$ be a seed in $\FF$, and $\inv \subseteq [1,N] \setminus \ex$. Assume there exists a skew-symmetric matrix $\rb \in M_N(\KK)$ such that $(\rb, \Btil)$ is a compatible pair and the entries $x_i$ of $\xbtil$ satisfy
$$
\{x_l, x_j\} = \Om_\rb(e_l,e_j) x_l x_j , \; \; \forall l,j \in [1,N].
$$
Then the cluster algebra $\AA(\xbtil, \Btil, \inv)_\KK$ is Poisson-compatible. 
\eth

\begin{proof} It suffices to prove the theorem in the case $\inv = \varnothing$.

Let $\wt{T}$ denote the collection of all triples $(\xbtil', \Btil', \rb')$ such that $(\xbtil', \Btil')$ is a seed in $\FF$ and $(\rb', \Btil')$ is a compatible pair. For any such triple and any $k \in \ex$, we have a mutated seed $(\mu_k(\xbtil'), \mu_k(\Btil'))$ and a mutated compatible pair $(\mu_k(\rb'), \mu_k(\Btil'))$ (by \prref{pair-mut}), so that $\mu_k(\xbtil', \Btil', \rb') := (\mu_k(\xbtil'), \mu_k(\Btil'), \mu_k(\rb'))$ is another triple in $\wt{T}$. This defines mutation operations $\mu_k$ on $\wt{T}$.

Let $\wt{T}_0$ be the collection of all triples in $\wt{T}$ mutation-equivalent to $(\xbtil, \Btil, \rb)$. It follows from \prref{brack-mut} that for each $(\xbtil', \Btil', \rb') \in \wt{T}_0$, the entries $x'_i$ of $\xbtil'$ satisfy
$$
\{x'_l, x'_j\} = \Om_{\rb'}(e_l,e_j) x'_l x'_j , \; \; \forall l,j \in [1,N].
$$
In particular, $\xbtil'$ is log-canonical.

For any seed $(\xbtil', \Btil') = \mu_{k_n} \cdots \mu_{k_1}(\xbtil, \Btil)$ mutation-equivalent to $(\xbtil, \Btil)$, we have a matrix $\rb' := \mu_{k_n} \cdots \mu_{k_1}(\rb)$ such that $(\xbtil', \Btil', \rb') \in \wt{T}_0$, and thus $\xbtil'$ is log-canonical. Therefore $\AA(\xbtil, \Btil)_\KK$ is Poisson-compatible.
\end{proof}

\sectionnew{Poisson-primes in Poisson-Ore extensions}
\label{PpPO}

\subsection{Equivariant Poisson unique factorization domains}
\label{p-ufd}
Let $R$ be a Poisson domain (i.e., a Poisson algebra which is also an integral domain). A \emph{Poisson-prime element} of $R$ is any 
Poisson-normal prime element $p \in R$, that is, $p$ is a nonzero Poisson-normal element and $R/(p)$ is a domain (where $(p)$ denotes the principle ideal $Rp$). 
It follows from \prref{Pvleaves} that for every affine Poisson domain $R$, the following 
are equivalent for $p \in R$:
\begin{enumerate}
\item $p$ is a Poisson-prime.
\item The ideal $(p)$ is a prime ideal and a Poisson ideal.
\item $p$ is a prime element of $R$ and the zero locus $V(p)$ is a union of symplectic leaves of the spectrum of $R$.
\item $p$ is a prime element of $R$ and $V(p)$ is a union of Zariski closures of symplectic leaves of $\Spec R$.
\end{enumerate}
In the special case when the spectrum of $R$ is smooth, the equivalence of (1) and (3) was proved in \cite[Remark 2.4(iii)]{NTY}. 
In case $\KK$ is uncountable and algebraically closed, \leref{P.Jacobson.cond} and \prref{Pvcores} show that $p$ is Poisson-prime if and only if 
\begin{enumerate}
\item[(5)] $p$ is a prime element of $R$ and $V(p)$ is a union of symplectic cores of $\Spec R$.
\end{enumerate}

We shall say that $R$ is a \emph{Poisson unique factorization domain}, abbreviated \emph{Poisson-UFD} or \emph{P-UFD}, if each nonzero Poisson-prime ideal of $R$ contains a Poisson-prime element. The hyphen is important here -- a Poisson algebra which is also a UFD need not be a P-UFD. For example, equip the polynomial algebra $R = \Cset[x,y,z]$ with the Poisson bracket such that
$$
\{x, y\} = 0, \qquad \{z, y\} = x+y^2, \qquad \{z, x\} = 2y.
$$
As shown in \cite[Example 5.12]{Jor}, $\langle x,y \rangle$ is a Poisson prime ideal of $R$ which does not properly contain any nonzero Poisson prime ideal. Hence, $\langle x,y \rangle$ does not contain any Poisson-prime element.

We shall also need an equivariant version of this concept, as in \cite{GY1, GYbig}. Suppose that $R$ is a Poisson domain equipped with a Poisson action of a group $\HH$. On replacing ``Poisson ideal" by ``$\HH$-stable Poisson ideal" in the definition of a Poisson-prime ideal, we obtain the definition of an \emph{$\HH$-Poisson-prime} ideal of $R$. We shall say that $R$ is an \emph{$\HH$-Poisson unique factorization domain}, abbreviated \emph{$\HH$-P-UFD}, if each nonzero $\HH$-Poisson-prime ideal of $R$ contains a Poisson-prime $\HH$-eigenvector. The arguments of  \cite[Proposition 2.2]{GY1} and
\cite[Proposition 6.18 (ii)]{Y-sqg} are easily adapted to give the following equivariant Poisson version of results of 
Chatters and Jordan \cite[Proposition 2.1]{Cha},
\cite[p. 24]{ChJo}. 

\bpr{factorHPUFD}
Let $R$ be a noetherian $\HH$-Poisson-UFD. Every Poisson-normal $\HH$-eigen\-vec\-tor in $R$ is either a unit or a product of Poisson-prime $\HH$-eigenvectors. The factors are unique 
up to reordering and taking associates.
\epr

\bre{HPUFDtoPUFD}
A Poisson version of \cite[Proposition 1.6 and Theorem 3.6]{LLR}, as adapted in \cite[Theorem 2.4]{GY1}, shows that if $R$ is a noetherian $\HH$-Poisson-UFD where $\HH$ is a $\KK$-torus acting rationally on $R$ by Poisson automorphisms, then $R$ is a Poisson-UFD. We will not need this result, however.
\ere

The next lemma clarifies the picture of $\HH$-Poisson-prime ideals in the case of a rational torus action.

\ble{torusPprime}
Let $R$ be a noetherian Poisson algebra, equipped with a rational Poisson action of a torus $\HH$. Then any prime ideal of $R$ minimal over an $\HH$-stable Poisson ideal is itself an $\HH$-stable Poisson ideal, and the $\HH$-Poisson-prime ideals of $R$ are exactly the $\HH$-stable, Poisson, prime ideals of $R$.
\ele

\begin{proof} The first statement follows from \cite[Lemma 1.1(d)]{GpDixMo} and \cite[Proposition II.2.9]{BrGo}.  Obviously $\HH$-stable, Poisson, prime ideals are $\HH$-Poisson-prime. Conversely, let $P$ be an $\HH$-Poisson-prime ideal of $R$. Since $R$ is noetherian, there are prime ideals $Q_1,\dots,Q_m$ minimal over $P$ such that $Q_1 Q_2 \cdots Q_m \subseteq P$. Each $Q_i$ is an $\HH$-stable Poisson ideal. The  $\HH$-Poisson-primeness of $P$ thus implies that some $Q_j \subseteq P$, whence $Q_j = P$. Therefore $P$ is a prime ideal.
\end{proof}

\leref{torusPprime} allows a shortening of terminology: In the setting of the lemma,
\begin{equation}
\label{pHprime}
\{ \HH\text{-Poisson-prime ideals of} \; R\} = \{ \text{Poisson} \; \HH\text{-prime ideals of} \; R\}.
\end{equation} 
This is immediate from the lemma together with the fact that the $\HH$-prime ideals of $R$ coincide with the $\HH$-stable prime ideals \cite[Proposition II.2.9]{BrGo}.

\subsection{Poisson-Cauchon extensions}
\label{p-Cauchon}
In this section, we study a Poisson analog of the \emph{Cauchon extensions} defined in \cite[Definition 2.5]{LLR}, and develop a version of \cite[Proposition 2.9]{LLR} which allows us to prove that suitable Poisson-Cauchon extensions are $\HH$-P-UFDs. 

\bde{p-Cauchon-def}
A \emph{Poisson-Cauchon extension} of a Poisson algebra $B$ is a Poisson-Ore extension $R = B[x; \sig, \de]_p$ which is equipped with a rational  Poisson action of a torus $\HH$ such that
\begin{enumerate}
\item[(i)] The subalgebra $B$ is $\HH$-stable, and $x$ is an $\HH$-eigenvector.
\item[(ii)] $\de$ is locally nilpotent.
\item[(iii)] There is some $h_\ci \in \h = \Lie\HH$ such that $\sig = (h_\ci\cdot)|_B$ and the $h_\ci$-eigenvalue of $x$, to be denoted $\la_\ci$, is nonzero.
\end{enumerate} 
Condition (iii) relies on $x$ being an $\h$-eigenvector, which follows from (i) and \eqref{hHeigenvectors}.
(This definition is not a precise parallel to \cite[Definition 2.5]{LLR}, since we have not assumed that $B$ is a noetherian domain, while we have restricted $\HH$ to being a torus.)

Whenever we refer to a Poisson-Cauchon extension $R = B[x; \sig, \de]_p$, we take the associated torus and its Lie algebra to be denoted $\HH$ and $\h$ unless otherwise specified.

Since $\HH$-stable subspaces of $B$ are also $\h$-stable (\S \ref{dHaction}), it follows from (iii) that
\begin{equation}
\label{Hsig.stable}
\text{Every} \; \HH\text{-stable subspace of} \; B \; \text{is} \; \sig\text{-stable}.
\end{equation}

Our assumptions also imply that
\begin{align}
\sig \ci \de &= \de \ci (\sig + \la_\ci)  \label{aldel}  \\
(h\cdot)|_B \ci \de &= \chi_x(h) \de \ci (h\cdot)|_B , \; \; \forall h \in \HH.  \label{hdotdel}
\end{align}
Namely, for $b \in B$ and $h \in \HH$ we compute that
\begin{align*}
\sig \de (b) &= h_\ci \cdot (\{x, b\} - \sig(b) x) = \{\la_\ci x, b\} + \{x, \sig(b)\} - \sig^2(b) x - \sig(b) \la_\ci x  \\
 &= \la_\ci \sig(b) x + \la_\ci \de(b) + \sig^2(b) x + \de \sig(b)- \sig^2(b) x - \sig(b) \la_\ci x = \de( \sig(b) + \la_\ci b)  \\
h\cdot \de(b) &= h\cdot (\{x, b\} - \sig(b) x) = \{ \chi_x(h) x, h\cdot b\} - h\cdot \sig(b) \chi_x(h) x  \\
 &= \chi_x(h) (\{x, h\cdot b\} - \sig(h\cdot b) x) = \chi_x(h) \de(h\cdot b),
\end{align*}
where we have used \eqref{hHcomm} to commute $(h\cdot)|_B$ and $\sig$. Further, this commutativity, together with \eqref{hdotdel}, implies that for any nonzero homogeneous element $b \in B$, the elements $\sig(b)$ and $\de(b)$ are homogeneous, with
\begin{equation}
\label{chisigdel}
\chi_{\sig(b)} = \chi_b \; \; \text{if} \; \sig(b) \ne 0, \; \; \text{and} \; \; \chi_{\de(b)} = \chi_b + \chi_x \; \; \text{if} \; \de(b) \ne 0.
\end{equation}
\ede

There is a Poisson version of Cauchon's derivation-deleting map \cite[Section 2]{Cau}, as follows.

\bpr{pCauchonmap}
Let $R = B[x; \sig, \de]_p$ be a Poisson-Cauchon extension.

{\rm (a)} There is a Poisson algebra isomorphism $\theta : B[y^{\pm1}; \sig]_p \longrightarrow \Rhat := B[x^{\pm 1}; \sig, \de]_p$ such that $\theta(y) = x$ and
\begin{equation}
\label{pCmap}
\theta(b) = \sum_{n=0}^\infty \, \frac1{n!} \left( \frac{-1}{\la_\ci} \right)^n \de^n(b) x^{-n}, \; \; \forall b \in B.
\end{equation}

{\rm (b)} $\{ x, \theta(b)\} = \theta\sig(b) x$ for all $b \in B$.

{\rm (c)} The given $\HH$-action on $B$ extends uniquely to a rational action of $\HH$ on $B[y^{\pm1}; \sig]_p$ by Poisson automorphisms such that $\chi_y = \chi_x$. With respect to this action, the isomorphism $\theta$ is $\HH$-equivariant.
\epr

\begin{proof} Parts (a) and (b) are proved in \cite[Lemmas 3.4, 3.6, Theorem 3.7]{GLaun}. 

(c) It is clear that the $\HH$-action on $B$ extends uniquely to a rational action of $\HH$ on $B[y^{\pm1}; \sig]_p$ by $\KK$-algebra automorphisms such that $\chi_y = \chi_x$. 

Any $h \in \HH$ already acts on $B$ by Poisson automorphisms. Moreover, $(h\cdot)|_B$ commutes with $(h_\ci \cdot)|_B = \sig$ by \eqref{hHcomm}, whence
$$
h \cdot \{y,b\} = h \cdot (\sig(b) y) = \sig(h\cdot b) \chi_y(h) y = \chi_y(h) \{y, h\cdot b\} = \{h\cdot y, h\cdot b\}, \; \; \forall b \in B.
$$
It follows that $h$ acts on $B[y^{\pm1}; \sig]_p$ by a Poisson automorphism.

The $\HH$-action on $B[y^{\pm1}; \sig]_p$ has been chosen so that $\theta(h\cdot y) = h\cdot x = h\cdot \theta(y)$ for any $h \in \HH$.  From \eqref{hdotdel}, we obtain $(h\cdot)|_B \de^n = \chi_x(h)^n \de^n (h\cdot)|_B$ for all $h \in \HH$ and $n \in \Znn$, and consequently $h\cdot (\de^n(b) x^{-n}) = \de^n(h\cdot b) x^{-n}$ for all $b \in B$. This implies that $h\cdot \theta(b) = \theta(h\cdot b)$ for all $b \in B$, and we conclude that $(h\cdot)\theta = \theta(h\cdot)$.
\end{proof}

\bco{more.pCauchon}
Let $R = B[x; \sig, \de]_p$, $\theta$, and $\Rhat$ be as in Proposition {\rm\ref{ppCauchonmap}}. 

{\rm (a)} $\theta$ restricts to a Poisson algebra isomorphism of $B[y; \sig]_p$ onto the subalgebra $\theta(B)[x; \al]_p$ of $\Rhat$, where $\al$ is the Poisson derivation on $\theta(B)$ such that $\al \theta(b) = \theta \sig(b)$ for $b \in B$. 

{\rm (b)} $\al = (h_\ci \cdot)|_{\theta(B)}$ and $\theta(B)[x; \al]_p$ is a Poisson-Cauchon extension of $\theta(B)$.

{\rm (c)} $\Rhat = \theta(B)[x^{\pm1}; \al]_p$ is a Poisson-Laurent extension of $\theta(B)$.

{\rm (d)} Every $\HH$-stable Poisson prime ideal of $\Rhat$ is induced from an $\HH$-stable Poisson prime ideal of $\theta(B)$.
\eco

\begin{proof} (a) Obviously $\theta$ must map the polynomial ring $B[y]$ isomorphically onto a polynomial ring $\theta(B)[\theta(y)] = \theta(B)[x]$. In view of \prref{pCauchonmap}, $\theta(B)$ is a Poisson subalgebra of $\Rhat$, and $\{x, a\} = \al(a) x$ for all $a \in \theta(B)$, by definition of $\al$. Thus, $\theta(B)[x]$ is a Poisson-Ore extension of the form $\theta(B)[x; \al]_p$.

(b) It is clear that the subalgebra $T := \theta(B)[x; \al]_p$ is $\HH$-stable, and that the $\HH$-action on $T$ is rational. Once we show that $\al = (h_\ci \cdot)|_{\theta(B)}$, we will have that $T$ is a Poisson-Cauchon extension.

If $b \in B$ is homogeneous and nonzero, then $\theta(b)$ is homogeneous of the same degree, so $b$ and $\theta(b)$ lie in the same $\HH$-eigenspace of $\Rhat$. Since this $\HH$-eigenspace is also an $\h$-eigenspace, $b$ and $\theta(b)$ are $\h$-eigenvectors with the same $\h$-eigenvalue. Hence, there is some $\mu \in \KK$ such that $h_\ci \cdot b = \mu b$ and $h_\ci \cdot \theta(b) = \mu \theta(b)$, and consequently
$$
h_\ci \cdot \theta(b) = \theta(\mu b) = \theta \sig(b) = \al \theta(b).
$$
Therefore $(h_\ci \cdot)|_{\theta(B)} = \al$, as required.

(c) This is clear from the isomorphism $\theta$ and part (a).

(d)  Since $\al = (h_\ci \cdot)|_{\theta(B)}$ and $h_\ci \cdot x = \la_\ci x$ with $\la_\ci \ne 0$, part (d) follows from \cite[Lemma 1.1]{GLaun}.
\end{proof}

\bth{PoissonLLR}
Let $R = B[x; \sig, \de]_p$ be a Poisson-Cauchon extension. If $B$ is a noetherian $\HH$-Poisson-UFD, then so is $R$.
\eth

\begin{proof} Let $\theta$ and $\Rhat$ be as \prref{pCauchonmap} and $\al$ as in \coref{more.pCauchon}. Set $A := \theta(B)$ and $T := A[x; \al]_p$, and recall from \coref{more.pCauchon}(c) that $\Rhat = A[x^{\pm1}; \al]_p$. Since $B$  is an $\HH$-P-UFD, so is $A$. We first establish the following

{\bf Claim}: Suppose $u$ is a homogeneous Poisson-prime element of $A$ and $s\ge 0$ is minimal such that $v := ux^s \in R$. Then $v$ is a homogeneous Poisson-prime element of $R$.

Homogeneity of $v$ is clear from the homogeneity of $u$ and $x$. Now $u = \theta(u_0)$ for some homogeneous element $u_0 \in B$, and $u_0$ is an $\h$-eigenvector by \eqref{hHeigenvectors}, hence also a $\sig$-eigenvector. Thus, $u$ is an $\al$-eigenvector, and so $\{x,u\} = \al(u) x \in \KK ux$. Since $u$ is a Poisson-normal element of $A$, it follows that $u$ is a Poisson-normal element of $T$. Moreover, $T/uT \cong (A/uA)[x]$, which is a domain. Therefore $u$ is a homogeneous Poisson-prime element of $T$.

Since the Poisson $\HH$-prime ideal $uT$ of $T$ does not contain $x$, we see that $u \Rhat$ is a Poisson $\HH$-prime ideal of $\Rhat$. But $u \Rhat = v \Rhat$, so $I := v \Rhat \cap R$ is a Poisson $\HH$-prime ideal of $R$. We show that $I = vR$, which yields the claim.

Obviously $vR \subseteq I$, so consider an element $w \in I$. Then $w \in v \Rhat$. Choose $t \ge 0$ minimal such that $w x^t \in vR$, and write $w x^t = v r$ with $r \in R$. Suppose $t > 0$, whence $x \mid vr$ in $R$. We cannot have $r = r_1 x$ with $r_1 \in R$, since that would imply $w x^{t-1} = v r_1 \in vR$, contradicting the minimality of $t$. Thus $x \nmid r$ in $R$, and consequently $x \mid v$ in $R$. If $s > 0$, this would imply $u x^{s-1} = v x^{-1} \in R$, contradicting the minimality of $s$. Hence, $s = 0$, and so $u \in R$ with $x \mid u$. On the other hand, since $u \in \theta(B)$, \eqref{pCmap} then implies $u \in B$. This is impossible, because $u \ne 0$ and $x \mid u$ in $R$. Therefore $t = 0$ and $w \in vR$, completing the proof of the claim.

To prove that $R$ is an $\HH$-P-UFD, consider a nonzero Poisson $\HH$-prime ideal $J$ of $R$. 

If $x \notin J$, then $J$ is disjoint from $\{ x^n \mid n \in \Znn \}$, and so $J \Rhat$ is a nonzero Poisson $\HH$-prime ideal of $\Rhat$, with $J \Rhat \cap R = J$. In view of \coref{more.pCauchon}(d), $J \Rhat \cap A$ is a nonzero Poisson $\HH$--prime ideal of $A$. Since $A$ is an $\HH$-P-UFD, $J \Rhat \cap A$ contains a homogeneous Poisson-prime element $u$ of $A$. Let $s \ge 0$ be minimal such that $v := u x^s \in R$, and note that $v \in J \Rhat \cap R = J$. By the Claim, $v$ is a homogeneous Poisson-prime element of $R$, so we are done in this case.

Now assume that $x \in J$. If $\de = 0$, then we are done because $x$ is a homogeneous Poisson-prime element of $R$, so we may assume that $\de \ne 0$. Since $\de(b) = \{x,b\} - \sig(b) x$ for $b \in B$, we see that $\de(B) \subseteq J$. Consequently, $J \cap B \ne 0$.

Note that we have an $\HH$-equivariant Poisson algebra homomorphism
$$
\psi : A \xrightarrow{\;\theta^{-1}\;} B \xrightarrow{\;\text{incl}\;} R \xrightarrow{\;\text{quo}\;} R/J.
$$
Hence, $J' := \ker \psi$ is a Poisson $\HH$-prime ideal of $A$. Since $J \cap B \ne 0$, we see that $J' \ne 0$. Thus, $J' $ contains a homogeneous Poisson-prime element $u$ of $A$. 
Let $s \ge 0$ be minimal such that $v := u x^s \in R$. Then $v$ is a homogeneous Poisson-prime element of $R$, and we will be done if $v \in J$.

Set $z := \theta^{-1}(u) \in B$, and let $d \ge 0$ be maximal such that $\de^d(z) \ne 0$. Then
$$
u = \theta(z) = \sum_{n=0}^d \, \frac1{n!} \left( \frac{-1}{\la_\ci} \right)^n \de^n(z) x^{-n},
$$
and since $\de^d(z) \ne 0$ we see that $s = d$. Consequently,
\begin{equation}
\label{vsum}
v = u x^d = \sum_{n=0}^d \, \frac1{n!} \left( \frac{-1}{\la_\ci} \right)^n \de^n(z) x^{d-n}.
\end{equation}
As $u \in \ker \psi$, we have $z \in J$. Since also $\de(B) \subseteq J$, we conclude from \eqref{vsum} that $v \in J$. This completes the proof.
\end{proof}

\subsection{Poisson-prime elements in Poisson-Cauchon extensions}
\label{PpPCext}

\ble{p-REinv}
Let $R = B[x; \sig, \de]_p$ be a Poisson-Cauchon extension, where $B$ is a noetherian $\HH$-Poisson-UFD. Let $E \subseteq B$ be the multiplicative set generated by all homogeneous Poisson-prime elements of $B$, and assume there exist proper nonzero $\HH$-stable Poisson ideals in $R[E^{-1}]$. Then there exists a unique homogeneous element $d \in B[E^{-1}]$ such that $d = 0$ or $d$ has the same degree as $x$ and
\begin{equation}
\label{delbd}
\de(b) = \{d,b\} - \sig(b) d, \; \; \forall b \in B[E^{-1}].
\end{equation}
Moreover, $R[E^{-1}](x - d)$ is the unique nonzero Poisson $\HH$-prime ideal of $R[E^{-1}]$.
\ele

\begin{proof} \leref{torusPprime} ensures that $R[E^{-1}]$ has at least one nonzero Poisson $\HH$-prime ideal $P$. Now $P \cap B$ is a Poisson $\HH$-prime ideal of $B$, disjoint from $E$. Since $B$ is an $\HH$-P-UFD, $P \cap B = 0$, whence $P \cap B[E^{-1}] = 0$.

Following the proof of \cite[Proposition 1.2]{GLaun}, we obtain an element $d \in B[E^{-1}]$ such that $h \cdot d = \chi_x(h) d$ for all $h \in \HH$ and \eqref{delbd} holds. Now $R[E^{-1}] = B[E^{-1}] [y; \sig]_p$, where $y := x - d$. Since $B[E^{-1}] [y^{\pm1}; \sig]_p$ has no nonzero Poisson $\HH$-prime ideals (by \cite[Lemma 1.1]{GLaun}), we conclude that $P = R[E^{-1}] y$. Thus, $R[E^{-1}] y$ is the  unique nonzero Poisson $\HH$-prime ideal of $R[E^{-1}]$.

If $d'$ is any homogeneous element of $B[E^{-1}]$ satisfying the stated properties of $d$, then $R[E^{-1}] = B[E^{-1}] [x - d'; \sig]_p$ and $P = R[E^{-1}] (x - d')$. Consequently, $d - d' \in P$, and therefore $d - d' = 0$.
\end{proof}

\bco{deg.le1} 
Let $R = B[x; \sig, \de]_p$ be a Poisson-Cauchon extension, where $B$ is a noetherian $\HH$-Poisson-UFD. Then all homogeneous Poisson-prime elements of $R$ have degree at most $1$ in $x$. Up to taking associates, there is at most one homogeneous Poisson-prime element of $R$ which does not lie in $B$ {\rm(}i.e., has degree $1$ in $x${\rm)}.
\eco

\begin{proof} If $v$ is a homogeneous Poisson-prime element of $R$ such that 
$v \notin B$, then $Rv$ is disjoint from the multiplicative set $E$ of \leref{p-REinv}. Consequently, $v$ is a homogeneous Poisson-prime element of $R[E^{-1}]$, and so $R[E^{-1}]v$ is a nonzero Poisson $\HH$-prime ideal of $R[E^{-1}]$. Hence, there exists $d\in B[E^{-1}]$ as in \leref{p-REinv}, and $v$ is an 
associate of the prime element $x - d \in R[E^{-1}]$ (as prime 
elements of $R[E^{-1}]$). This implies that $v$ has degree at most $1$ in $x$. Since $Rv$ is a prime ideal of $R$ disjoint from $E$, we have
$$
Rv = (R[E^{-1}]v) \cap R = 
\big( R[E^{-1}] (x-d) \big) \cap R.
$$

If $w$ is any other homogeneous Poisson-prime element of $R$ that is not in $B$, the same argument as above shows that $Rw = \big( R[E^{-1}] (x-d) \big) \cap R$, and therefore $Rw= Rv$.
\end{proof}

For an ideal $Q$ of a polynomial ring $R = B[x]$, denote the ideal of its leading coefficients
$$
\lc(Q) := \{ b \in B \mid \exists a \in Q, \; m \in \Znn \; \text{such that} \; a - bx^m \in Bx^{m-1} +\cdots+ B \}.
$$

\ble{thetaJRhat}
Let $R = B[x; \sig, \de]_p$ be a Poisson-Cauchon extension, where $B$ is a noetherian $\HH$-Poisson-UFD, let $\theta : B[y^{\pm1}; \sig]_p \rightarrow \Rhat$ be the Poisson isomorphism of Proposition {\rm\ref{ppCauchonmap}}, and let $J$ be a Poisson $\HH$-prime ideal of $B$. If $J$ is a height one prime ideal of $B$, then $\lc(\theta(J) \Rhat \cap R) = J$, and $\theta(J) \Rhat \cap R$ is a Poisson $\HH$-prime ideal of $R$ as well as a height one prime ideal of $R$.
\ele

\begin{proof} Set $Q := \theta(J) \Rhat \cap R$, and note  that $Q$ is a nonzero prime ideal of $R$. Since $\theta$ is an $\HH$-equivariant Poisson isomorphism, we see that $Q$ is $\HH$-stable and that $\{ \theta(b), \theta(J) \Rhat \} \subseteq \theta(J) \Rhat$ for all $b \in B$. Moreover, $\{ x, \theta(J) \} \subseteq \theta(J) x$ because of \prref{pCauchonmap} and \eqref{Hsig.stable}. It follows that $\theta(J) \Rhat$ is a Poisson ideal of $\Rhat$, whence $Q$ is a Poisson ideal of $R$. Thus, $Q$ is a Poisson $\HH$-prime ideal of $R$.

We have $J = uB$ for some homogeneous Poisson-prime element $u \in B$, and so
$$
\theta(J) \Rhat = \theta(u) \Rhat = \sum_{n \in \Zset} \theta(u) Bx^n .
$$
Since $\theta(u) \in u + \sum_{i < 0} Bx^i$, it is clear that $\lc(Q) = uB = J$.

Because $R$ is an $\HH$-P-UFD (\thref{PoissonLLR}), $Q$ contains a homogeneous Poisson-prime element $p$ of $R$, and hence it contains the nonzero Poisson $\HH$-prime ideal $P := pR$. Obviously $x \notin Q$, since $1 \notin J = \lc(Q)$. As a result, we have prime ideals
\begin{equation}
\label{PQRhat}
0 \subsetneq P \Rhat \subseteq Q \Rhat = \theta(J) \Rhat
\end{equation} in $\Rhat$. By \coref{more.pCauchon}(d), $P \Rhat = P' \Rhat$ for some nonzero Poisson $\HH$-prime ideal $P'$ of $\theta(B)$, and \eqref{PQRhat} implies that $P' \subseteq \theta(J)$. But $\theta(J)$ is a height one prime ideal of $\theta(B)$, so we must have $P' = \theta(J)$, whence $P \Rhat = Q \Rhat$. We conclude that $Q = P$, and therefore $Q$ has height one.
\end{proof}

\ble{locnilp}
{\rm (a)} Let $B$ be a domain of characteristic zero and $\de$ a locally nilpotent derivation on $B$. If $a \in B$ and $\de(a) \in aB$ or $\de(a) \in Ba$, then $\de(a) = 0$. 

{\rm (b)} Assume that $R = B[x; \sig, \de]_p$ is a Poisson-Cauchon extension of a domain $B$, and that $\de$ is inner as a Poisson $\sig$-derivation, that is, there exists $c \in B$ such that $\de(b) = \{c,b\} - \sig(b) c$ for all $b \in B$. 

Then $\de = 0$. If, in addition, $c$ is a homogeneous element of $B$ of the same $\xh$-degree as $x$, then $c = 0$.
\ele

\begin{proof} (a) is proved by the argument of \cite[Lemme 7.2.3.2]{Ric}, as follows. Suppose that $\de(a) = ab$ for some $b\in B$ and that $\de(a) \ne 0$. Let $m,n\in \Zset_{>0}$ be minimal such that $\de^m(a) = \de^n(b) = 0$. By the Leibniz Rule,
\begin{align*}
\de^{m+n-1}(a) &= \de^{m+n-2}(ab) = \sum_{i=0}^{m+n-2} \tbinom{m+n-2}{i} \de^{m+n-2-i}(a) \de^i(b)  \\
&= \tbinom{m+n-2}{n-1} \de^{m-1}(a) \de^{n-1}(b) \ne 0,
\end{align*}
because $\de^i(b) = 0$ for $i \ge n$ and $\de^{m+n-2-i}(a) = 0$ for $i \le n-2$. This forces $m+n-1 \le m-1$ and so $n=0$, contradicting our assumptions. Thus, $\de(a) = 0$.

(b) Write $c = c_0 + \cdots + c_m$ where the $c_i$ are homogeneous elements (possibly zero) of $B$ with distinct $\xh$-degrees $\mu_i$ and $\mu_0 = \chi_x$. If $b \in B$ is homogeneous of degree $\rho$, then $\de(b)$ is homogeneous of degree $\chi_x + \rho$ by \eqref{chisigdel}. On the other hand, each $\{c_i,b\} - \sig(b) c_i$ is homogeneous of degree $\mu_i + \rho$, so $\de(b) = \{c_0,b\} - \sig(b) c_0$. Since this holds for all homogeneous $b \in B$, it holds for all $b \in B$.

The elements $c_0$ and $x$ lie in the same $\HH$-eigenspace of $R$, so they also lie in the same $\h$-eigenspace by \eqref{hHeigenspaces}. Hence, $\sig(c_0) = h_\ci \cdot c_0 = \la_\ci c_0$. Now $\de(c_0) =  \{c_0,c_0\} - \sig(c_0) c_0 = - \la_\ci c_0^2$. Part (a) implies that $\de(c_0) = 0$, and so $c_0 = 0$. Thus, $\de = 0$ and the lemma is proved.
\end{proof}

\bth{u-i,ii}
Let $R = B[x; \sig, \de]_p$ be a Poisson-Cauchon extension, where $B$ is a noetherian $\HH$-Poisson-UFD. Let $u$ be a homogeneous Poisson-prime element of $B$, and let $\mu \in \KK$ such that
$$
\sig(u) = h_\ci \cdot u = \mu u.
$$
Then exactly one of the following two situations occurs:

{\rm (i)} The element $u$ remains a Poisson-prime element of $R$.

In this case, $\{u,x\} = - \mu ux$ and $\de(u) = 0$.

{\rm (ii)} There exists a unique element $c_\ci \in B$ such that $v := ux - c_\ci$ is a homogeneous Poisson-prime element of $R$ {\rm(}in particular, $u \nmid c_\ci${\rm)}. Moreover, $c_\ci$ is homogeneous.

In this case, $\de$ is given by
\begin{equation}
\label{iide}
\de(b) = \{ u^{-1} c_\ci, b\} - \sig(b) (u^{-1} c_\ci), \; \; \forall b \in B,
\end{equation}
and $\{v,-\}$ is determined by
\begin{equation}
\label{vPnormalize}
\{v,x\} = - \mu xv, \qquad \{v,b\} = (\partial_u + \sig)(b) v, \; \; \forall b \in B,
\end{equation}
where $\partial_u$ is the derivation on $B$ given in \eqref{defdelc}. Furthermore,
$$
\{u,c_\ci\} = - (\mu + \la_\ci) c_\ci u, \qquad \de(u) = \la_\ci c_\ci \ne 0, \qquad \de(c_\ci) = 0.
$$
\eth

We note that by \coref{deg.le1}, the situation (ii) cannot simultaneously occur for two homogeneous Poisson-prime elements $v$ of $R$, since such elements $v$ must be associates of each other.

\begin{proof} We first show that $\de(u) = 0$ in situation (i) while $\de(u) \ne 0$ in situation (ii), and so the two situations cannot occur simultaneously.

In situation (i), $u$ is Poisson-normal in $R$, whence the element $\{x,u\} = \mu ux + \de(u)$ is divisible by $u$ in $R$. Consequently, $\de(u)$ is divisible by $u$ in $B$, and so $\de(u) = 0$ by \leref{locnilp}(a).

In situation (ii), the element $\{x,v\}$ is divisible by $v$ in $R$, and we obtain
$$
(\mu ux + \de(u)) x - (\sig(c_\ci) x + \de(c_\ci)) = \{x,v\} = (ax + b)(ux - c_\ci)
$$
for some $a,b \in B$, whence
$$
au = \mu u, \qquad bc_\ci = \de(c_\ci), \qquad bu - ac_\ci = \de(u) - \sig(c_\ci).
$$
Then $a = \mu$ and $b = \de(c_\ci) = 0$ by \leref{locnilp}(a), so
$$
\de(u) = \sig(c_\ci) - \mu c_\ci .
$$
Since $v$ and $ux$ are homogeneous, they have the same $\HH$-eigenvalue, and $c_\ci$ is homogeneous with that $\HH$-eigenvalue. Hence, by \eqref{hHeigenvectors}, $v$, $ux$, and $c_\ci$ are $h_\ci$-eigenvectors with the same $h_\ci$-eigenvalue, namely $\mu + \la_\ci$, and so $\sig(c_\ci) = h_\ci \cdot c_\ci = (\mu + \la_\ci) c_\ci$. Thus,
$$
\de(u) = \la_\ci c_\ci \ne 0.
$$
This verifies that situations (i) and (ii) are disjoint. 

We now show that the cases $\de(u) = 0$ and $\de(u) \ne 0$ lead to situations (i) and (ii), respectively. Let $\theta : B[y^{\pm1}; \sig]_p \rightarrow \Rhat$ be the Poisson isomorphism of \prref{pCauchonmap}. 

If $\de(u) = 0$, then $\{x,u\} = \mu ux$, and the Poisson-normality of $u$ in $B$ implies that $u$ is Poisson-normal in $R$. Moreover, $\theta(u) = u$, whence
$$
\theta(uB) \Rhat \cap R = uB[x^{\pm1}] \cap B[x] = u B[x] = uR.
$$
\leref{thetaJRhat} now implies that $uR$ is a prime ideal of $R$. Therefore $u$ is a Poisson-prime element of $R$, and situation (i) holds.

Assume that $\de(u) \ne 0$ for the remainder of the proof. By \leref{thetaJRhat}, $\theta(uB) \Rhat \cap R$ is a height one Poisson $\HH$-prime ideal of $R$, and by \thref{PoissonLLR}, $R$ is an $\HH$-P-UFD. Therefore there exists a homogeneous Poisson-prime element $v$ of $R$ such that
$$
\theta(uB) \Rhat \cap R = vR.
$$
Denote the leading coefficient of $v$ (as a polynomial in $B[x]$) by $u'$. Then \leref{thetaJRhat} implies that
$$
uB = \lc( \theta(uB) \Rhat \cap R) = \lc(vR) = u' B.
$$
Consequently, $u'$ is a homogeneous Poisson-prime element of $B$ which is an associate of $u$ (in $B$). Thus, after multiplying $v$ by a homogeneous unit of $B$ we can assume that
$$
u' = u.
$$

By what we proved at the beginning, $v \ne u$, since that would imply $\de(u) = 0$. On the other hand, $v$ has degree at most $1$ in $x$, by \coref{deg.le1}. Therefore
$$
v = u x - c_\ci
$$
for some $c_\ci \in B$. Uniqueness of $c_\ci$ holds by the comment ahead of the proof.

We are now in situation (ii), and it remains to verify the associated conditions stated in the theorem. We have already seen that $\de(c_\ci) = 0$ and $\de(u) = \la_\ci c_\ci \ne 0
$, while $\{x,v\} = \mu xv$.

For any $b \in B$, the element
$$
\{v,b\} = \partial_u(b) ux + u( \sig(b) x + \de(b)) - \{c_\ci,b\}
$$
is divisible by $v$ in $R$, whence
$$
\{v,b\} = (\partial_u(b) + \sig(b)) v \qquad \text{and} \qquad u \de(b) - \{c_\ci,b\} = - (\partial_u(b) + \sig(b)) c_\ci.
$$
The first equation gives us the remainder of \eqref{vPnormalize}, and the second equation together with a straightforward calculation yields \eqref{iide}.

We saw above that $\sig(c_\ci) = (\mu + \la_\ci) c_\ci$. Applying \eqref{iide} to $c_\ci$, we obtain
$$
0 = \de(c_\ci) = \{ u^{-1}c_\ci, c_\ci\} - \sig(c_\ci) u^{-1} c_\ci = - u^{-2} \{u,c_\ci\} c_\ci - \sig(c_\ci) u^{-1} c_\ci,
$$
and consequently $\{ u,c_\ci\} = - \sig(c_\ci) u = - (\mu + \la_\ci) c_\ci u$. This completes the proof.
\end{proof}

\bth{BlistR}
Let $R = B[x; \sig, \de]_p$ be a Poisson-Cauchon extension, where $B$ is a noetherian $\HH$-Poisson-UFD, and let $\{ u_i \mid i \in I\}$ be a list of the homogeneous Poisson-prime elements of $B$ up to taking associates. Then there are the following three possibilities for a list of the homogeneous Poisson-prime elements of $R$ up to taking associates:

{\rm (i)} $\{ u_i \mid i \in I, \; i \ne i_0\} \sqcup \{ v_{i_0} := u_{i_0} x - c_\ci \}$, for some $i_0 \in I$ and $c_\ci \in B$ such that $u_{i_0}^{-1} c_\ci$ is a nonzero homogeneous element of $B[u_{i_0}^{-1}]$ with the same $\xh$-degree as $x$.

{\rm (ii)} $\{ u_i \mid i \in I\} \sqcup \{x\}$.

{\rm (iii)} $\{ u_i \mid i \in I\}$.
\eth

Before proving \thref{BlistR}, we record the following information for the three cases of the theorem. This proposition is a direct consequence of \thref{u-i,ii}.

\bpr{info-i,ii,iii}
Assume the setting of Theorem {\rm\ref{tBlistR}}. For $i \in I$, let $\mu_i \in \KK$ such that $\sig(u_i) = h_\ci \cdot u_i = \mu_i u_i$.

{\rm (a)} In case {\rm (i)} of Theorem {\rm\ref{tBlistR}}, we have
$$
\de(b) = \{ u_{i_0}^{-1} c_\ci, b\} - \sig(b) (u_{i_0}^{-1} c_\ci), \; \; \forall b \in B, 
$$
$$
\begin{aligned}
\de(u_i) &= 0, \; \; \forall i \in I \setminus \{i_0\},  &\de(u_{i_0}) &= \la_\ci c_\ci \ne 0,  \\
\de(c_\ci) &= 0,  &\qquad\{u_{i_0}, c_\ci\} &= - (\mu_{i_0} + \la_\ci) c_\ci u_{i_0} .
\end{aligned}
$$
Furthermore,
\begin{align*}
\{u_i,x\} &= - \mu_i u_i x, \; \; \forall i \in I \setminus \{i_0\}, \qquad \{v_{i_0},x\} = - \mu_{i_0} v_{i_0} x, \\
\{v_{i_0},b\} &= (\partial_{u_{i_0}} + \sig)(b) v_{i_0}, \; \; \forall b \in B.
\end{align*}

{\rm (b)}  In case {\rm (ii)} of Theorem {\rm\ref{tBlistR}}, we have $\de = 0$,
$$
\{x,b\} = \sig(b) x, \; \; \forall b \in B, \qquad \text{and} \qquad \{u_i,x\} = - \mu_i u_i x, \; \; \forall i \in I.
$$

{\rm (c)} In case {\rm (iii)} of Theorem {\rm\ref{tBlistR}}, we have $\de \ne 0$,
$$
\de(u_i) = 0 \quad \text{and} \quad \{u_i,x\} = - \mu_i u_i x, \; \; \forall i \in I.
$$
\epr

\begin{proof}[Proof of Theorem {\rm\ref{tBlistR}}] \coref{deg.le1} implies that there is at most one index $i \in I$ for which $u_i$ falls into situation (ii) of \thref{u-i,ii}. Thus, we have two cases:
\begin{enumerate}
\item For all $i \in  I$, the element $u_i$ remains Poisson-prime in $R$.
\item There is an index $i_0 \in I$ such that $u_i$ remains Poisson-prime in $R$ for all $i \ne i_0$ and there exists $c_\ci \in B$ such that $v_{i_0} := u_{i_0} x - c_\ci$ is a homogeneous Poisson-prime element of $R$.
\end{enumerate}
Note that it is possible that $I = \varnothing$, in which case (1) holds.

Case (1). In this case, $\de(u_i) = 0$ for all $i \in I$, by \thref{u-i,ii}.  Any homogeneous Poisson-prime element of $R$ which has degree $0$ in $x$ is clearly Poisson-prime in $B$ and so is an associate of one of the elements $u_i$. If there are no other homogeneous Poisson-prime elements in $R$, we are in the situation (iii).

Suppose that $R$ has a homogeneous Poisson-prime element $v$ that is not an associate of any of the $u_i$. By \coref{deg.le1}, $v$ has degree $1$ in $x$ and any homogeneous Poisson-prime element of $R$ is an associate of either $v$ or one of the $u_i$. Write $v = ux - c_\ci$ for some homogeneous elements $u,c_\ci \in B$, with $u \ne 0$. For any $b \in B$, the element
$$
\{v,b\} = \{u,b\} x + u( \sig(b) x + \de(b)) - \{c_\ci,b\}
$$
is divisible by $v$ in $R$, whence $\{u,b\} + u \sig(b)$ is divisible by $u$ in $B$, and so $\{u,b\}$ is divisible by $u$. Thus, $u$ is a Poisson-normal element of $B$. 

By \prref{factorHPUFD}, $u = wu'$ where $w$ is a homogeneous unit in $B$ and $u'$ is either $1$ or a product $u_{i_1} u_{i_2} \cdots u_{i_m}$ for some $i_j \in I$. After replacing $v$ by $w^{-1} v$, we may assume that $u = u'$. Since $\de(u_i) = 0$ for all $i \in I$, we now have $\de(u) = 0$. Let $\mu \in \KK$ be such that $\sig(u) = h_\ci \cdot u = \mu u$, and observe that $v$ is an $h_\ci$-eigenvector with eigenvalue $\mu + \la_\ci$, whence $\sig(c_\ci) = h_\ci \cdot c_\ci = (\mu + \la_\ci) c_\ci$. Now the element
$$
\{x,v\} = \sig(u) x^2 - \sig(c_\ci) x - \de(c_\ci) = \mu ux^2 - (\mu + \la_\ci) c_\ci x - \de(c_\ci)
$$
is divisible by $v$ in $R$, and hence
$$
\mu ux^2 - (\mu + \la_\ci) c_\ci x - \de(c_\ci) = (ux - c_\ci) (\mu x + c)
$$
for some $c \in B$. Consequently, $- \la_\ci c_\ci = u c$ and so $u \mid c_\ci$ in $B$, whence $u \mid v$ in $R$. It follows that $u$ must be a unit in $B$. By our choices above, $u = 1$.

Therefore $v = x - c_\ci$. The homogeneity of $v$ implies $\chi_{c_\ci} = \chi_x$. For any $b \in B$, the element
$$
\{v,b\} = \sig(b) x + \de(b) - \{c_\ci,b\}
$$
is divisible by $v$, whence $\{v,b\} = \sig(b) (x - c_\ci)$ and consequently
$$
\de(b) = \{c_\ci,b\} - \sig(b) c_\ci .
$$
\leref{locnilp}(b) implies that $c_\ci = 0$ and $\de = 0$. Therefore $v = x$, and we are in the situation (ii).

Case (2). As before, any homogeneous Poisson-prime element of $R$ which has degree $0$ in $x$ is an associate of one of the elements $u_i$. By \coref{deg.le1}, any other homogeneous Poisson-prime element of $R$ is an associate of $v_{i_0}$. It follows from \thref{u-i,ii} that $u_{i_0}$ cannot be Poisson-prime in $R$, since the cases in that theorem are mutually exclusive. Thus, in this case we are in situation (i).
\end{proof}

\sectionnew{Iterated Poisson-Ore extensions}
\label{iterp-o}

We now introduce the class of iterated Poisson-Ore extensions that we call \emph{Poisson nilpotent algebras}, and we classify the homogeneous Poisson-prime elements in such algebras.

\subsection{Poisson-CGL extensions}
\label{PCGLext}
We focus on iterated Poisson-Ore extensions
\begin{equation} 
\label{itpOre}
R := \KK[x_1]_p [x_2; \sig_2, \delta_2]_p \cdots [x_N; \sig_N, \delta_N]_p ,
\end{equation}
where it is taken as implied that $\KK[x_1]_p = \KK[x_1; \sig_1, \delta_1]_p$ with $\sig_1 = \de_1 = 0$. The integer $N$ will be called the \emph{length} of the extension $R$. For $k \in [0,N]$, set
$$
R_k := \KK [ x_1, \dots, x_k ] = \KK[x_1]_p [x_2; \sig_2, \delta_2]_p \cdots [x_k; \sig_k, \delta_k]_p .
$$
In particular, $R_0 = \KK$.

The Poisson analog of the concept of a CGL extension introduced in \cite[Definition 3.1]{LLR} is defined as follows.

\bde{PCGL} An iterated Poisson-Ore extension \eqref{itpOre}
is called a \emph{Poisson-CGL} (or \emph{P-CGL}) \emph{extension} 
 if it is equipped with a rational Poisson action of a torus $\HH$ such that
\begin{enumerate}
\item[(i)] The elements $x_1, \ldots, x_N$ are $\HH$-eigenvectors.
\item[(ii)] For every $k \in [2,N]$, $\de_k$ is a locally nilpotent Poisson 
$\sig_k$-derivation of the algebra $R_{k-1}$. 
\item[(iii)] For every $k \in [1,N]$, there exists $h_k \in \h = \Lie \HH$ such that 
$\sig_k = (h_k \cdot)|_{R_{k-1}}$ and the $h_k$-eigenvalue of $x_k$, to be denoted by $\la_k$, is nonzero.
\end{enumerate}
Note that each $x_k$ is an $\h$-eigenvector, by virtue of (i) and \eqref{hHeigenvectors}.

Whenever we refer to a Poisson-CGL extension $R$, we take the associated torus and its Lie algebra to be denoted $\HH$ and $\h$ unless otherwise specified. For any iterated Poisson-Ore extension \eqref{itpOre}, there is a canonical (and maximal in a suitable sense) choice of a $\KK$-torus acting rationally on $R$ by Poisson automorphisms such that the variables $x_k$ are eigenvectors. (See \S\ref{maxtori}.)

Conditions (i) and (iii) above imply that
\begin{equation}
\label{sigkxj}
\sig_k(x_j) = h_k \cdot x_j = \la_{kj} x_j \; \; \text{for some} \; \la_{kj} \in \KK, \; \; \forall 1 \le j < k \le N.
\end{equation}
We then set $\la_{kk} := 0$ and $\la_{jk} := -\la_{kj}$ for $j < k$. This gives rise to a skew-symmetric matrix $\lab := (\la_{kj}) \in M_N(\KK)$ and the corresponding skew-symmetric bicharacter $\Om_\lab : \Zset^N \times \Zset^N \rightarrow \KK$ from \eqref{Om}.

For each $k \in [1,N]$, the algebra $R_k$ is a Poisson-Cauchon extension of $R_{k-1}$. In particular, it follows from \eqref{aldel} and \eqref{hdotdel} that
\begin{gather}
\sig_k \ci \de_k = \de_k \ci (\sig_k + \la_k), \; \; \forall k \in [2,N]  \label{sigkdelk}  \\
(h\cdot)|_{R_{k-1}} \ci \de_k = \chi_{x_k}(h) \de_k \ci (h\cdot)|_{R_{k-1}}, \; \; \forall h \in \HH, \; k \in [2,N].  \label{hdotdelk}
\end{gather}
\ede

By induction on \thref{PoissonLLR}, we obtain

\bth{PCGL-UFD}
Every Poisson-CGL extension $R$ is a noetherian $\HH$-Poisson-UFD.
\eth

Any Poisson affine space algebra $R = \KK[x_1,\dots,x_N]$ (Examples \ref{PaffPtorus}) is a Poisson-CGL extension with all $\de_k = 0$, relative to the standard action of $\HH = (\kx)^N$ on $R$.
Many examples of Poisson-CGL extensions occur as semiclassical limits of noncommutative CGL extensions, as shown, for instance, in \cite[\S\S 2.2--2.7, 4.1--4.6]{GLaun}. The semiclassical limit of the standard (uniparameter) quantized coordinate ring of $m\times n$ matrices yields the following.

\bex{OMmn}
For positive integers $m$ and $n$, let $\OMmn$ be the ring of polynomial functions on the matrix variety $M_{m,n}(\KK)$, written as a polynomial ring 
$$
\OMmn = \KK[ t_{ij} \mid 1 \le i \le m, \; 1 \le j \le n] .
$$
There is a standard Poisson structure on $\OMmn$ such that
\begin{align*}
\{t_{ij}, t_{kj}\} &= t_{ij} t_{kj},  &&\quad \text{for} \; i < k,  \\
\{t_{ij}, t_{il}\} &= t_{ij} t_{il},  &&\quad \text{for} \; j < l,  \\
\{t_{ij}, t_{kl}\} &= 0,  &&\quad \text{for} \; i < k, \; j > l,  \\
\{t_{ij}, t_{kl}\} &= 2 t_{il} t_{kj},  &&\quad \text{for} \; i < k, \; j < l.
\end{align*}
This Poisson algebra is an iterated Poisson-Ore extension
$$
\OMmn =  \KK[x_1]_p [x_2; \sig_2, \delta_2]_p \cdots [x_N; \sig_N, \delta_N]_p ,
$$
where $N = mn$ and $x_{(r-1)n+c} = t_{rc}$. It is easy to write explicit formulas for the maps $\sig_k$ and $\de_k$, and to check that each $\de_k$ is locally nilpotent. For later reference, we note that the scalars $\la_{kj}$ are given as follows:
$$
\la_{(r-1)n+c, (r'-1)n+c'} = \begin{cases} 
\sign(r'-r), &\text{if} \; \; c = c', \\
\sign(c'-c), &\text{if} \; \; r = r',  \\
0, &\text{otherwise},
\end{cases}
\; \; \forall r,r' \in [1,m], \; \; c,c' \in [1,n].
$$

There is a rational Poisson action of the torus $\HH = (\kx)^{m+n}$ on $\OMmn$ such that
$$
(\xi_1, \dots, \xi_{m+n}) \cdot t_{rc} = \xi_r \xi_{m+c}^{-1} t_{rc}
$$
for all $(\xi_1, \dots, \xi_{m+n}) \in \HH$, $r \in [1,m]$, $c \in [1,n]$. Define
$$
h_{rc} := (0, \dots, 0, -1, 0,\dots, 0, 1, 0,\dots, 0) \; \in \; \h = \KK^{m+n} = \Lie \HH,
$$
where the entries $-1$ and $1$ reside in positions $r$ and $m+c$, respectively. Then $\sig_{(r-1)n+c} = (h_{rc} \cdot)$ and $h_{rc} \cdot t_{rc} = -2 t_{rc}$. In other words, for $k = (r-1)n+c$ we have $h_k = h_{rc}$ and the $h_k$-eigenvalue of $x_k$ is $\la_k = -2$. Thus, $\OMmn$ is a P-CGL extension.
\eex

\subsection{Poisson-prime elements in Poisson-CGL extensions}
\label{PprimePCGL}

The next theorem describes the homogeneous Poisson-prime elements of a Poisson-CGL extension $R = R_N$ iteratively from those of $R_{N-1}$. It shows that the situation (iii) from \thref{BlistR} never arises in the framework of Poisson-CGL extensions. Note that the group of units of $R$ is reduced to scalars. Thus, two prime elements of $R$ are associates if and only if they are scalar multiples of each other.

\bth{PCGL-list}
Let $R$ be an arbitrary Poisson-CGL extension of length $N$ as in \eqref{itpOre}.

{\rm (a)} Let $\{ u_i \mid i\in I\}$ be a list of the homogeneous Poisson-prime elements of $R_{N-1}$ up to scalar multiples. There are two possibilities for a list of the homogeneous Poisson-prime elements of $R$ up to scalar multiples:
\begin{enumerate}
\item[(i)] $\{u_i \mid i \in I \setminus \{i_0\} \} \sqcup \{ u_{i_0} x_N - c_\ci \}$ for some $i_0 \in I$ and $c_\ci \in R_{N-1}$ such that $u_{i_0}^{-1} c_\ci$ is a nonzero homogeneous element of $R_{N-1}[u_{i_0}^{-1}]$ with the same $\xh$-degree as $x_N$.
\item[(ii)] $\{u_i \mid i \in I\} \sqcup \{x_N\}$.
\end{enumerate}

{\rm (b)} Let $\de$ be a locally nilpotent derivation on $R$. If
\begin{equation}
\label{de0Pprime}
\de(u) = 0 \; \; \text{for all homogeneous Poisson-prime elements} \; u \; \text{of} \; R,
\end{equation}
then $\de = 0$.
\eth

The proof of \thref{PCGL-list} will be given in subsection \ref{proofT3.4}. We next derive a number of consequences from part (a) of \thref{PCGL-list}.

It follows from \thref{BlistR} (or \thref{PCGL-list}) that a Poisson-CGL extension $R$ has only a finite number of pairwise nonproportional homogeneous Poisson-prime elements. We will call this number the \emph{rank} of $R$. It also equals the number of Poisson $\HH$-prime ideals of height $1$ in $R$. For each $k \in [1,N]$, \thref{PCGL-list} in combination with \prref{info-i,ii,iii} implies that $\rank R_k = \rank R_{k-1}$ if $\de_k \ne 0$, while $\rank R_k = \rank R_{k-1} + 1$ if $\de_k = 0$. Thus,
\begin{equation}
\label{rankRk}
\rank R_k = \big| \{ j \in [1,k] \mid \de_j = 0 \} \big|, \; \; \forall k \in [1,N].
\end{equation}

Using \thref{PCGL-list}, we can describe the homogeneous Poisson-prime elements in each stage of any Poisson-CGL extension $R$, as follows.
The statement of the result involves the standard predecessor and successor functions, $p = p_\eta$ and $s = s_\eta$, for the level sets of a function $\eta : [1,N] \to \Zset$, defined as follows:
\begin{equation}
\label{pred.succ}
\begin{aligned}
p(k) &= \max \{ j <k \mid \eta(j) = \eta(k) \}, \\
s(k) &= \min \{ j > k \mid \eta(j) = \eta(k) \}, 
\end{aligned}
\end{equation}
where $\max \varnothing = -\infty$ and $\min \varnothing = +\infty$. Define corresponding order functions $O_\pm : [1,N] \rightarrow \Znn$ by
\begin{equation}
\label{O-+}
\begin{aligned}
O_-(k) &:= \max \{ m \in \Znn \mid p^m(k) \ne -\infty \},  \\
O_+(k) &:= \max \{ m \in \Znn \mid s^m(k) \ne +\infty \}.
\end{aligned}
\end{equation}

\bth{mainPprime}
Let $R$ be a Poisson-CGL extension of length $N$ as in \eqref{itpOre}. There exist a function $\eta : [1,N] \rightarrow \Zset$ and homogeneous elements
$$
c_k \in R_{k-1} \; \; \text{for all} \; \; k \in [2,N] \; \; \text{with} \; \; p(k) \ne - \infty
$$
such that the elements $y_1, \dots, y_N \in R$, recursively defined by
\begin{equation}
\label{ykdef}
y_k := \begin{cases}
y_{p(k)} x_k - c_k,  &\quad \text{if} \; p(k) \ne - \infty  \\
x_k,  &\quad \text{if} \; p(k) = - \infty,
\end{cases}
\end{equation}
are homogeneous and have the property that for every $k \in [1,N]$,
\begin{equation}
\label{yset1k}
\{ y_j \mid j \in [1,k], \; s(j) > k \}
\end{equation}
is a list of the homogeneous Poisson-prime elements of $R_k$ up to scalar multiples. In particular, $y_k$ is a homogeneous Poisson-prime element of $R_k$ as well as a prime element of $R$, for all $k \in [1,N]$.

The elements $y_1, \dots, y_N \in R$ with these properties are unique. The function $\eta$ satisfying the above conditions is not unique, but the partition of $[1,N]$ into the disjoint union of the level sets of $\eta$ is uniquely determined by $R$, as are the predecessor and successor functions $p$ and $s$. The function $p$ has the property that $p(k) = -\infty$ if and only if $\de_k = 0$.
\eth

The uniqueness statements in the final paragraph of the theorem are dependent on the given P-CGL extension presentation \eqref{itpOre} of $R$, not on $R$ as an algebra. Typically, $R$ will have many different P-CGL extension presentations, as we discuss in subsection \ref{symmP-CGL}.

Uniqueness of the elements $y_k$, $k \in [1,N]$, and of the level sets of $\eta$, follows at once from \thref{PCGL-list}. This uniqueness immediately implies the uniqueness of $p$ and $s$. From the final statement of the theorem together with \eqref{rankRk}, we see that the rank of $R$ is given by
\begin{equation}
\label{rankRs>N}
\rank R = | \eta([1,N]) | = \big| \{ j \in [1,N] \mid p(j) = -\infty \} \big| = \big| \{ j \in [1,N] \mid s(j) = +\infty \} \big| .
\end{equation}

\begin{proof}[Proof of Theorem {\rm\ref{tmainPprime}}] We define $\eta : [1,k] \rightarrow \Zset$ and elements $c_k \in R_{k-1}$ (when $p(k) \ne -\infty$) for $k = 1, \dots, N$. At each step, the new function $\eta$ will be an extension of the previous one, and so the corresponding new predecessor function will also be an extension of the previous one. However, the successor functions may change, so we will write $s_k$ for the successor function going with $\eta$ on $[1,k]$.

To start, set $\eta(1) := 1$. Note that $p(1) = -\infty$ and $s_1(1) = +\infty$. Moreover, $y_1 := x_1$ is the unique homogeneous Poisson-prime element of $R_1$ up to scalar multiples. This relies on the fact that the homogeneous elements of $R_1$ are just the scalar multiples of the monomials $x_1^m$, which follows from \eqref{hHeigenvectors} and the assumption that the $h_1$-eigenvalue of $x_1$ is nonzero.

Now let $1 < k \le N$, and assume that $\eta$ has been defined on $[1, k-1]$, together with elements $c_j \in R_{j-1}$ for $j \in [1, k-1]$ with $p(j) \ne -\infty$ and homogeneous elements $y_j \in R_j$ for $j \in [1, k-1]$ such that \eqref{ykdef} and \eqref{yset1k} hold. In particular, 
$$
\{ y_j \mid j \in [1,k-1], \; s_{k-1}(j) \ge k \}
$$
is a list of the homogeneous Poisson-prime elements of $R_{k-1}$ up to scalar multiples. There are two cases to consider, corresponding to situations (i), (ii) of \thref{PCGL-list}(a).

In the first case, there is some $j_0 \in [1, k-1]$ such that $s_{k-1}(j_0) \ge k$ and some homogeneous $c_k \in R_{k-1}$ such that
\begin{equation}
\label{P4.8}
\{ y_j \mid j \in [1, k-1], \; j \ne j_0, \; s_{k-1}(j) \ge k \} \sqcup \{ y_{j_0} x_k - c_k \}
\end{equation}
is a list of the homogeneous Poisson-prime elements of $R_k$ up to scalar multiples. In this case, we extend $\eta$ to a function $[1,k] \rightarrow \Zset$ by setting $\eta(k) = \eta(j_0)$. Since $s_{k-1}(j_0) \ge k$, we see that $p(k) = j_0$, and so $y_k := y_{j_0} x_k - c_k$. It is easily checked that the set \eqref{P4.8} equals \eqref{yset1k}.

In the second case, $\{ y_j \mid j \in [1, k-1], \; s_{k-1}(j) \ge k \} \sqcup \{ x_k \}$ is a list of the homogeneous Poisson-prime elements of $R_k$ up to scalar multiples. In this case, we set $\eta(k)$ equal to any integer not in $\eta([1,k-1])$ and readily check the desired properties.

It remains to verify the final statement of the theorem. If $k \in [1,N]$ and $p(k) = -\infty$, then $x_k = y_k$ is a Poisson-prime element of $R_k$. Consequently, $\sig_k(b) x_k + \de_k(b) = \{x_k,b\}$ is divisible by $x_k$ for all $b \in R_{k-1}$, which forces $\de_k = 0$. On the other hand, if we have $\de_k = 0$, then $x_k$ is a homogeneous Poisson-prime element of $R_k$, whence $y_k$ must be a scalar multiple of $x_k$. That is not possible when $p(k) \ne -\infty$, by \eqref{ykdef}, so we must have $p(k) = -\infty$. This concludes the proof.
\end{proof}

To illustrate \thref{mainPprime}, we continue \exref{OMmn}.

\bex{OMmn2}
Let $R = \OMmn$ be the P-CGL extension in \exref{OMmn}. For any two subsets $I \subseteq [1,m]$ and $J \subseteq [1,n]$ of the same cardinality $d$, let $\De_{I,J}$ denote the $d \times d$ minor in $R$ with row index set $I$ and column index set $J$. The sequence of homogeneous Poisson-prime elements of $R$ from \thref{mainPprime} consists of solid minors; more precisely,
$$
y_{(r-1)n+c} = \De_{ [r- \min(r,c) +1, r] , [c- \min(r,c) +1 , c] } , \; \; \forall r \in [1,m], \; c \in [1,n] .
$$
Furthermore, the function $\eta : [1,N] \rightarrow \Zset$ can be chosen as
$$
\eta( (r-1)n+c ) := c- r , \; \; \forall r \in [1,m], \; c \in [1,n] .
$$
It is easily checked that for the P-CGL extension presentation of $R$ in \exref{OMmn}, we have $\de_k = 0$ if and only if $k \in [1,n]$ or $k = (r-1)n+1$ for some $r \in [2,m]$. Hence, $\rank R = m+n-1$, by \eqref{rankRk}.
\eex

The next result provides a constructive method for finding the sequence of Poisson-prime elements $y_1, \dots, y_N$ for a given Poisson-CGL extension $R$. Note that the elements $y'_k$ in the statement of the proposition are not a priori assumed to be prime.

\bpr{constr-y}
Let $R$ be a Poisson-CGL extension of length $N$. Assume that $y'_1, \dots, y'_N$ and $c'_1, \dots, c'_N$ are two sequences of elements of $R$ such that
\begin{enumerate}
\item[(i)] $y'_1, \dots, y'_N$ are homogeneous Poisson-normal elements of $R_1, \dots, R_N$, respectively.
\item[(ii)] $c'_k \in R_{k-1}$, for all $k \in [1,N]$.
\item[(iii)] For every $k \in [1,N]$, either $y'_k = x_k - c'_k$ or there exists $j \in [1, k-1]$ such that $y'_k = y'_j x_k - c'_k$.
\item[(iv)] If $k \in [1,N]$ and $p(k) = -\infty$, then the first equality in {\rm (iii)} holds.
\end{enumerate}
Then $y'_1, \dots, y'_N$ is precisely the sequence $y_1, \dots, y_N$ of homogeneous Poisson-prime elements  from Theorem {\rm\ref{tmainPprime}}, and the function $p$ satisfies $p(k) = j$ if the second equality in {\rm (iii)} holds and $p(k) = -\infty$ otherwise.
\epr

\begin{proof} The given assumptions imply that $c'_1 \in \KK$ and $y'_1 = x_1 - c'_1$. Since $y'_1$ is homogeneous, we must have $y'_1 = x_1 = y_1$.

Now let $k \in [2,N]$. We will prove that if $y'_i = y_i$ for all $i \in [1, k-1]$, then $y'_k = y_k$. This implies the first statement of the proposition by induction. By \prref{factorHPUFD} and Theorems \ref{tPCGL-UFD}, \ref{tmainPprime}, 
$$
y'_k = \xi \prod \{ y_i^{m_i} \mid i \in [1,k], \; s(i) > k \}
$$
for some $\xi \in \kx$ and $m_i \in \Znn$. Comparing the coefficients of $x_k$ and using the form of $y_k$ from \thref{mainPprime}, we obtain that $m_k = 1$. 

The $x_k$-coefficient of $y'_k$ is either $1$ or $y'_j$, hence either a unit or a prime element, by induction. Thus, one of the following three situations holds:
\begin{enumerate}
\item[(a)] $y_k = x_k$ and $y'_k = \xi y_k$.
\item[(b)] $y_k = x_k$ and $y'_k = \xi y_j y_k$, for some $j \in [1, k-1]$.
\item[(c)] $y_k = y_{p(k)} x_k - c_k$ and $y'_k = \xi y_k$.
\end{enumerate}
Because of assumption (iv), the situation (b) cannot occur.
In the situation (a), we have $\xi = 1$, $c'_k = 0$, and $y'_k = y_k$. In the situation (c), we obtain $\xi = 1$, $j = p(k)$, $c'_k = c_k$, and $y'_k = y_k$. This argument also proves the second statement of the proposition.
\end{proof}

We next pin down the scalars involved in certain leading coefficients and Poisson brackets.
Assume that $R$ is a Poisson-CGL extension of length $N$, with elements $y_1, \dots, y_N$ as in \thref{mainPprime}.

The algebra $R$ has the $\KK$-basis
\begin{equation}
\label{Rbasis}
\{ x^f := (x_1,\dots,x_N)^f \mid f  \in \Znn^N \}
\end{equation}
(recall \eqref{x^f}).
Denote by $\prec$ the reverse lexicographic order on $\Znn^N$; namely, 
\begin{multline*}
(m'_1, \dots, m'_N)^T \prec (m_1, \dots, m_N)^T \; \iff \; \exists\, j \in [1,N] \; \text{such that}  \\
m'_j < m_j \; \text{and} \; m'_k = m_k \; \forall\, k \in [j+1, N].
\end{multline*}
We will say that $b \in R \setminus \{0\}$ has \emph{leading term} $\lt(b) := \xi x^f$, where $\xi \in \kx$ and $f \in \Znn^N$, if
$$
b = \xi x^f + \sum_{g \in \Znn^N, \, g \prec f} \xi_g x^g
$$
for some $\xi_g \in \KK$. It follows from \eqref{sigkxj} that
\begin{equation}
\label{ltbracket}
\begin{aligned}
\lt( \{x^f, x^{f'} \} ) &= \biggl( \, \sum_{k > j} m_k m'_j \la_{kj} \biggr) x^{f+f'} = \Om_\lab(f,f') x^{f+f'}, 
\\
 \forall f &= (m_1,\dots, m_N)^T, \; f' = (m'_1,\dots, m' _N)^T \in \Znn^N ,
 \end{aligned}
 \end{equation}
 recall \eqref{Om}.
 Denote 
\begin{equation}
\label{ekbar}
 \ol{e}_k := \sum_{l = 0}^{O_-(k)} e_{p^l(k)}, \; \; \forall k \in [1,N]. 
 \end{equation}
 Equation \eqref{ykdef} implies
\begin{equation}
\label{ltyk}
\lt(y_k) = x^{\ol{e}_k}, \; \; \forall k \in [1,N].
\end{equation}

For $k,j \in [1,N]$, set
\begin{align}
\al_{kj} &:= \Om_\lab(e_k, \ol{e}_j) = \sum_{m=0}^{O_-(j)} \la_{k, p^m(j)} \in \KK,  \label{alkj}  \\
q_{kj} &:= \Om_\lab(\ol{e}_k, \ol{e}_j) = \sum_{l=0}^{O_-(k)} \sum_{m=0}^{O_-(j)} \la_{p^l(k),  p^m(j)} \in \KK,  \label{qkj}  \\
\al_{k,-\infty} &= q_{k,-\infty} := 0.  \label{alqkinf}
\end{align} 
It will also be useful to set
\begin{equation}
\label{yinf}
y_{-\infty} := 1.
\end{equation}
Since $\lab$ is a skew-symmetric matrix, so is $\qb := (q_{kj}) \in M_N(\KK)$. It follows from \eqref{ltyk}, the homogeneity of $y_j$, and \eqref{hHeigenvectors} that
\begin{equation}
\label{sigkyj}
\sig_k( x^{\ol{e}_j} ) = \al_{kj} x^{\ol{e}_j} \qquad \text{and} \qquad \sig_k(y_j) = \al_{kj} y_j , \; \; \forall 1 \le j < k \le N.
\end{equation}

\bpr{PtorusR}
Let $R$ be a Poisson-CGL extension of length $N$. The elements $y_1, \dots, y_N$ from Theorem {\rm\ref{tmainPprime}} are algebraically independent over $\KK$, and
\begin{equation}
\label{Rsqueeze}
\KK[y_1, \dots, y_N] \subseteq R \subseteq \KK[y_1^{\pm1}, \dots, y_N^{\pm1}] \subset \KK(x_1, \dots, x_N) = \KK(y_1, \dots, y_N).
\end{equation}
The algebras $\KK[y_1, \dots, y_N]$ and $\KK[y_1^{\pm1}, \dots, y_N^{\pm1}]$ are Poisson algebras, with
\begin{equation}
 \label{bracketykyj} 
\{ y_k, y_j \} = q_{kj} y_j y_k, \; \; \forall k,j \in [1,N] ,
\end{equation}
and the inclusions in \eqref{Rsqueeze} are inclusions of Poisson algebras.
\epr

\begin{proof} Algebraic independence of $y_1, \dots, y_N$ is clear from the fact that each $y_k$ is a polynomial in $\KK[x_1, \dots, x_k]$ of degree $1$ in $x_k$. For $k \in [1,N]$, either $x_k = y_k$ or $x_k = y_{p(k)}^{-1} (y_k + c_k)$ with $c_k \in R_{k-1}$. By induction, we obtain $R_k \subseteq \KK[y_1^{\pm1}, \dots, y_N^{\pm1}]$ for all $k \in [1,N]$, which yields  the second inclusion of \eqref{Rsqueeze}. The first and third inclusions are clear, and the final equality follows.

From \thref{mainPprime} and equations \eqref{ltbracket}, \eqref{ltyk}, we obtain
$$
\lt( \{ y_k, y_j \}) = q_{kj} \lt(y_j y_k), \; \; \forall j,k \in [1,N].
$$
On the other hand, using \thref{mainPprime} and \prref{info-i,ii,iii} (or \prref{ykfacts} below), we see by induction on $k$ that $\{y_k, y_j\}$ is a scalar multiple of $y_j y_k$ for $1 \le j < k \le N$. Therefore, \eqref{bracketykyj} holds for $1 \le j < k \le N$. The remaining cases follow because of the skew-symmetry of the matrix $(q_{kj})$.
\end{proof}

Applying \prref{info-i,ii,iii} to the situation of \thref{mainPprime} leads to the following facts.

\bpr{ykfacts}
Keep the assumptions and notation of Theorem {\rm\ref{tmainPprime}}. Let $k \in [2,N]$.

{\rm (a)} If $p(k) = -\infty$, then $\de_k = 0$, and the derivations of $R_k$ corresponding to its homogeneous Poisson-prime elements are determined by
$$
\{ y_j, x_k \} = - \al_{kj} y_j x_k, \; \; \forall j \in [1, k-1], \qquad\qquad \{ y_k, a \} = \sig_k(a) y_k, \; \; \forall a \in R_{k-1} ,
$$
together with the actions of $\partial_{y_j}$ on $R_{k-1}$ for $j \in [1, k-1]$ such that $s(j) \ge k$ {\rm(}obtainable by recursion{\rm)}.

{\rm (b)} If $p(k) \ne -\infty$, then the derivation $\de_k$ is nonzero and is given by
$$
\de_k(a) = \{ y_{p(k)}^{-1} c_k , a \} - \sig_k(a) y_{p(k)}^{-1} c_k , \; \; \forall a \in R_{k-1} .
$$
This derivation also satisfies
\begin{equation}
\label{delkypk}
\begin{aligned}
\de_k(y_{p(k)}) &= \la_k c_k \ne 0,  &\de_k(c_k) &= 0,  \\
\de_k(y_j) &= 0, \; \; \forall j \in [1, k-1] \; \; \text{such that} \; \; s(j) > k.
\end{aligned}
\end{equation}
The derivations of $R_k$ corresponding to  its homogeneous Poisson-prime elements are determined by
\begin{equation}
\label{brackyx}
\begin{aligned}
\{ y_j, x_k \} &= - \al_{kj} y_j x_k, \; \; \forall  j \in [1, k-1] \; \; \text{such that} \; \; s(j) > k,  \\
\{ y_k, x_k \} &= - \al_{k p(k)} y_k x_k ,  \\
\{ y_k, a \} &= ( \partial_{y_{p(k)}} + \sig_k )(a) y_k, \; \; \forall a \in R_{k-1} ,
\end{aligned}
\end{equation}
together with the actions of $\partial_{y_j}$ on $R_{k-1}$ for $j \in [1, k-1]$ such that $s(j) \ge k$ {\rm(}obtainable by recursion{\rm)}. Furthermore, the components $y_{p(k)}$ and $c_k$ of $y_k$ satisfy
$$
\{ y_{p(k)}, c_k \} = - (\al_{k p(k)} + \la_k) c_k y_{p(k)} .
$$
\epr

\bco{yjPcommxk}
If $R$ is a Poisson-CGL extension of length $N$, then every homogeneous Poisson-prime element of $R$ quasi-Poisson-commutes with $x_1, \dots, x_N$. More precisely,
\begin{equation}
\label{bracketyjxk}
\{ y_j, x_k \} = - \al_{kj} y_j x_k, \; \; \forall j,k \in [1,N] \; \; \text{with} \; \; s(j) > k.
\end{equation}
\eco

\begin{proof} We proceed by induction on $l \in [1,N]$, to prove that \eqref{bracketyjxk} holds for $j,k \in [1,l]$. The case $l = 1$ is clear, since $y_1 = x_1$ and $\al_{11} = \la_{11} = 0$.

Now let $l > 1$, and assume \eqref{bracketyjxk} holds for $j,k \in [1, l-1]$. If $j \in [1, l-1]$ and $s(j) > l$, then both cases of \prref{ykfacts} yield $\{ y_j, x_l \} = - \al_{lj} y_j x_l$. Hence, it just remains to consider $\{ y_l, - \}$.

If $p(l)  = -\infty$, then $y_l = x_l$ and $\de_l = 0$. In this case,
$$
\{ y_l, x_k \} = \sig_l(x_k) y_l = \la _{lk} y_l x_k = - \al_{kl} y_l x_k
$$
for $k \in [1, l-1]$, while $\{ y_l, x_l \} = - \al_{ll} y_l x_l$ because $\al_{ll} = \la_{ll} = 0$. Finally, suppose that $p(l) \ne -\infty$, and note that $\al_{ll} = \la_{ll} + \al_{l, p(l)} = \al_{l, p(l)}$. Hence, it follows from \prref{ykfacts}(b) that $\{ y_l, x_l \} = - \al_{ll} y_l x_l$. Since $s(p(l)) = l$, our induction hypothesis implies that $\{ y_{p(l)}, x_k \} = - \al_{k, p(l)} y_{p(l)} x_k$ for $k \in [1, l-1]$. Appealing again to \prref{ykfacts}(b), we conclude that
$$
\{ y_l, x_k \} = ( \partial_{y_{p(l)}} + \sig_l )(x_k) y_l = ( - \al_{k, p(l)} + \la_{lk}) y_l x_k = - \al_{kl} y_l x_k
$$
for $k \in [1, l-1]$.

This completes the induction.
\end{proof}

From the first equation in \eqref{delkypk}, we immediately obtain the following.

\bco{ykrecursive} Keep the assumptions and notation of Theorem {\rm\ref{tmainPprime}}. The homogeneous Poisson-prime elements $y_k$ are explicitly given by
\begin{equation}  \label{ykformula}
y_k = \begin{cases}
x_k  &\text{if} \; p(k) = -\infty  \\
y_{p(k)} x_k - \la_k^{-1} \de_k(y_{p(k)})  &\text{if} \; p(k) \ne -\infty
\end{cases}
\end{equation}
for $k \in [1,N]$.
\qed\eco

\subsection{Proof of Theorem \ref{tPCGL-list}}
\label{proofT3.4}

For any positive integer $L$ and statement $X$ about Poisson-CGL extensions, we denote by $X[L]$ the validity of the statement $X$ for all Poisson-CGL extensions $R$ of length at most $L$. The work of subsection \ref{PprimePCGL} shows that
\begin{equation}
\label{implic1}
\text{\thref{PCGL-list}(a)}[N] \; \implies \text{\thref{mainPprime}}[N] \; \text{and \prref{PtorusR}}[N].
\end{equation}
The proof of \thref{PCGL-list} will be completed once the following implications are established:
\begin{align}
\text{\thref{PCGL-list}(b)}[N-1] &\; \implies \; \text{\thref{PCGL-list}(a)}[N],  \label{implic2}  \\
 \text{\thref{PCGL-list}(a)}[N] &\; \implies \text{\thref{PCGL-list}(b)}[N].  \label{implic3}
 \end{align}
 
\begin{proof}[Proof of the implication \eqref{implic2}] We apply \thref{BlistR} to the Poisson-Ore extension $R = R_{N-1}[ x_N; \sig_N, \de_N]$. All we need to show is that in this setting, the situation (iii) of \thref{BlistR} can never occur. Suppose that situation (iii) does obtain. Then by \prref{info-i,ii,iii}, $\de_N \ne 0$ but $\de_N(u) = 0$ for all homogeneous Poisson-prime elements $u$ of $R_{N-1}$. However, this contradicts \thref{PCGL-list}(b)$[N-1]$.
\end{proof}

Our proof of \eqref{implic3} involves some analysis of derivations on the Laurent polynomial ring $\TT := \KK[y_1^{\pm1}, \dots, y_N^{\pm1}]$. For $f = (m_1, \dots, m_N)^T \in \Zset^N$, define the Laurent monomial
\begin{equation}
\label{defyf}
y^f := y_1^{m_1} \cdots y_N^{m_N} \in \TT.
\end{equation}
The algebra $\TT$ is $\Zset^N$-graded, with
$$
\deg y^f := f, \; \; \forall f \in \Zset^N .
$$
We will say that a $\KK$-linear map $\varphi : \TT \rightarrow \TT$ is \emph{$\Zset^N$-homogeneous of degree $g \in \Zset^N$} if $\varphi(y^f) \in \KK y^{f+g}$ for all $f \in \Zset^N$, where we use the term $\Zset^N$-homogeneous to distinguish the above condition from homogeneity with respect to the $\xh$-grading. Given a general $\KK$-linear map $\varphi$ on $\TT$, there are unique $\Zset^N$-homogeneous $\KK$-linear maps $\varphi^g$ of degree $g$ on $\TT$ such that
$$
\varphi = \sum_{g \in \Zset^N} \varphi^g .
$$
If $\de$ is a derivation of $\TT$, then each component $\de^g$, for $g \in \Zset^N$, is a derivation of $\TT$. Since $\TT$ is a finitely generated algebra, $\de^g \ne 0$ for at most finitely many $g \in \Zset^N$.

Let $\prec$ be the reverse lexicographic order on $\Zset^N$ (defined as it was above on $\Znn^N$). Any nonzero element $u \in \TT$ can be uniquely written in the form
\begin{equation}
\label{uinT}
u = \zeta_1 y^{f_1} + \cdots + \zeta_r y^{f_r} \; \; \text{where} \; f_1 \prec \cdots \prec f_r \; \text{in} \; \Zset^N \; \text{and} \; \zeta_1,  \dots, \zeta_r \in \kx.
\end{equation}
We will say that $\zeta_r y^{f_r}$ is the \emph{leading term} of $u$ and denote it $\lt_\TT(u)$, to distinguish it from our previous usage of leading terms. For future reference, observe that
\begin{equation}
\label{lt(ayk+b)}
\lt_\TT( a y_k + b ) = \lt_\TT(a) y_k, \; \; \forall k \in [2,N], \; \; a,b \in \KK[y_1^{\pm1}, \dots, y_{k-1}^{\pm1}], \; \; a \ne 0.
\end{equation}

If $\de$ is a nonzero derivation of $\TT$, we have
\begin{equation}
\label{comp.de}
\de = \de^{g_1} + \cdots + \de^{g_t} \; \; \text{with} \; \; g_1 \prec \cdots \prec g_t \; \text{in} \; \Zset^N, \; \; \text{all} \; \de^{g_i} \ne 0.
\end{equation}
For $m > 0$ and $u$ as in \eqref{uinT}, the component of $\de^m(u)$ in degree $f_r + mg_t$ is $\zeta_r (\de^{g_t})^m (y^{f_r})$. Hence,
\begin{equation}
\label{demu}
\de^m(u) = 0 \; \implies \; (\de^{g_t})^m (\lt_\TT(u)) = 0
\end{equation}
for all $m > 0$ and nonzero $u \in \TT$.

\ble{deyyinv}
Let $\TT$ be a domain of characteristic zero and $\de$ a derivation on $\TT$. Suppose $y \in \TT$ is a unit such that $\de^m(y) = \de^m(y^{-1}) = 0$ for some $m > 0$. Then $\de(y) = 0$.
\ele

\begin{proof} Let $r, s \ge 0$ be maximal such that $\de^r(y) \ne 0$ and $\de^s(y^{-1}) \ne 0$. Then
$$
\de^{r+s}(yy^{-1}) = \sum_{l=0}^{r+s} \tbinom{r+s}{l} \de^l(y) \de^{r+s-l}(y^{-1}) = \tbinom{r+s}{r} \de^r(y) \de^s(y^{-1}) \ne 0.
$$
Since $\de(1) = 0$, we must have $r + s = 0$. Hence, $r = 0$, and therefore $\de(y) = 0$.
\end{proof}

\begin{proof}[Proof of the implication \eqref{implic3}] Let $R$ be a Poisson-CGL extension of length $N$. By \eqref{implic1}, the statements of \thref{mainPprime} and \prref{PtorusR} hold for $R$. Let $\de$ be a locally nilpotent derivation of $R$ satisfying \eqref{de0Pprime}, and suppose that $\de \ne 0$. The assumption \eqref{de0Pprime} implies
\begin{equation}
\label{deyk0}
\de(y_k) = 0, \; \; \forall k \in [1,N] \; \; \text{such that} \; \; s(k) = +\infty.
\end{equation}

Continue to denote by $\de$ the canonical (unique) extension of $\de$ to a derivation of the quotient field $\KK(x_1, \dots, x_N)$.
Set  $\TT := \KK[y_1^{\pm1}, \dots, y_N^{\pm1}]$ as above. Since $\de(y_k) \in R \subseteq \TT$ for all $k \in [1,N]$, we see that $\de(\TT) \subseteq \TT$. Now view $\de$ as a derivation of $\TT$, and decompose $\de$ as in \eqref{comp.de}. Since $\de^{g_t} \ne 0$, we must have $\de^{g_t}(y_j) \ne 0$ for some $j \in [1,N]$. Moreover, $\de(y_j) \ne 0$, and so $s(j) \ne +\infty$.
We will prove the following implication:
\begin{equation}
\label{implic4}
\text{If} \; \de^{g_t}(y_k) = 0 \; \text{for some} \; k \in [1,N] \; \text{with} \; p(k) \ne -\infty, \; \text{then} \; \de^{g_t}(y_{p(k)}) = 0.
\end{equation}
By \eqref{deyk0} and \eqref{demu}, $\de^{g_t}(y_k) = 0$ for all $k \in [1,N]$ such that $s(k) = +\infty$. A downward recursive application of \eqref{implic4} leads to $\de^{g_t}(y_k) = 0$ for all $k \in [1,N]$, contradicting the fact that $\de^{g_t}(y_j) \ne 0$. This contradiction proves the implication \eqref{implic3}.

We are left with proving \eqref{implic4}. Assume that $\de^{g_t}(y_k) = 0$ for some $k \in [1,N]$ with $p(k) \ne -\infty$. There exists $m > 0$ such that
$$
0 = \de^m(x_k) = \de^m( y_{p(k)}^{-1} y_k + y_{p(k)}^{-1} c_k ).
$$
Since $c_k \in R_{k-1}$, the elements $c_k$ and $y_{p(k)}^{-1} c_k$ belong to $\KK[y_1^{\pm1}, \dots, y_{k-1}^{\pm1}]$. Observations \eqref{lt(ayk+b)} and \eqref{demu} then imply $(\de^{g_t})^m( y_{p(k)}^{-1} y_k ) = 0$. But $\de^{g_t}(y_k) = 0$, so
\begin{equation}
\label{degtm}
0 = (\de^{g_t})^m( y_{p(k)}^{-1} y_k ) = \bigl( (\de^{g_t})^m( y_{p(k)}^{-1} ) \bigr) y_k, \; \; \text{i.e.,} \; \; (\de^{g_t})^m( y_{p(k)}^{-1} ) = 0.
\end{equation}
On the other hand, the restriction of $\de$ to $R$ is locally nilpotent. Hence, $\de^{m'}(y_{p(k)}) = 0$ for some $m' > 0$. It follows from \eqref{demu} that $(\de^{g_t})^{m'}(y_{p(k)}) = 0$. We combine this and the last equality in \eqref{degtm}, and apply \leref{deyyinv} to obtain $\de^{g_t}(y_{p(k)}) = 0$, as required.
\end{proof}

\sectionnew{Symmetry and maximal tori for Poisson-CGL extensions}
\label{maxtorisymm}

We first introduce \emph{symmetric} Poisson-CGL extensions, which are ones that also have a Poisson-CGL extension presentation with the variables in reverse order.
We then analyze the torus action on an arbitrary Poisson-CGL extension $R$, and show that there is a unique torus which is maximal among tori acting faithfully on $R$ by rational Poisson actions.

\subsection{Symmetric Poisson-CGL extensions}
\label{symmP-CGL}

If $R = \KK[x_1, \dots, x_N]$ is a polynomial ring, with the variables fixed in the order $x_1, \dots, x_N$, we denote
\begin{equation}
\label{R[jk]}
R_{[j,k]} := \KK[ x_i \mid j \le i \le k], \; \; \forall j, k \in [1,N].
\end{equation}
So, $R_{[j,k]} = \KK$ if $j \nleq k$.

\bde{symmPCGL}
We call a Poisson-CGL extension $R$ of length $N$ as in \deref{PCGL} \emph{symmetric} if the following conditions hold:
\begin{enumerate}
\item[(i)] $\de_k(x_j) \in R_{[j+1, k-1]}$, for all $1 \le j < k \le N$.
\item[(ii)] For all $j \in [1,N]$, there exists $h^*_j \in \h$ such that
$$
h^*_j \cdot x_k = - \la_{kj} x_k = \la_{jk} x_k, \; \; \forall k \in [j+1, N]
$$
and $h^*_j \cdot x_j = \la^*_j x_j$ for some nonzero $\la^*_j \in \KK$.
\end{enumerate}

Under these conditions, set
$$
\sig^*_j := (h^*_j \cdot), \; \; \forall j \in [1, N-1].
$$
Then $\sig^*_j$ is a Poisson derivation on $R$, and there is an inner Poisson $\sig^*_j$-derivation $\de^*_j$ on $R$ given by $\de^*_j(a) := \{ x_j, a \} - \sig^*_j(a) x_j$. The restriction of $\de^*_j$ to $R_{[j+1, N]}$ is determined by
$$
\de^*_j(x_k) = \{ x_j, x_k \} - \la_{jk} x_j x_k = - \de_k(x_j), \; \; \forall k \in [j+1, N].
$$
For all $1 \le j < k \le N$, the maps $\sig_k$ and $\de_k$ preserve $R_{[j,k-1]}$, and $\sig^*_j$ and $\de^*_j$ preserve $R_{[j+1, N]}$. This gives rise to the Poisson-Ore extensions
\begin{equation}
\label{Rjkpres}
R_{[j,k]} = R_{[j,k-1]} [x_k; \sig_k, \de_k]_p \qquad \text{and} \qquad R_{[j,k]} = R_{[j+1, k]} [x_j; \sig^*_j, \de^*_j]_p .
\end{equation}
In particular, it follows that $R$ has an iterated Poisson-Ore extension presentation with the variables $x_k$ in descending order:
\begin{equation}
\label{revitpOre}
R = \KK[x_N]_p [x_{N-1}; \sig^*_{N-1}, \de^*_{N-1}]_p \cdots [x_1; \sig^*_1, \de^*_1]_p ,
\end{equation}
which is the reason for the name ``symmetric". The presentation \eqref{revitpOre} is, in fact, a Poisson-CGL extension presentation of $R$; we will record this as part of \prref{tauPCGL}.
\ede

\bex{OMmn4}
For example, the Poisson-CGL extension $\OMmn$ of \exref{OMmn} is symmetric. It is clear that in this example, condition (i) of \deref{symmPCGL} holds. Condition (ii) can be verified for the following elements $h^*_j \in \h$. Given $j \in [1,N]$, write $j = (r-1)n+c$ for some $r\in [1,m]$ and $c\in [1,n]$, and take
$$
h^*_j := (0, \dots, 0, 1, 0, \dots, 0, -1, 0, \dots, 0)
$$
where the entries $1$ and $-1$ reside in positions $r$ and $m+c$, respectively. In particular, $h^*_j \cdot x_j = 2x_j$, so that $\la^*_j = 2$.
\eex

\bde{XiN}
Denote the following subset of the symmetric group $S_N$:
\begin{equation}
\label{tau}
\begin{aligned}
\Xi_N := \{ \tau \in S_N \mid 
\tau(k) &= \max \, \tau( [1,k-1]) +1 \; \;
\mbox{or}  \\
\tau(k) &= \min \, \tau( [1,k-1]) - 1, 
\; \; \forall k \in [2,N] \}.
\end{aligned}
\end{equation}
In other words, $\Xi_N$ consists of those $\tau \in S_N$ 
such that $\tau([1,k])$ is an interval for all $k \in [2,N]$. 
For each $\tau \in \Xi_N$, we have the iterated 
Poisson-Ore extension presentation
\begin{equation}
\label{taupOre}
R = \KK [x_{\tau(1)}]_p [x_{\tau(2)}; \sig''_{\tau(2)}, \de''_{\tau(2)}]_p 
\cdots [x_{\tau(N)}; \sig''_{\tau(N)}, \de''_{\tau(N)}]_p ,
\end{equation}
where $\sig''_{\tau(k)} := \sig_{\tau(k)}$ and 
$\de''_{\tau(k)} := \de_{\tau(k)}$ if 
$\tau(k) = \max \, \tau( [1,k-1]) +1$, while 
$\sig''_{\tau(k)} := \sig^\sy_{\tau(k)}$ and 
$\de''_{\tau(k)} := \de^\sy_{\tau(k)}$ if 
$\tau(k) = \min \, \tau( [1,k-1]) -1$.
\ede

The following proposition is straightforward.

\bpr{tauPCGL} 
For every symmetric P-CGL extension $R$ of length $N$ and $\tau \in \Xi_N$,
the iterated Poisson-Ore extension presentation \eqref{taupOre} of $R$ 
is a Poisson-CGL extension presentation for the 
same choice of torus $\HH$, and the associated elements 
$h''_{\tau(1)}, \ldots, h''_{\tau(N)} \in \h$ required by Definition {\rm\ref{dPCGL}(iii)} 
are given by $h''_{\tau(k)} = h_{\tau(k)}$ if $\tau(k) = \max \, \tau( [1,k-1]) +1$  
and $h''_{\tau(k)} = h^*_{\tau(k)}$ if $\tau(k) = \min \, \tau( [1,k-1]) -1$.
\epr

\thref{Hmax=G} implies that the group $\Hmax(R)$ does not change in passing from the Poisson-CGL extension presentation \eqref{itpOre} to \eqref{taupOre}, modulo identifying $(\kx)^N$ with the copy obtained by permuting coordinates with $\tau$. Consequently,  the theorem shows that the rank of $R$ does
not depend on the choice of Poisson-CGL extension presentation \eqref{taupOre}.

When describing permutations $\tau \in S_N$ as functions, we will use the one-line notation,
\begin{equation}
\label{one-line}
\tau = [\tau(1),\tau(2),\dots,\tau(N)] := \begin{bmatrix} 1 &2 &\cdots &N\\ \tau(1) &\tau(2) &\cdots &\tau(N) \end{bmatrix} .
\end{equation}
A special role is played by the longest element of $S_N$,
\begin{equation}
\label{tau-ci}
w_\circ := [N, N-1, \ldots, 1].
\end{equation}
The corresponding Poisson-CGL extension presentation from \prref{tauPCGL}, namely \eqref{revitpOre}, is symmetric, 
while the ones for the other elements of $\Xi_N$ do not possess 
this property in general.

\subsection{Maximal tori}
\label{maxtori}

Let $R$ be a Poisson-CGL extension of length $N$ as in \eqref{itpOre}. Equip $R$ with the rational action of the torus $(\kx)^N$ by $\KK$-algebra automorphisms such that
\begin{equation}
\label{K*Naction}
(\al_1, \dots, \al_N) \cdot x_i =  \al_i x_i , \; \; \forall i \in [1,N].
\end{equation}
The differential of this action is the derivation action of $\KK^N = \Lie \, (\kx)^N$ on $R$ given by
\begin{equation}
\label{KNaction}
(\be_1, \dots, \be_N) \cdot (x_1^{m_1} \cdots x_N^{m_N}) =  \biggl( \, \sum_{i=1}^N m_i \be_i \biggr) x_1^{m_1} \cdots x_N^{m_N} , \; \; \forall (m_1, \dots, m_N) \in \Znn^N .
\end{equation}

The given action of $\HH$ on $R$ factors through the above $(\kx)^N$-action via a morphism (of algebraic groups) $\HH \rightarrow (\kx)^N$. Since nothing is lost by reducing $\HH$ modulo the kernel of its action, we may assume that the action of $\HH$ is faithful and then identify $\HH$ with its image in $(\kx)^N$. Thus,
\begin{equation}
\label{WLOGH}
\text{Without loss of generality,} \; \HH \; \text{is a closed subgroup of} \; (\kx)^N.
\end{equation}
For the fact that the image of $\HH$ in $(\kx)^N$ must be closed, see, e.g., \cite[Corollary 1.4]{Bor}. 

Next, set
\begin{equation}
\label{groupG}
G := \{ \psi \in (\kx)^N \mid (\psi \cdot), \; \; \text{acting as in \eqref{K*Naction}, is a Poisson automorphism of} \; R \},
\end{equation}
and observe that $G$ is a closed subgroup of $(\kx)^N$. Since $G$ is diagonalizable, its connected component of the identity, $G^\ci$, is a torus (e.g., \cite[Corollary 8.5]{Bor}). This subgroup is the unique maximal torus of $G$, and so it contains $\HH$. Let us set
\begin{equation}
\label{defHmax}
\Hmax(R) := G^\ci.
\end{equation}
(The definition of this torus, and its position within $(\kx)^N$, depend on the Poisson-CGL extension presentation \eqref{itpOre} of $R$. However, we do not indicate this dependence in the notation.) Since $\Hmax(R)$ contains $\HH$, the algebra $R$ is also a Poisson-CGL extension with respect to $\Hmax(R)$. We shall prove that $G$ is connected and so $\Hmax(R) = G$.

The torus $\Hmax(R)$ is unique (up to isomorphism) because of the following universal property, assuming that we fix the Poisson-CGL extension presentation \eqref{itpOre}. If $\HH_1$ is any torus with a rational Poisson action on $R$ such that $(R,\HH_1)$ is Poisson-CGL for the presentation \eqref{itpOre}, then the action of $\HH_1$ on $R$ factors uniquely through the action of $\Hmax(R)$, via an algebraic group morphism $\HH_1 \rightarrow \Hmax(R)$.

Analogously to $G$ and $G^\ci$, define the Lie subalgebra
\begin{equation}
\label{hmax}
\hmax(R) := \{ t \in \KK^N \mid (t \cdot), \; \; \text{acting as in \eqref{KNaction}, is a Poisson derivation of} \; R \},
\end{equation}
and observe that $\h = \Lie \HH \subseteq \hmax(R)$. 
We shall prove that $G$ is connected and $\Lie G = \hmax(R)$. To do so, we first pin down $\hmax(R)$.

\ble{tinhmax}
Let $R$ be a P-CGL extension of length $N$ as in \eqref{itpOre}. For any $j \in [1,N]$ with $\de_j = 0$, there exists $t \in \hmax(R)$ such that $t\cdot x_j = x_j$ and $t\cdot x_k = 0$ for all $k \in [1,j-1]$.
\ele

\begin{proof} By \thref{mainPprime}, $p(j) = -\infty$. There is some $l \in [1,N]$ with $s(l)= +\infty$ and $p^{O_-(l)}(l) = j$. Then $y_l$ is Poisson-normal in $R$, and we have a Poisson derivation $\partial_{y_l}$ on $R$ (recall \eqref{defdelc}) such that $\{ y_l, a \} = \partial_{y_l}(a) y_l$ for all $a\in R$. In view of \coref{yjPcommxk}, $\partial_{y_l}(x_k) = - \al_{kl} x_k$ for all $k \in [1,N]$.

Set $\theta := - \partial_{y_l} + \sum_{m=0}^{O_-(l)} (h_{p^m(l)} \cdot)$. Since $\theta$ is a Poisson derivation of $R$ for which all the $x_k$ are eigenvectors, $\theta = (t' \cdot)$ for some $t' \in \hmax(R)$. Since $h_{p^m(l)} \cdot x_k = \la_{p^m(l),k} x_k$ when $k < p^m(l)$, we find that $t' \cdot x_k = 0$ for $k < j$ and $t' \cdot x_j = \la_j x_j$. Thus, the element $t := \la_j^{-1} t'$ of $\hmax(R)$ has the desired properties.
\end{proof}

\bpr{hmaxform}
Let $R$ be a Poisson-CGL extension of length $N$ and rank $n$ as in Definition {\rm\ref{dPCGL}}, and define $\hmax(R)$ as in \eqref{hmax}. Set 
$$
D' := \{ k \in [2,N] \mid \de_k \ne 0 \} = \{ k \in [2,N] \mid p(k) \ne -\infty \}.
$$
For each $k \in D'$, choose $j_k \in [1, k-1]$ such that $\de_k(x_{j_k}) \ne 0$, and choose $f_k \in \Znn^{k-1}$ such that the monomial $x^{f_k}$ appears in $\de_k(x_{j_k})$ when $\de_k(x_{j_k})$ is expressed in the basis \eqref{Rbasis}. Then
\begin{equation}
\label{hmaxdescr}
\hmax(R) = \bigl\{ \xi \in \KK^N \bigm| \xi_k = - \xi{j_k} + \sum_{i=1}^{k-1}  f_{ki} \xi_i, \; \; \forall k \in D' \},
\end{equation}
and $\dim \hmax(R) = n$.
\epr

\begin{proof} Let $\h_1$ denote the subspace of $\KK^N$ described on the right hand side of \eqref{hmaxdescr}, and note that $\dim \h_1 = N - |D'| = n$. For $j \in [1,N] \setminus D'$, \leref{tinhmax} provides us with an element $t_j \in \hmax(R)$ such that $t_j \cdot x_j = x_j$ and $t_j \cdot x_k = 0$ for all $k \in [1,j-1]$. These $t_j$ are obviously linearly independent, and therefore $\dim \hmax(R) \ge n$.

For any $\xi \in \hmax(R)$, the assumption that $(\xi \cdot)$ is a Poisson derivation implies
$$
\xi \cdot \{ x_k, x_{j_k} \} = \{ \xi_k x_k ,x_{j_k} \} + \{ x_k, \xi_{j_k} x_{j_k} \} = (\xi_k + \xi_{j_k}) \bigl( \la_{k,j_k} x_{j_k} x_k + \de_k(x_{j_k}) \bigr)
$$
for $k \in D'$. On the other hand,
$$
\xi \cdot \{ x_k, x_{j_k} \} = \xi \cdot \bigl( \la_{k,j_k} x_{j_k} x_k + \de_k(x_{j_k}) \bigr) =  \la_{k,j_k} (\xi_{j_k} + \xi_k)  x_{j_k} x_k + \xi \cdot \de_k(x_{j_k}),
$$
from which we see that $\xi \cdot \de_k(x_{j_k}) = (\xi_k + \xi_{j_k}) \de_k(x_{j_k})$, and hence  $\xi_k + \xi_{j_k} = \sum_{i=1}^{k-1}  f_{ki} \xi_i$. Since this holds for all $k \in D'$, we conclude that $\xi \in \h_1$.

Thus, $\hmax(R) \subseteq \h_1$. Comparison of dimensions then forces $\hmax(R) = \h_1$ and $\dim \hmax(R) = n$.
\end{proof}

In order to establish an analog of \prref{hmaxform} for $\Hmax(R)$, we will need to transport Poisson properties from $\hmax(R)$ to a torus whose Lie algebra is $\hmax(R)$. The following lemma accomplishes this.

\ble{hHPoisson}
Let $R$ be a Poisson algebra equipped with a rational action of a torus $\HH$ by $\KK$-algebra automorphisms, and let $\h = \Lie \HH$ act on $R$ by the differential of the $\HH$-action. Then $\HH$ acts on $R$ by Poisson automorphisms if and only if $\h$ acts on $R$ by Poisson derivations.
\ele

\begin{proof} We have already noted, in \leref{rat-p-action}, that if $\HH$ acts by Poisson automorphisms, then $\h$ acts by Poisson derivations \cite[Lemma 1.4]{GLaun}. Now assume that $\h$ acts by Poisson derivations.

The algebra $R$ is graded by the character group $X(\HH)$, which we will write additively in this proof. We may assume that $\HH = (\kx)^n$ and $\h = \KK^n$, where $n := \rank \HH$. Then $X(\HH)$ is naturally identified with $\Zset^n$, so that
$$
u(h) = h_1^{u_1} h_2^{u_2} \cdots h_n^{u_n}, \; \; \forall u = (u_1,\dots,u_n) \in \Zset^n, \; \; h = (h_1,\dots,h_n) \in \HH.
$$
For $u = (u_1,\dots,u_n) \in \Zset^n$, let $R_u$ denote the $\HH$-eigenspace of $R$ for the $\HH$-eigenvalue $u$, and define $u' \in \h^*$ by
$$
u'(\xi) = u_1 \xi_1 + \cdots + u_n \xi_n, \; \; \forall \xi = (\xi_1, \dots, \xi_n) \in \h.
$$
The map $u \mapsto u'$ is an injective group homomorphism $\Zset^n \rightarrow \h^*$.

Recall from \eqref{hHeigenspaces} that the $\HH$-eigenspaces of $R$ coincide with the $\h$-eigenspaces. As indicated in \cite[Lemma 1.3]{GLaun}, the action of $\h$ on the $\HH$-eigenspaces of $R$ is given by
\begin{equation}
\label{xidota}
\xi \cdot a = u'(\xi) a, \; \; \forall \xi \in \h, \; \; a \in R_u, \; \; u \in \Zset^n.
\end{equation}

Let $h \in \HH$. To prove that $(h \cdot)$ is a Poisson automorphism of $R$, it is enough to show that it preserves Poisson brackets of homogeneous elements. Thus, let $u,v \in \Zset^n$ and $a\in R_u$, $b\in R_v$. Then
$$
\{ h\cdot a, h\cdot b \} = \{ u(h) a, v(h) b \} = u(h) v(h) \{ a, b \}.
$$
If $\{a,b\} = 0$, then immediately $h \cdot \{ a,b \} = \{ h\cdot a, h\cdot b \}$, so we may assume that $\{ a,b \} \ne 0$. For any $\xi \in \h$, we have
\begin{equation}
\label{xibrackab}
\xi \cdot \{ a,b \} = \{ \xi \cdot a, b \} + \{ a, \xi \cdot b \} = (u'(\xi) + v'(\xi)) \{ a,b \}.
\end{equation}
Then $\{ a,b \}$ is an $\h$-eigenvector, hence also an $\HH$-eigenvector, so $\{ a,b \} \in R_w$ for some $w \in \Zset^n$. Moreover, it follows from \eqref{xidota} and \eqref{xibrackab} that $w' = u' + v'$, whence $w = u + v$. Therefore
$$
h \cdot \{ a,b \} = w(h) \{ a,b \} = u(h) v(h) \{ a, b \} = \{ h\cdot a, h\cdot b \},
$$
as required.
\end{proof}

\bth{Hmax=G}
Let $R$ be a Poisson-CGL extension of length $N$ and rank $n$ as in Definition {\rm\ref{dPCGL}}, and define the group $G \subseteq (\kx)^N$ as in \eqref{groupG}. Then
$$
\Hmax(R) = G, \qquad \rank \Hmax(R) = n, \qquad \text{and} \qquad \Lie \Hmax(R) = \hmax(R).
$$

Now set $D' := \{ k \in [2,N] \mid \de_k \ne 0 \} = \{ k \in [2,N] \mid p(k) \ne -\infty \}$. For each $k \in D'$, choose $j_k \in [1, k-1]$ such that $\de_k(x_{j_k}) \ne 0$, and choose $f_k \in \Znn^{k-1}$ such that the monomial $x^{f_k}$ appears in $\de_k(x_{j_k})$ when $\de_k(x_{j_k})$ is expressed in the basis \eqref{Rbasis}. Then
\begin{equation}
\label{Hmax.descript}
\Hmax(R) = \bigl\{ \psi \in (\kx)^N \bigm| \psi_k = \psi_{j_k}^{-1} \prod_{i=1}^{k-1} \psi_i^{f_{ki}}, \; \; \forall k \in D' \}.
\end{equation}
\eth

\begin{proof} Throughout the proof, assume that $\HH = \Hmax(R)$. Let $G_2$ denote the closed subgroup of $(\kx)^N$ described on the right hand side of \eqref{Hmax.descript}, and note that $G_2$ is a torus of rank $N - |D'| = n$, recalling \eqref{rankRk}. We shall prove that $\HH = G = G_2$ and $\Lie G_2 = \hmax(R)$.

By construction, $\HH \subseteq G$. Let $\psi \in G$, and let $k \in D'$. On applying the Poisson automorphism $(\psi \cdot)$ to the relation
$$
\{ x_k, x_{j_k} \} = \la_{k, j_k} x_{j_k} x_k + \de_k(x_{j_k}),
$$
we see that $\psi \cdot \de_k(x_{j_k}) = \psi_k \psi_{j_k} \de_k(x_{j_k})$. Consequently, all the monomials appearing in the expansion of $\de_k(x_{j_k})$ in the basis \eqref{Rbasis} must have $\psi$-eigenvalue $\psi_k \psi_{j_k}$. One of these monomials is $x^{f_k}$, whose $\psi$-eigenvalue also equals $\prod_{i=1}^{k-1} \psi_i^{f_{ki}}$. Hence, $\psi_k = \psi_{j_k}^{-1} \prod_{i=1}^{k-1} \psi_i^{f_{ki}}$. This proves that $\HH \subseteq G \subseteq G_2$.

Write $L := \OO((\kx)^N)$ as a Laurent polynomial ring $L = \KK[ Z_1^{\pm1}, \dots, Z_N^{\pm1}]$, where $Z_k$ is the $k$-th coordinate projection $(\kx)^N \rightarrow \KK$. The standard identification of $\KK^N$ with the Lie algebra of $(\kx)^N$ is via the map
$$
(\xi_1, \dots, \xi_N) \longmapsto \sum_{k=1}^N \xi_k \frac{\partial}{\partial Z_k} \bigg|_e ,
$$
where $e := (1,\dots,1)$. Now $\OO(G_2) = L/I_2$ where
$$
I_2 := \bigl\langle Z_k - Z_{j_k}^{-1} \prod_{i-1}^{k-1} Z_i^{f_{ki}} \bigm| k \in D' \bigr\rangle.
$$
Observe that a point derivation $\sum_{k=1}^N \xi_k \frac{\partial}{\partial Z_k} \big|_e$ vanishes on $I_2$ if and only if 
$$
\xi_k + \xi_{j_k} - \sum_{i=1}^{k-1} f_{ki} \xi_i = 0, \; \; \forall k \in D'.
$$
Thus, $\Lie G_2 = \hmax(R)$, by \eqref{hmaxdescr}.

By definition of $\hmax(R)$, we now have that $\Lie G_2$ acts on $R$ by Poisson derivations. \leref{hHPoisson} then implies that $G_2$ acts on $R$ by Poisson automorphisms. Therefore $G_2 \subseteq G$. Since $G_2$ is a torus, it follows that $\HH = G = G_2$. This verifies \eqref{Hmax.descript}, and also yields $\Lie \HH = \hmax(R)$. Finally, $\rank \HH = \dim \hmax(R) = n$ by \prref{hmaxform}.
\end{proof}

\bex{OMmn3}
In the case of $R = \OMmn$ (\exref{OMmn}), \thref{Hmax=G} can be applied as follows. First, $D' = \{ (r-1)n+c \mid r\in [2,m], \; c \in [2,n] \}$. For $k = (r-1)n+c \in D'$, we can choose $j_k := 1$ and $f_k := e_c + e_{(r-1)n+1}$. Then
$$
\Hmax(R) = \{ \psi \in (\kx)^{mn} \mid \psi_{(r-1)n+c} = \psi_1^{-1} \psi_c \psi_{(r-1)n+1} \;\; \forall r \in [2,m], \; c \in [2,n] \}.
$$
The relation of this torus to the torus $\HH = (\kx)^{m+n}$ of \exref{OMmn} is as follows. There is a surjective morphism of algebraic groups $\pi : \HH \rightarrow \Hmax(R)$ given by
$$
\pi(\xi)_{(r-1)n+c} = \xi_r \xi_{m+c}^{-1} , \; \; \forall r \in [1,m], \; c \in [1,n],
$$
and $\pi$ transports the given action of $\HH$ to that of $\Hmax(R)$, that is, $\pi(\xi) \cdot r = \xi \cdot r$ for all $\xi \in \HH$ and $r \in R$.
\eex

\sectionnew{One-step mutations in Poisson-CGL extensions}
\label{mCGL}
In this section we obtain a very general way of 
constructing mutations of potential cluster variables in Poisson-CGL extensions.
The key idea is that, if an algebra $R$ has 
two different P-CGL extension presentations
obtained by reversing the order in which two adjacent variables 
$x_k$ and $x_{k+1}$ are adjoined, then 
the corresponding sequences of Poisson-prime elements 
from \thref{mainPprime} are obtained by a type of 
mutation formula. This is realized in \S \ref{genmut}. 
In Section \ref{mut-sym}, we treat more general mutations. 
One problem arises along the way: In the mutation-type formula \eqref{thm2-disp} for the current situation, 
the last term has a nonzero coefficient which 
does not equal one in general. 
For a one-step mutation such a coefficient 
can be always made 1 after rescaling, but for the purposes of constructing  cluster 
algebras one needs to be able to synchronize those rescalings to obtain a chain of mutations.
This delicate issue is 
resolved in the next two sections.

We investigate a P-CGL extension
\begin{equation}
\label{firstPCGL}
R := \KK[x_1]_p [x_2; \sig_2, \delta_2]_p \cdots 
[x_k; \sig_k, \delta_k]_p
[x_{k+1}; \sig_{k+1}, \delta_{k+1}]_p \cdots
[x_N; \sig_N, \delta_N]_p
\end{equation}
of length $N$ as in \deref{PCGL} such that, for some $k \in [1,N-1]$, $R$ has a second P-CGL extension presentation of the form 
\begin{multline}
\label{secondPCGL}
R := \KK[x_1]_p [x_2; \sig_2, \delta_2]_p \cdots 
[x_{k-1}; \sig_{k-1}, \delta_{k-1}]_p  \\
[x_{k+1}; \sig'_k, \delta'_k]_p
[x_k; \sig'_{k+1}, \delta'_{k+1}]_p 
[x_{k+2}; \sig_{k+2}, \delta_{k+2}]_p
\cdots
[x_N; \sig_N, \delta_N].
\end{multline}
Corresponding to the presentation \eqref{secondPCGL}, we write $y'_l$ for the elements from \thref{mainPprime}, $\la'_{lj}$, $\al'_{lj}$, $q'_{lj}$ for the scalars from \eqref{sigkxj}, \eqref{alkj}, \eqref{qkj}, and $h'_l$, $\la'_l$ for the elements of $\h$ and $\kx$ from \deref{PCGL}(iii).

\subsection{A general mutation formula}
\label{genmut}

\ble{PCGLtranspose}
Assume that $R$
is a P-CGL extension as in \eqref{firstPCGL}, and that 
$R$ has a second P-CGL extension presentation of the form \eqref{secondPCGL}.
Then $\de_{k+1}$ and $\de'_{k+1}$ map $R_{k-1}$ to itself, and
\begin{equation}
\label{sig'de'}
\begin{aligned}
\sig'_k &= \sig_{k+1}|_{R_{k-1}}  &\de'_k &= \de_{k+1}|_{R_{k-1}}  \\
\sig'_{k+1}|_{R_{k-1}} &= \sig_k  &\de'_{k+1}|_{R_{k-1}} &= \de_k .
\end{aligned}
\end{equation}
Moreover, 
\begin{align}
\label{moresig'de'}
\sig'_{k+1}(x_{k+1}) &= \la_{k,k+1} x_{k+1}  &\de'_{k+1}(x_{k+1}) &= - \de_{k+1}(x_k).
\end{align}
\ele

\begin{proof} Note first that $R_{k-1}$ is stable under $\sig_{k+1}$ and $\sig'_{k+1}$.
 For $a \in R_{k-1}$, we have $\{x_{k+1} , a\} = \sig'_k(a) x_{k+1} + \de'_k(a)$ with $\sig_k'(a), \de'_k(a) \in R_{k-1}$. Comparing this relation with $\{x_{k+1} , a\} = \sig_{k+1}(a) x_{k+1} + \de_{k+1}(a)$, and using the fact that $1$, $x_{k+1}$ are linearly independent over $R_k$, we conclude that $\sig'_k(a) = \sig_{k+1}(a)$ and $\de'_k(a) = \de_{k+1}(a)$. Thus, $R_{k-1}$ is stable under $\de_{k+1}$, and the first line of \eqref{sig'de'} holds. By symmetry (since we may view \eqref{secondPCGL} as the initial P-CGL extension presentation of $R$ and \eqref{firstPCGL} as the second one), $R_{k-1}$ is stable under $\de'_{k+1}$, and the second line of \eqref{sig'de'} holds.

Now $\{x_k , x_{k+1}\} = \la'_{k+1,k} x_{k+1} x_k + \de'_{k+1}(x_{k+1})$, and so we have
$$- \la'_{k+1,k} x_k x_{k+1} - \de'_{k+1}(x_{k+1}) = \{x_{k+1} , x_k\} = \la_{k+1,k} x_k x_{k+1} + \de_{k+1}(x_k),$$
with $\de'_{k+1}(x_{k+1}) \in R'_k = \bigoplus_{l=0}^\infty R_{k-1} x_{k+1}^l$ and $\de_{k+1}(x_k) \in R_k = \bigoplus_{l=0}^\infty R_{k-1} x_k^l$. Moreover,
\begin{equation}
\label{Rk+1basis}
R_{k+1} \; \; \mbox{is a free} \; \; R_{k-1}\mbox{-module with basis} \; \; \{x_k^{l_k} x_{k+1}^{l_{k+1}} \mid l_k, l_{k+1} \in \Zset_{\geq 0 } \}.
\end{equation}
Hence, we conclude that $- \la'_{k+1,k} = \la_{k+1,k}$ and
$- \de'_{k+1}(x_{k+1}) = \de_{k+1}(x_k) \in R_{k-1}$,
from which \eqref{moresig'de'} follows.
\end{proof}

\bth{swapkk+1} Assume that $R$ is a Poisson-CGL extension of length $N$ as in {\rm\eqref{firstPCGL}},
and $k \in [1,N-1]$ such that $\de_{k+1}(x_k) \in \KK$.
Denote by $y_1, \ldots, y_N$ and $\eta : [1,N] \to \Zset$ 
the sequence and function from Theorem {\rm\ref{tmainPprime}}. Assume that 
$R$ has a second Poisson-CGL extension presentation of the form {\rm\eqref{secondPCGL}}, and let $y'_1, \ldots, y'_N$ 
be the corresponding sequence from Theorem {\rm\ref{tmainPprime}}.

{\rm(a)} If $\eta(k) \neq \eta(k+1)$, then $y'_j = y_j$ for $j \neq k, k+1$ 
and $y'_k = y_{k+1}$, $y'_{k+1} = y_k$. 

{\rm(b)} If $\eta(k) = \eta(k+1)$, then 
$$
y_k y'_k - y_{p(k)} y_{k+1}
$$ 
is a homogeneous Poisson-normal element of $R_{k-1}$, recall \eqref{yinf}.
It Poisson-normalizes the elements of $R_{k-1}$ in exactly the same way as $y_{p(k)} y_{k+1}$.
Furthermore,
\begin{equation}
\label{y-equal}
y'_j = y_j , \; \; \forall j\in [1,N], \; j\ne k.
\end{equation}

In both cases {\rm(a)} and {\rm(b)}, the function $\eta' : [1,N] \to \Zset$ from
Theorem {\rm\ref{tmainPprime}} associated to the second presentation can be chosen 
to be $\eta' = \eta (k, k+1)$, 
where $(k,k+1)$ denotes a transposition in $S_N$. In particular, 
the ranges of $\eta$ and $\eta'$ coincide 
and the rank of $R$ is the same 
for both CGL extension presentations.
\eth

\begin{remark*} The assumption $\de_{k+1}(x_k) \in \KK$ holds in case $R$ is a symmetric P-CGL extension, as part of condition (i) in \deref{symmPCGL}. In any case, it is equivalent to $\de'_{k+1}(x_{k+1}) \in \KK$ in view of \eqref{moresig'de'}.
\end{remark*}

\begin{proof} Let $R'_j$ be the $j$-th algebra in the 
chain \eqref{secondPCGL} for $j \in [0,N]$. Obviously 
$R'_j = R_j$ for $j \neq k$ and $y'_j = y_j$
for $j \in [1,k-1]$, and we may choose $\eta'(j) = \eta(j)$ for $j \in [1,k-1]$. Since $R_{k+1}$ is a free $R_{k-1}$-module with basis as in \eqref{Rk+1basis},
and $R_k$ and $R'_k$ equal the $R_{k-1}$-submodules with bases 
$\{x_k^{l_k} \mid l_k \in \Zset_{\geq 0 } \}$ 
and 
$\{x_{k+1}^{l_{k+1}} \mid l_{k+1} \in \Zset_{\geq 0 } \}$
respectively, we have
\begin{equation}
\label{Rkk'intersect}
R_k \cap R'_k = R_{k-1}.
\end{equation}
In particular, $y'_k \notin R_k$.

Denote by $L$ the number of homogeneous Poisson-prime elements of 
$R_{k+1}$ that do not belong to $R_{k-1}$, up to taking associates.

(a) The condition $\eta(k) \neq \eta(k+1)$ implies $L =2$ and thus 
$\eta'(k) \neq \eta'(k+1)$. Moreover, $y'_k$ and $y'_{k+1}$ 
are scalar multiples of either $y_k$ and $y_{k+1}$ or $y_{k+1}$ and $y_k$. 
Since $y'_k \notin R_k$, we must have $y'_k = \xi_{k+1} y_{k+1}$ and $y'_{k+1} = \xi_k y_k$
for some $\xi_k, \xi_{k+1} \in \kx$. Invoking \thref{mainPprime} 
and looking at leading terms gives $\xi_k = \xi_{k+1} = 1$ and allows us to choose $\eta'(k) = \eta(k+1)$ and $\eta'(k+1) = \eta(k)$. It follows
from $R'_j = R_j$ for $j>k+1$ that $y'_j = y_j$ for such $j$, and that we can choose $\eta'(j) = \eta(j)$ 
for all $j > k+1$. With these choices, $\eta' = \eta (k,k+1)$.

(b) In this case, $L=1$ and $y'_{k+1}$ is a scalar multiple
of $y_{k+1}$. Furthermore, $\eta'(k)$ and $\eta'(k+1)$ must agree, and they need to equal $\eta'(p(k))$ if $p(k) \ne -\infty$, so we can choose $\eta'(k) = \eta'(k+1) = \eta(k) = \eta(k+1)$. 
By \thref{mainPprime},
\begin{equation}
\label{ykyk'}
\begin{aligned}
y_k &= \begin{cases} y_{p(k)} x_k - c_k, &\quad\mbox{if} \; \; p(k) \ne -\infty\\  x_k, &\quad\mbox{if} \; \; p(k) = -\infty \end{cases}  &\qquad y_{k+1} &= y_k x_{k+1} - c_{k+1}  \\
y'_k &= \begin{cases} y_{p(k)} x_{k+1} - c'_k, &\mbox{if} \; \; p(k) \ne -\infty\\  x_{k+1}, &\mbox{if} \; \; p(k) = -\infty \end{cases}  &y'_{k+1} &= y'_k x_k - c'_{k+1}
\end{aligned}
\end{equation}
for some $c_k, c'_k \in R_{k-1}$, $c_{k+1} \in R_k$, and $c'_{k+1} \in R'_k$. Write the above elements in terms of the basis in \eqref{Rk+1basis}. The coefficients of $x_k x_{k+1}$ in $y_{k+1}$ and $y'_{k+1}$ are both equal to $y_{p(k)}$, recall \eqref{yinf}. Since $y'_{k+1}$ is a scalar multiple of $y_{k+1}$, we thus see that
\begin{equation}
\label{y'k+1}
y'_{k+1} = y_{k+1}.
\end{equation}
Eq. \eqref{y-equal} and the fact that we may choose $\eta'(j) = \eta(j)$ for $j > k+1$ now follow easily from \thref{mainPprime}. In particular, $\eta' = \eta = \eta (k, k+1)$.

Next, we verify that 
\begin{equation}
\label{in}
y_k y'_k - y_{p(k)} y_{k+1} \in R_{k-1}.
\end{equation}
Assume first that $p(k) = -\infty$.
By \eqref{delkypk}, $c_{k+1}$ is a nonzero scalar multiple of $\de_{k+1}(x_k)$, and so $c_{k+1} \in \kx$ by hypothesis.
Then from \eqref{ykyk'} we obtain
\begin{equation}
\label{eltp(k)-inf}
y_k y'_k - y_{p(k)} y_{k+1} = x_k x_{k+1} - (x_k x_{k+1} - c_{k+1}) = c_{k+1} \in \kx .
\end{equation}
 This verifies \eqref{in} (and also shows that $y_k y'_k - y_{p(k)} y_{k+1}$ is a homogeneous Poisson-normal element of $R_{k-1}$) in the present case.

Now assume that $p(k) \ne -\infty$. First, we obtain
\begin{equation}
\label{theelt}
\begin{aligned}
y_k y'_k - y_{p(k)} y_{k+1} 
 &= y_k ( y_{p(k)} x_{k+1} - c'_k) -  
y_{p(k)} (y_k x_{k+1} - c_{k+1}) \\
 &= - y_k c'_k + y_{p(k)} c_{k+1} \in R_k .
\end{aligned} 
\end{equation}
Using \eqref{y'k+1}, we obtain
\begin{align*}
y_k y'_k - y_{p(k)} y_{k+1} &= (y_{p(k)} x_k - c_k) y'_k - y_{p(k)} ( y'_k x_k - c'_{k+1} )  \\
 &= - c_k y'_k + y_{p(k)} c'_{k+1} \in R'_k .
\end{align*}
This equation, combined with \eqref{theelt} and \eqref{Rkk'intersect}, yields \eqref{in}.

We now use \eqref{theelt} to verify that $y_k y'_k - y_{p(k)} y_{k+1}$ is homogeneous. Note first that $y_k c'_k$ and $y_{p(k)} c_{k+1}$ are homogeneous. By \eqref{delkypk}, $c'_k$ and $c_{k+1}$ are scalar multiples of $\de'_k(y_{p(k)}) = \de_{k+1}(y_{p(k)})$ and $\de_{k+1}(y_k)$, respectively. Hence, it follows from \eqref{hdotdelk} that
$$\xh\mbox{-deg}(y_k c'_k) = \chi_{y_k} + \chi_{x_{k+1}} + \chi_{y_{p(k)}} = \xh\mbox{-deg} (y_{p(k)} c_{k+1} ).$$
Thus, $-y_k c'_k + y_{p(k)} c_{k+1}$ is homogeneous, as desired.

Finally, whether $p(k)= -\infty$ or not, it follows from \eqref{brackyx} that 
\begin{equation}
\begin{aligned}
\{y_k y'_k , x_j\} &= - (\al_{jk} + \al'_{jk}) x_j (y_k y'_k), \\
\{y_{p(k)} y_{k+1} , x_j\} &= - (\al_{j,p(k)} + \al_{j,k+1}) x_j (y_{p(k)} y_{k+1}),
\end{aligned}
\; \; \forall j \in [1,k-1].
\end{equation}
Since
\begin{align*}
\al_{jk} + \al'_{jk} &= \sum_{m=0}^{O_-(k)} (\la_{j,p^m(k)} + \la'_{j,p^m(k)})  \\
 &= 
 \la_{j,k+1} + \sum_{m=0}^{O_-(k)} \la_{j,p^m(k)} + \sum_{m=1}^{O_-(k)} \la_{j,p^m(k)}  = 
\al_{j,k+1} + \al_{j,p(k)} ,
\end{align*}
 we obtain $\{ y_k y'_k - y_{p(k)} y_{k+1} , x_j\} = - (\al_{j,p(k)} + \al_{j,k+1}) x_j \bigl( y_k y'_k - y_{p(k)} y_{k+1} \bigr)$. This shows that $y_k y'_k - y_{p(k)} y_{k+1}$ is 
a Poisson-normal element of $R_{k-1}$ which Poisson-normalizes the elements of $R_{k-1}$ in exactly the same way as $y_{p(k)} y_{k+1}$, and completes the proof of the theorem.
\end{proof}

Our next result turns the conclusion of \thref{swapkk+1}(b) into a cluster 
mutation statement. We define
\begin{equation}
\label{Pkset}
P(k) := \{ j \in [1,k] \mid s(j) > k \}, \; \; \forall k \in [1,N].
\end{equation}
Then $\{ y_j \mid j \in P(k) \}$ is a list of the homogeneous Poisson-prime elements of $R_k$ up to scalar multiples, and $|P(N)| = \rank R$ \eqref{rankRs>N}.

\bth{swap2} In the setting of Theorem {\rm\ref{tswapkk+1}(b)}, there exist $\kappa \in \kx$ and
a collection of nonnegative integers
$\{ m_i \mid i \in P(k-1), \, i \neq p(k)\}$ 
such that 
\begin{equation}
\label{thm2-disp}
y'_k = y_k^{-1} \Big( y_{p(k)} y_{k+1} + \kappa 
\prod_{i \in P(k-1),\,  i \neq p(k)} y_i^{m_i} \Big).
\end{equation}
If $p(k) = -\infty$, then all $m_i = 0$.
\eth

\begin{proof} By \thref{swapkk+1}(b), $y_k y'_k - y_{p(k)} y_{k+1}$
is a homogeneous Poisson-normal element of $R_{k-1}$. Applying \prref{factorHPUFD} and Theorems \ref{tPCGL-UFD}, \ref{tmainPprime}, we obtain 
$$
y_k y'_k - y_{p(k)} y_{k+1} =
\kappa \prod_{i \in P(k-1)} y_i^{m_i}
$$
for some $\kappa \in \KK$ and a collection of nonnegative integers
$\{ m_i \mid i \in P(k-1) \}$. Recall from \eqref{eltp(k)-inf} that if $p(k) = -\infty$, then $y_k y'_k - y_{p(k)} y_{k+1}$ is a nonzero scalar. Hence, $\kappa \ne 0$ and $m_i = 0$ for all $i \in P(k-1)$ in this case.
We need to prove in general that $\kappa \neq 0$, 
and that $m_{p(k)} = 0$ if $p(k) \neq - \infty$.

Suppose that $\kappa = 0$. Then  
$$
y_k y'_k - y_{p(k)} y_{k+1} = 0,
$$ 
which is a contradiction since 
$y_{k+1}$ is a prime element of $R_{k+1}$ which does not divide either $y_k$ or $y'_k$ (for $y_{k+1} \nmid y'_k$, recall \eqref{ykyk'}).

Now suppose that $p(k) \neq - \infty$ and $m_{p(k)} \neq 0$. Then
$y_{p(k)}$ is a prime element of $R_{k-1}$ and 
$$
y_k y'_k - y_{p(k)} y_{k+1}  \in y_{p(k)} R_{k-1}.
$$
Hence,
$$
y_k y'_k \in y_{p(k)} R_{k+1}.
$$
Furthermore, by \thref{mainPprime}, $y_k = y_{p(k)} x_k - c_k$ 
and $y'_k = y_{p(k)} x_{k+1} - c'_k$ for some 
$c_k,  c'_k \in R_{k-1}$ such that $y_{p(k)}$ does not divide $c_k$ or 
$c'_k$. But since
$$
c_k c'_k = y_k y'_k - y_{p(k)} x_k y'_k + c_k y_{p(k)} x_{k+1} 
\in y_{p(k)} R_{k+1} \cap R_{k-1} = y_{p(k)} R_{k-1},
$$
this contradicts the fact that $y_{p(k)}$ 
is a prime element of $R_{k-1}$.
\end{proof}

\subsection{Almost cluster mutations between Poisson-CGL extension presentations}
\label{almostmut}
Let $R$ be a $\KK$-algebra with a P-CGL extension presentation
as in \eqref{firstPCGL}. It is easy to see that the assumption that \eqref{secondPCGL} is a second P-CGL extension 
presentation of $R$ is equivalent to the following condition:
\begin{enumerate}
\item[(i)] $\de_{k+1} (x_j) \in R_{k-1}$ for all $j \in [1,k]$, and there exists 
$h'_{k+1} \in \Lie \HH$ satisfying $h'_{k+1} \cdot x_k = \la'_{k+1} x_k$ for some $\la'_{k+1} \in \kx$ and
$h'_{k+1} \cdot x_j = \la_{kj} x_j$ for all $j \in [1,k-1]\cup \{ k+1 \}$. 
\end{enumerate}
The rest of the data for the second P-CGL extension presentation of $R$
($\KK$-torus $\HH'$, scalars $\la'_{kj} \in \KK$, and elements $h'_k \in \Lie \HH'$) 
is given by \leref{PCGLtranspose} and the following:  
\begin{enumerate}
\item[(ii)] The torus $\HH'$ acting on the second P-CGL extension presentation 
can be taken as the original torus $\HH$. The corresponding elements $h'_j \in \Lie \HH$ 
are given by (i) for $j = k+1$ and $h'_j = h_j$ for $j \neq k+1$.
\item[(iii)] $\sig'_{k+1} = (h'_{k+1} \cdot)|_{R'_k}$, 
where $R'_k$ is the unital $\KK$-subalgebra of $R$ 
generated by $x_1, \ldots, x_{k-1}$ and $x_{k+1}$.
\item[(iv)] $\la'_{lj} = \la_{(k,k+1)(l), (k, k+1) (j)}$ for $l,j \in [1,N]$.
\end{enumerate}

{\em{Here and below, $(k,k+1)$ denotes a transposition 
in the symmetric group $S_N$, and $S_N$ is embedded in $GL_N(\Zset)$ 
via permutation matrices.}}

We continue with the hypotheses and notation of Theorems \ref{tswapkk+1}, \ref{tswap2}. Define maps $Y,Y' : \Zset^N \rightarrow \Fract(R)$ by
\begin{equation}
\label{YY'def}
Y(f) := \prod_{i=1}^N y_i^{f_i} \; \; \text{and} \; \; Y'(f) := \prod_{i=1}^N (y'_i)^{f_i}, \; \; \forall f = (f_1,\dots,f_N) \in \Zset^N.
\end{equation}
We also adopt the convention that
\begin{equation}
\label{e-inf}
e_{-\infty} = \ol{e}_{-\infty} := 0.
\end{equation}

\bth{swap3} Assume the setting of Theorem {\rm\ref{tswapkk+1}}.

{\rm(a)} If $\eta(k) \neq \eta(k+1)$, then $y'_j = y_{(k, k+1) j}$ for all $j \in [1,N]$, i.e., $Y' = Y(k,k+1)$.

{\rm(b)} If $\eta(k) = \eta(k+1)$, then $Y'(e_j) = y'_j = y_j= Y(e_j)$ for all $j \neq k$ and 
\begin{equation}
\label{YY'}
Y'(e_k)= y'_k = Y ( - e_k + e_{p(k)} + e_{k+1}) + \kappa
Y \Big( - e_k + \sum_{i \in P(k-1), \, i \neq p(k)} m_i e_i \Big)
\end{equation}
for the collection of nonnegative integers $\{m_i \mid i \in P(k-1), \; i \neq p(k)\}$
and the scalar $\kappa \in \kx$ from Theorem {\rm\ref{tswap2}}.
\eth

\begin{proof} This follows from Theorems \ref{tswapkk+1} and \ref{tswap2}. 
\end{proof}

\subsection{Scalars associated to mutations of Poisson-prime elements}
\label{scalars}
Next, we derive formulas for certain scalars which will appear in the seeds which we 
construct in Section \ref{main}. 

Recall the skew-symmetric matrix $\qb = (q_{kj}) \in M_N(\KK)$ from \eqref{qkj}, and define the corresponding skew-symmetric bicharacter $\Om_\qb : \Zset^N \times \Zset^N \rightarrow \KK$ by \eqref{Om}, so that
\begin{equation}
\label{Ombfq}
\Om_\qb(e_k, e_j) = q_{kj} = \Om_\lab( \ol{e}_k, \ol{e}_j), \; \; \forall j,k \in [1,N].
\end{equation}
In view of \prref{PtorusR}, we have
\begin{equation}
\label{brackYfYg}
\{ Y(f), Y(g) \} = \Om_\qb(f,g) Y(f) Y(g), \; \; \forall f,g \in \Zset^N.
\end{equation}

Recall from \S \ref{notaconv} that for an $\HH$-eigenvector $u \in R$, $\chi_u \in \xh$ denotes
its eigenvalue.

\bth{Omqval} Assume the setting of Theorem {\rm\ref{tswap2}}.
Then
\begin{equation}
\label{Omqval1}
\Om_\qb \Big(e_{p(k)} + e_{k+1} - \sum_{i \in P(k-1),\, i \neq p(k)} m_i e_i, \,
e_j \Big) = 0, \; \; \forall j \neq k
\end{equation}
and
\begin{equation}
\label{Omqval2}
\begin{aligned}
\Om_\qb \Big(e_{p(k)} + e_{k+1} - \sum_{i \in P(k-1),\, i \neq p(k)} m_i e_i, \,
e_k \Big) = -\la_{k+1}= \la'_{k+1}
 \end{aligned}
\end{equation}
for the collection of nonnegative integers $\{m_i \mid i \in P(k-1), \, i \neq p(k)\}$
from Theorems {\rm\ref{tswap2}} and {\rm\ref{tswap3}(b)}.
\eth

\begin{proof} Denote for brevity the elements
$$
g := \sum_{i \in P(k-1), \, i \neq p(k)} m_i e_i \quad \text{and} \quad
\ol{g} := \sum_{i \in P(k-1), \, i \neq p(k)} m_i \ol{e}_i \quad \text{in} \quad \Zset^N.
$$

For $j\in [1,N]$ with $j\ne k$, we have $y_j = y'_j$ by \eqref{y-equal}, and so
\begin{align*}
q'_{kj} y'_k y_j &= q'_{kj} \bigl( Y(-e_k +e_{p(k)} +e_{k+1}) + \kappa Y(-e_k +g) \bigr) y_j  \\
q'_{kj} y'_k y_j &= \{ y'_k, y_j\} = \big\{ Y(-e_k +e_{p(k)} +e_{k+1}) + \kappa Y(-e_k +g), y_j \big\}  \\
 &= \Om_\qb (-e_k +e_{p(k)} +e_{k+1}, e_j) Y(-e_k +e_{p(k)} +e_{k+1}) y_j  \\
  &\qquad\qquad+  \kappa \Om_\qb (-e_k +g, e_j) Y(-e_k +g) y_j .
 \end{align*}
Consequently,
$$
\Om_\qb (- e_k + e_{p(k)} + e_{k+1}, e_j) = \Om_\qb(- e_k + g, e_j) = q'_{kj} ,
$$
which implies \eqref{Omqval1}.

Applying \eqref{Omqval1} for 
$j = k+1$ leads to
\begin{align}
\Om_\qb (e_{p(k)} +\, &e_{k+1} - g, e_k)  \notag  \\
 &= \Om_\qb (e_{p(k)} + e_{k+1} - g, e_k) - \Om_\qb (e_{p(k)} + e_{k+1} - g, e_{k+1})  \notag  \\
&= \Om_\lab (\ol{e}_{p(k)} + \ol{e}_{k+1} - \ol{g}, \ol{e}_k) - 
\Om_\lab (\ol{e}_{p(k)} + \ol{e}_{k+1} - \ol{g}, \ol{e}_{k+1})  \label{Omeq}  \\
&= - \Om_\lab (\ol{e}_{p(k)} + \ol{e}_{k+1} - \ol{g}, e_{k+1}) = - \Om_\lab (\ol{e}_{p(k)} + \ol{e}_k - \ol{g}, e_{k+1}).  \notag
\end{align}
Since $y'_k$ is an $\HH$-eigenvector, it follows from \eqref{thm2-disp} that $Y(e_{p(k)} +e_{k+1})$ and $Y(g)$ have the same $\xh$-degree, whence $\chi_{Y(e_{p(k)} + e_{k+1} -g)} = 0$. Using that $\chi_{Y(e_{p(k)} + e_{k+1} - g)} = \chi_{Y(e_{p(k)} + e_k -g)} + \chi_{x_{k+1}}$,
we obtain
$$
h_{k+1} . ( Y(e_{p(k)} + e_k -g) x_{k+1} ) = 0
$$
from \eqref{h.0}. Eq.~\eqref{sigkyj} and the definition of $\la_{k+1}$ then yield
\begin{align*}
0 &= \al_{k+1,p(k)} +\al_{k+1, k} - \sum_{i \in P(k-1),\, i \neq p(k)} m_i \al_{k+1, i} + \la_{k+1}  \\
 &= \Om_\lab( e_{k+1}, \ol{e}_{p(k)} +\ol{e}_k - \ol{g}) + \la_{k+1} .
\end{align*}
Combining this with \eqref{Omeq} proves the first equality in \eqref{Omqval2}.

It follows from \eqref{Omqval1} and \eqref{brackYfYg} that
$$
\{ Y(e_{p(k)} + e_{k+1} - g), y_j \} = 0, \; \; \forall j \ne k,
$$
whence $Y(e_{p(k)} + e_{k+1} - g)$ Poisson-commutes with $Y(g)$, and so $Y(e_{p(k)} + e_{k+1} - g)$ also Poisson-commutes with $Y( e_{p(k)}+ e_{k+1}) = y_{p(k)} y_{k+1}$. Since $y_k y'_k$ is a linear combination of $y_{p(k)} y_{k+1}$ and $Y(g)$, we get
\begin{equation}
\label{brackstuff}
\{ Y(e_{p(k)} + e_{k+1} - g) , y_k y'_k \} = 0.
\end{equation}

By \eqref{brackYfYg} and the first equality in \eqref{Omqval2},
$$
\{  Y(e_{p(k)} + e_{k+1} - g) , y_k \} = - \la_{k+1}  Y(e_{p(k)} + e_{k+1} - g) y_k .
$$
Interchanging the roles of the P-CGL extension presentations \eqref{firstPCGL} and \eqref{secondPCGL} and using the symmetric nature 
of the assumptions of \thref{swapkk+1}(b) shows that 
$$
\{  Y'(e_{p(k)} + e_{k+1} - g) , y'_k \} = - \la'_{k+1}  Y'(e_{p(k)} + e_{k+1} - g) y'_k .
$$
Since $Y'(e_{p(k)} + e_{k+1} - g) = Y(e_{p(k)} + e_{k+1} - g)$, we obtain
$$
\{ Y(e_{p(k)} + e_{k+1} - g) , y_k y'_k \} = - (\la_{k+1} + \la'_{k+1}) Y(e_{p(k)} + e_{k+1} - g) y_k y'_k .
$$
Therefore $\la_{k+1} + \la'_{k+1} = 0$ because of \eqref{brackstuff}, which proves the final equality in \eqref{Omqval2}.
\end{proof}

\sectionnew{Homogeneous Poisson-prime elements for subalgebras of symmetric Poisson-CGL extensions}
\label{Pprime-sym}

Each symmetric Poisson-CGL extension $R$ of length $N$ has many different 
P-CGL extension presentations given by \eqref{taupOre}. They are parametrized by 
the elements of the subset $\Xi_N$ of $S_N$, cf. \eqref{tau}.
In order to phrase \thref{swap2} into a mutation statement 
between cluster variables associated to the elements of $\Xi_N$ and 
to make the scalars $\kappa$ from \thref{swap2} equal to one, 
we need a good picture of the sequences of homogeneous Poisson-prime 
elements $y_1, \ldots, y_N$ from \thref{mainPprime} 
associated to each presentation \eqref{taupOre}. 
This is obtained in \thref{sym-Pprime}. \thref{y-int} contains
a description of the homogeneous Poisson-prime elements that enter
into this result. Those Poisson-prime elements (up to rescaling)
comprise the cluster variables that 
will be used in Section \ref{main} to construct 
cluster algebra structures on symmetric 
P-CGL extensions. Along the way, we explicitly describe the elements 
of $\Xi_N$ and prove an invariance property of the scalars 
$\la_l$ and $\la_l^*$ from Definitions \ref{dPCGL} and \ref{dsymmPCGL}.
\thref{u-prod}, which appears at the end of the section, contains
a key result used in the next section to normalize the generators $x_j$ 
of symmetric P-CGL extensions so that all scalars $\kappa$ in \thref{swap2} 
become equal to one.

Throughout the section, $\eta$ will denote a function $[1,N] \rightarrow \Zset$ satisfying the conditions of \thref{mainPprime}, with respect to the original P-CGL extension presentation \eqref{itpOre} of $R$, and $p$ and $s$ will denote the corresponding predecessor and successor functions. We will repeatedly use the one-line notation \eqref{one-line} for permutations.

\subsection{The elements $y_{[i,s^m(i)]}$}
\label{yismi-elts}
Recall from \S\ref{symmP-CGL} that for a symmetric P-CGL extension $R$ of rank $N$
and $1 \leq j \leq k \leq N$, $R_{[j,k]}$ denotes the unital subalgebra 
of $R$ generated by $x_j , \ldots, x_k$. It is a Poisson-Ore extension of both 
$R_{[j,k-1]}$ and $R_{[j+1,k]}$ \eqref{Rjkpres}. All such subalgebras 
are (symmetric) P-CGL extensions and \thref{mainPprime} applies 
to them. 

For $i \in [1,N]$ and $0 \le m \le O_+(i)$, recall \eqref{O-+} (i.e., 
$s^m(i) \in [1,N]$), set 
\begin{equation}
\label{e-int}
e_{[i, s^m(i)]} := e_i + e_{s(i)} + \cdots + e_{s^m(i)} \in \Zset^N.
\end{equation}
The vectors \eqref{ekbar} are special cases of these:
$$
\ol{e}_k = e_{[p^{O_-(k)}(k), k]}, \quad
\forall k \in [1,N]. 
$$
We also set $e_\varnothing := 0$.
The next theorem treats the Poisson-prime elements that will appear as 
cluster variables for symmetric CGL extensions. It will be proved in \S \ref{proof.yint}.

\bth{y-int} Assume that $R$ is a symmetric Poisson-CGL extension of length $N$, and $i \in [1,N]$ and 
$m \in \Zset_{\geq 0}$ are such that $s^m(i) \in [1,N]$, i.e., 
$s^m(i) \neq + \infty$. Then the following hold:

{\rm(a)} All homogeneous Poisson-prime elements of $R_{[i, s^m(i)]}$ 
that do not belong to $R_{[i, s^m(i)-1]}$ are associates of each other.

{\rm(b)} All homogeneous Poisson-prime elements of $R_{[i, s^m(i)]}$ 
that do not belong to $R_{[i+1, s^m(i)]}$ are associates of each other.
In addition, the set of these homogeneous Poisson-prime elements coincides 
with the set of homogeneous Poisson-prime elements in part {\rm(a)}.

{\rm(c)} The homogeneous Poisson-prime elements in {\rm(a)} and {\rm(b)} have leading terms of the form 
$$
\xi x^{e_{[i,s^m(i)]}} = \xi x_{i} x_{s(i)} \ldots x_{s^m(i)}
$$ 
for some $\xi \in \kx$, see \S{\rm\ref{PprimePCGL}}. For each $\xi \in \kx$, there 
is a unique homogeneous Poisson-prime element of $R_{[i,s^m(i)]}$ with such a leading term.
Denote by $y_{[i,s^m(i)]}$ the homogeneous
Poisson-prime element of $R_{[i, s^m(i)]}$ with leading term $x_{i} x_{s(i)} \ldots x_{s^m(i)}$.
Let $y_\varnothing :=1$.

{\rm(d)} We have 
$$
y_{[i,s^m(i)]} = y_{[i, s^{m-1}(i)]} x_{s^m(i)} - c_{[i,s^m(i)-1]}
= x_i y_{[s(i), s^m(i)]} - c'_{[i+1, s^m(i)]}
$$
for some $c_{[i, s^m(i)-1]} \in R_{[i,s^m(i)-1]}$ 
and $c'_{[i+1, s^m(i)]} \in R_{[i+1, s^m(i)]}$.

{\rm(e)} For all $k \in [1,N]$ such that $p(i)< k< s^{m+1}(i)$, we have
$$
\{ y_{[i,s^m(i)]} , x_k \} = \Om_\lab( e_{[i,s^m(i)]}, e_k ) x_k y_{[i,s^m(i)]} .$$
\eth

The case $m=0$ of this theorem is easy to verify.
In that case, $y_{[i,i]}= x_i$, and statement (e) follows by applying \eqref{bracketyjxk} in $R_{[\min\{i,k\}, \max\{i,k\}]}$.

\bex{OMmn5}
In the case of $R = \OMmn$, the elements $y_{[i,s^m(i)]}$ of \thref{y-int} are solid minors, just as in \exref{OMmn2}. More precisely, if $i = (r-1)n+c$ for some $r \in [1,m]$ and $c \in [1,n]$ and $s^l(i) \ne +\infty$, then $s^l(i) = (r+l-1)n + c+l$ with $r+l \le m$ and $c+l \le n$, and $y_{[i, s^l(i)]} = \Delta_{[r, r+l], [c, c+l]}$.
\eex

The following theorem describes the $y$-sequences from \thref{mainPprime} associated to the 
P-CGL extension presentations \eqref{taupOre} in terms of the Poisson-prime elements 
from \thref{y-int}. It will be proved in \S \ref{proof.yint}. Recall that 
for every $\tau \in \Xi_N$ and $k \in [1,N]$, the set $\tau([1,k])$ is an 
interval.

\bth{sym-Pprime} Assume that $R$ is a symmetric Poisson-CGL extension of length $N$
and $\tau$ an element of the subset $\Xi_N$ of $S_N$, cf. \eqref{tau}.
Let $y_{\tau, 1}, \ldots, y_{\tau, N}$ be the sequence 
in $R$ from Theorem {\rm\ref{tmainPprime}} applied to the P-CGL extension presentation 
\eqref{taupOre} of $R$ corresponding to $\tau$. Let $k \in [1,N]$. 

{\rm (a)} If $\tau(k) \geq \tau(1)$, then $y_{\tau, k} = y_{[p^m(\tau(k)), \tau(k)]}$, where 
\begin{equation}
\label{m-p}
m = \max \{ n \in \Znn \mid p^n( \tau(k)) \in \tau([1, k]) \}.
\end{equation}

{\rm (b)} If $\tau(k) \leq \tau(1)$, then $y_{\tau, k} = y_{[\tau(k), s^m(\tau(k))]}$, where 
\begin{equation}
\label{m-s}
m = \max \{ n \in \Znn \mid s^n( \tau(k)) \in \tau([1, k]) \}.
\end{equation}
In both cases, the predecessor and successor functions are with 
respect to the original P-CGL extension presentation 
\eqref{itpOre} of $R$.
\eth

\subsection{The elements of $\Xi_N$}
\label{eltsXiN}
In this and the next subsection, we investigate the 
elements of the subset $\Xi_N$ of $S_N$ defined 
in \eqref{tau}. It follows from \eqref{tau} that every element $\tau \in \Xi_N$ 
has the property that either $\tau(N) = 1$ or $\tau(N) = N$. This implies the following recursive 
description of $\Xi_N$.

\ble{SigN-ind} For each $\tau \in \Xi_N$, there exists $\tau' \in \Xi_{N-1}$ 
such that either
$$
\tau(i) = \tau'(i), \; \; \forall i \in [1,N-1] \quad 
\mbox{and} \quad \tau(N) = N
$$
or 
$$
\tau(i) = \tau'(i) + 1, \; \;  \forall i \in [1,N-1]
\quad \mbox{and} \quad \tau(N) =1.
$$
For all $\tau' \in \Xi_{N-1}$, the above define elements of $\Xi_N$.
\ele

Given $k \in [1,N]$ and a sequence $k \leq j_k \leq \cdots \leq j_1 \leq N$, 
define
\begin{equation}
\label{tau-seq}
\tau_{(j_k, \ldots, j_1)} 
:= (k (k+1) \ldots j_k) 
\ldots
(2 3 \ldots j_2)
(1 2 \ldots j_1) \in S_N,
\end{equation}
where in the right hand side we use the standard notation 
for cycles in $S_N$. Using \leref{SigN-ind}, the elements of $\Xi_N$ are easily characterized as follows.

\ble{charSigN} {\rm \cite[Lemma 5.4]{GYbig}}  The subset $\Xi_N\subset S_N$ consists 
of the elements of the form $\tau_{(j_k, \ldots, j_1)}$, 
where $k \in [1,N]$ and $k \leq j_k \leq \cdots \leq j_1 \leq N$.  
\ele

The representation of an element of $\Xi_N$ 
in the form \eqref{tau-seq} is not unique. As observed in \cite{GYbig}, one way to visualize 
$\tau_{(j_k, \ldots, j_1)}$ is that the sequence $\tau(1), \ldots, \tau(N)$ 
is obtained from the sequence $1, \ldots, N$ by the following procedure:

$(*)$ {\em{The number $1$ is pulled to the right to position $j_1$ (preserving the order of the other numbers), 
then the number $2$ is pulled to the right to position 
$j_2-1$, ..., at the end the number $k$ is pulled to the right to position $j_k-k+1$.}}

For example, for $k=2$ the following illustrates how $\tau_{(j_2, j_1)}$ is obtained from the 
identity permutation:
\begin{align*}
&[\circled{1}, \circled{2},3, 4, \ldots, j_2, j_2+1, \ldots, j_1, j_1+1, \ldots, N] \mt
\\
&[\circled{2},3, 4, \ldots, j_2, j_2+1, \ldots, j_1, \circled{1}, j_1+1, \ldots, N] \mt
\\
&[3, 4, \ldots, j_2, \circled{2}, j_2+1, \ldots, j_1, \circled{1}, j_1+1, \ldots, N],
\end{align*}
where the numbers that are pulled ($1$ and $2$) are circled.

If we perform the above procedure one step at a time, so in each step we only interchange the positions 
of two adjacent numbers, then the elements of $S_N$ from all intermediate steps will 
belong to $\Xi_N$. For example, this requires factoring the cycle $(1,2,\dots,j_1)$ as $(1,j_1)(1,j_1-1) \cdots (1,2)$ rather than as $(1,2) (2,3) \cdots (j_1-1,j_1)$. This implies at once the first part of the next corollary \cite[Corollary 5.5(a)]{GYbig}.
 
\bco{steps} Let $R$ be a symmetric Poisson-CGL extension of length $N$ and $\tau \in \Xi_N$. 

{\rm(a)} There exists a sequence $\tau_0 = \id, \tau_1, \ldots, \tau_n = \tau$ in $\Xi_N$ such that 
for all $l \in [1,n]$,
$$
\tau_l = ( \tau_{l-1}(k_l), \tau_{l-1} (k_l+1)) \tau_{l-1} = \tau_{l-1}(k_l, k_l+1)
$$
for some $k_l \in [1,N-1]$ such that $\tau_{l-1}(k_l) < \tau_{l-1}(k_l +1)$.

{\rm(b)} If $\eta : [1,N] \to \Zset$ is a function satisfying the conditions of 
Theorem {\rm\ref{tmainPprime}} for the original P-CGL presentation of $R$, then 
\begin{equation}
\label{etatau}
\eta_\tau:=\eta \tau : [1,N] \to \Zset 
\end{equation}
satisfies the conditions of Theorem {\rm\ref{tmainPprime}} for the $\eta$-function of the 
P-CGL extension presentation \eqref{taupOre} of $R$ corresponding to $\tau$.
\eco

The sequence described in part (a) of the corollary for the element
$\tau_{(j_k, \ldots, j_1)} \in \Xi_N$ has length $n= j_1 + \cdots + j_k - k(k+1)/2$. 
The second part of the corollary follows by recursively applying \thref{swapkk+1} to the 
P-CGL extension presentations \eqref{taupOre} for the elements
$\tau_{l-1}$ and $\tau_l$. \coref{steps}(b) gives a second proof of the 
fact that the rank of a symmetric P-CGL extension does not depend on the choice 
of P-CGL extension presentation of the form \eqref{taupOre}, see \S \ref{symmP-CGL}.

\bco{restr-eta} Assume that $R$ is a symmetric Poisson-CGL extension of length $N$ 
and $\eta : [1,N] \to \Zset$ is a function satisfying the conditions of 
Theorem {\rm\ref{tmainPprime}}. Then for all $1 \leq j \leq k \leq N$, 
the function $\eta_{[j,k]} \colon [1,k-j+1] \to \Zset$ given by 
$$
\eta_{[j,k]}(l) = \eta(j+l-1), \; \; \forall l \in [1, k-j+1]
$$
satisfies
the conditions of Theorem {\rm\ref{tmainPprime}} for the symmetric P-CGL extension 
$R_{[j,k]}$.
\eco

Thus, the $\eta$-functions for ``interval 
subalgebras'' of symmetric P-CGL extensions can be chosen to be
restrictions of the original $\eta$-function up to shifts.

\begin{proof}
Apply \coref{steps}(b) to 
$$
\tau = [j, \ldots, k, j-1, \ldots, 1, k+1, \ldots, N] \in \Xi_N 
$$
and consider the $(k-j+1)$-st subalgebra of the corresponding Poisson-Ore extension,
which equals $R_{[j,k]}$.
\end{proof}

\subsection{A subset of $\Xi_N$}
\label{GammaN}
 We will work with a certain subset $\Ga_N$ of $\Xi_N$ which was introduced in \cite[\S5.3]{GYbig} and which will 
play an important role in Section \ref{main}. Recall \eqref{tau-ci} that $w_\ci = [N, N-1, \ldots, 1]$ denotes the longest element of $S_N$.
For $1 \leq i \leq j \leq N$, define the following elements of $\Xi_N$:
\begin{equation}
\label{tauij}
\tau_{i,j} := [i+1, \ldots, j, i, j+1, \ldots, N, i-1, i-2, \ldots, 1] \in \Xi_N.
\end{equation}
They satisfy
$$
\tau_{1,1} = \id, \quad \tau_{i,N} = \tau_{i+1, i+1}, \; \forall i \in [1,N-1], \quad
\tau_{N,N} = w_\ci.
$$
Denote by $\Ga_N$ the subset of $\Xi_N$ consisting 
of all $\tau_{i,j}$'s and consider the 
following linear ordering on it
\begin{multline}
\label{sequence}
\Ga_N := \{
\id = \tau_{1,1} \prec \ldots \prec \tau_{1,N} = \tau_{2,2} \prec \ldots \prec \tau_{2, N} = \tau_{3,3} \prec \ldots \prec \\ 
\prec \tau_{3,N}= \tau_{4,4} \prec \ldots \prec
\tau_{N-2, N} = \tau_{N-1,N-1} \prec \tau_{N-1,N} = \tau_{N,N} = w_\ci \}.
\end{multline}
In the notation of \eqref{tau-seq}, the elements of $\Ga_N$ are given by
$$
\tau_{i,j} := \tau_{(j, N, \ldots, N)}, \, \, \text{where} \; N \; \text{is repeated} \; i-1 \; \text{times}, \, \, \forall 1\le i \le j \le N.
$$
The sequence of elements \eqref{sequence} is nothing but the sequence of the intermediate steps 
of the procedure $(*)$ from \S \ref{eltsXiN} applied to the longest element $w_\ci \in S_N$ (in which 
case $k=N-1$ and $j_1= j_2 = \cdots = j_{N-1} = N$). Note that
$$\tau_{i,l+1} = \tau_{i,l} (l+2-i, \, l+1-i), \; \; \forall 1 \le i \le l < N.$$

Assume that $R$ is a symmetric P-CGL extension.  
To each element of $\Ga_N$, \prref{tauPCGL} associates a P-CGL extension presentation of $R$.  
Each two consecutive presentations are associated to a pair $\tau_{i,j}, \tau_{i,j+1} \in \Xi_N$ 
for some $1 \leq i \leq j < N$. They have the forms
\begin{equation}
\label{present1}
\begin{aligned}
R = 
\KK[x_{i+1}]_p \cdots [x_j; \sig_j, \de_j]_p &[x_i; \sig^\sy_i, \de^\sy_i]_p 
[x_{j+1}; \sig_{j+1}, \de_{j+1}]_p \cdots 
\\
&[x_N; \sig_N, \de_N]_p 
[x_{i-1}; \sig^\sy_{i-1}, \de^\sy_{i-1}]_p \cdots [x_1; \sig^\sy_1, \de^\sy_1]_p 
\end{aligned}
\end{equation}
and
\begin{equation}
\label{present2}
\begin{aligned}
R= 
\KK[x_{i+1}]_p \cdots [x_j; \sig_j, \de_j]_p 
&[x_{j+1}; \sig_{j+1}, \de_{j+1}]_p [x_i; \sig^\sy_i, \de^\sy_i]_p [x_{j+2}; \sig_{j+2}, \de_{j+2}]_p 
\cdots 
\\
&[x_N; \sig_N, \de_N]_p 
[x_{i-1}; \sig^\sy_{i-1}, \de^\sy_{i-1}]_p \cdots [x_1; \sig^\sy_1, \de^\sy_1]_p, 
\end{aligned}
\end{equation}
respectively. These two presentations satisfy the assumptions of \thref{swapkk+1}.
Recall from Definitions \ref{dPCGL} and \ref{dsymmPCGL}
that $\sig_l = (h_l \cdot)$, $\sig^\sy_l = (h^\sy_l \cdot)$ 
and $h_l \cdot x_l = \la_l x_l$, $h^\sy_l \cdot x_l = \la^\sy_l x_l$
with $h_l, h^\sy_l \in \Lie \HH$ and $\la_l, \la^\sy_l \in \kx$. By \coref{steps}(b), the $\eta$-function for the P-CGL extension presentation \eqref{present1} can be taken to be $\eta \tau_{i,j}$. In particular, the values of this function on $j+1-i$ and $j+2-i$ are $\eta(i)$ and $\eta(j+1)$, respectively.
If $\eta(i)=\eta(j+1)$, then Eq. \eqref{Omqval2} of \thref{Omqval} applied to the 
presentations \eqref{present1} and \eqref{present2} implies 
\begin{equation}
\label{lala*}
\la^\sy_i + \la_{j+1} = 0.
\end{equation}

\bpr{la-equal} Let $R$ be a symmetric Poisson-CGL extension of length $N$
and $a \in \Zset$ be such that $|\eta^{-1}(a)|>1$. Denote
$$
\eta^{-1}(a) = \{l, s(l), \ldots, s^m(l) \} 
$$
where $l \in [1,N]$ and $m = O_+(l) \in \Zset_{> 0}$. Then 
\begin{equation}
\label{la-eq}
\la^*_l = \la^*_{s(l)} = \cdots = \la^*_{s^{m-1}(l)} = 
- \la_{s(l)} = - \la_{s^2(l)}= \cdots = - \la_{s^m(l)}.
\end{equation}
\epr

\begin{proof} It follows from \eqref{lala*} that 
$$
\la^\sy_{s^{m_1}(l)} = - \la_{s^{m_2}(l)}, \quad \forall 0 \leq m_1 < m_2 \leq m
$$
which is equivalent to the statement of the proposition. 
\end{proof}

\subsection{Sequences of homogeneous Poisson-prime elements} 
\label{proof.yint}
We proceed with the proofs of the two theorems formulated in \S \ref{yismi-elts}. 

\begin{proof}[Proof of Theorem {\rm\ref{ty-int}}] As already noted, the case $m=0$ is easily verified. Assume now that $m > 0$.

(a) Consider the P-CGL extension 
presentation of $R$ associated to 
\begin{multline}
\label{y-int1}
\tau_{(s^m(i)-1, s^m(i), \ldots, s^m(i))} = \\ 
[i +1 , \ldots, s^m(i)-1, i, s^m(i), i-1, \ldots, 1, s^m(i)+1, \ldots, N] \in \Xi_N,
\end{multline}
where $s^m(i)$ is repeated $i-1$ times.
The $(s^m(i)-i)$-th and $(s^m(i)-i+1)$-st algebras in the chain 
are precisely $R_{[i, s^m(i)-1]}$ and $R_{[i,s^m(i)]}$.
\thref{mainPprime} implies that the homogeneous Poisson-prime elements of $R_{[i,s^m(i)]}$ 
that do not belong to $R_{[i, s^m(i)-1]}$ are associates of each other.
Denote by $z_{[i, s^m(i)]}$ one such element. Set $z_\varnothing = 1$.
(One can prove part (a) using 
the simpler presentation of $R$ associated to the permutation 
$$[i, i + 1 , \ldots, s^m(i), 
i-1, \ldots, 1, s^m(i)+1, \ldots, N] \in \Xi_N$$
but 
the first presentation will also play a role in the proof of part (b).) The equivariance 
property of the $\eta$-function from \coref{steps}(b) and \thref{mainPprime} imply that 
\begin{equation}
\label{y-int2}
z_{[i,s^m(i)]} = \xi_{i,m}( z_{[i,s^{m-1}(i)]} x_{s^m(i)} - c_{[i,s^m(i)-1]})
\end{equation}
for some $c_{[i,s^m(i)-1]} \in R_{[i,s^m(i)-1]}$ and $\xi_{i,m} \in \kx$.

(b) Now consider the P-CGL extension 
presentation of $R$ associated to 
\begin{equation}
\label{y-int3}
\tau_{(s^m(i), \ldots, s^m(i))} =
[i +1 , \ldots, s^m(i), i, \ldots, 1, s^m(i)+1, \ldots, N] \in \Xi_N,
\end{equation}
where $s^m(i)$ is repeated $i$ times.
The $(s^m(i)-i)$-th and $(s^m(i)-i+1)$-st algebras in the chain 
are $R_{[i+1, s^m(i)]}$ and $R_{[i,s^m(i)]}$.
\thref{mainPprime} implies that the homogeneous Poisson-prime elements of $R_{[i,s^m(i)]}$ 
that do not belong to $R_{[i+1, s^m(i)]}$ are associates of each other.
They are also associates of $z_{[i,s^m(i)]}$. This follows from Eq.~\eqref{y-equal} of
\thref{swapkk+1}(b) applied to the P-CGL extension presentations of $R$ associated to 
the elements \eqref{y-int1} and \eqref{y-int3}. The fact that 
these two P-CGL extension presentations satisfy the assumptions of 
\thref{swapkk+1}(b) follows from  
the equivariance of the $\eta$-function from \coref{steps}(b).
This equivariance and \thref{mainPprime} applied to the presentation for \eqref{y-int3} 
also imply that 
\begin{equation}
\label{y-int4}
z_{[i,s^m(i)]} = \xi'_{i,m}( x_{i} z_{[s(i),s^m(i)]} - c'_{[i+1,s^m(i)]})
\end{equation}
for some $c'_{[i+1,s^m(i)]} \in R_{[i+1,s^m(i)]}$ and $\xi'_{i,m} \in \kx$.

Parts (c) and (d) follow at once by comparing the leading terms in 
Eqs. \eqref{y-int2} and \eqref{y-int4}, and using the fact that the group of units 
of an iterated Poisson-Ore extension over $\KK$ is reduced to scalars.

For part (e), we apply \eqref{bracketyjxk} in $R_{[\min\{k,i\},\max\{k,s^m(i)\}]}$ with $j = s^m(i)$.
\end{proof}

\begin{proof}[Proof of Theorem {\rm\ref{tsym-Pprime}}] We first show that in situation (a), $y_{\tau,k}$ is a scalar multiple of $y_{[p^m(\tau(k)), \tau(k)]}$.

For $k \in [0,N]$, 
denote by $R_{\tau, k}$ the $k$-th algebra in the chain
\eqref{taupOre}. Since $\tau([1,j])$ is an interval for all $j \leq k$,
$$
\tau([1, k]) = [\tau(i), \tau(k)] \quad
\mbox{for some} \; \; i \in [1,k].
$$
Therefore $R_{\tau, k} = R_{[\tau(i), \tau(k)]}$ and 
$R_{\tau, k-1} = R_{[\tau(i), \tau(k)-1]}$.
For $m \in \Zset_{\geq 0}$ given by \eqref{m-p} 
we have 
$$
\tau(i) \leq p^m(\tau(k)) \leq \tau(k).
$$
\thref{mainPprime} implies that $y_{[p^m(\tau(k)), \tau(k)]}$
is a homogeneous Poisson-prime element of $R_{[\tau(i),\tau(k)]} = R_{\tau, k}$. 
It does not belong to $R_{\tau, k-1} = R_{[\tau(i), \tau(k)-1]}$
because of \thref{y-int}(c). It follows from \thref{mainPprime} that 
the homogeneous Poisson-prime elements of $R_{\tau, k}$ that do not belong to $R_{\tau, k-1}$ are associates of  $y_{\tau, k}$.  
Hence, $y_{\tau, k}$ is a scalar multiple of 
$y_{[p^m(\tau(k)), \tau(k)]}$.

Analogously, in situation (b), $y_{\tau,k}$ is a scalar multiple of $y_{[\tau(k), s^m(\tau(k))]}$.

To prove the stated equalities, we proceed by induction on $k$. The case $k=1$ is clear, since $y_{\tau,1} = x_{\tau(1)} = y_{[\tau(1), \tau(1)]}$. Now assume that $k>1$, and suppose we are in situation (a). If $m=0$, then $y_{\tau,k} = x_{\tau(k)} = y_{[\tau(k), \tau(k)]}$, so we may assume that $m>0$.

Let $l \in [1,k-1]$ be maximal such that $\eta(\tau(l)) = \eta(\tau(k))$, and note that
$$
p(\tau(k)), \dots, p^m(\tau(k)) \in \tau([1,l]).
$$
By construction of $y_{\tau,k}$ and \thref{y-int}(d), we have
\begin{equation}
\label{y-compare}
y_{\tau,k} = y_{\tau,l} x_{\tau(k)} - c  \qquad \text{and} \qquad y_{[p^m(\tau(k)), \tau(k)]} = y_{[p^m(\tau(k)), p(\tau(k))]} x_{\tau(k)} - c'
\end{equation}
for some $c,c' \in R_{[1, \tau(k) -1]}$. If $\tau(l) \ge \tau(1)$, then $\tau([1,l]) \subseteq [1,\tau(l)]$, whence $\tau(l) = p(\tau(k))$ and 
$$
m-1 = \max \{ n \in \Znn \mid p^n(\tau(l)) \in \tau([1,l]) \}.
$$
By induction, $y_{\tau,l} = y_{[p^{m-1}(\tau(l)), \tau(l)]} = y_{[p^m(\tau(k)), p(\tau(k))]}$. On the other hand, if $\tau(l) \le \tau(1)$, then $\tau([1,l]) \subseteq [\tau(l), \tau(k)-1]$, whence $\tau(l) = p^m(\tau(k))$ and 
$$
m-1 = \max \{ n \in \Znn \mid s^n(\tau(l)) \in \tau([1,l]) \}.
$$
In this case, induction yields $y_{\tau,l} = y_{[\tau(l), s^{m-1}(\tau(l))]} = y_{[p^m(\tau(k)), p(\tau(k))]}$. In both cases, we conclude from \eqref{y-compare} and the equality $y_{\tau,l} = y_{[p^m(\tau(k)), p(\tau(k))]}$ that $y_{\tau,k} = y_{[p^m(\tau(k)), \tau(k)]}$. 

Situation (b), when $\tau(k) < \tau(1)$, is analogous and is left to the reader.
\end{proof}

\bco{int-brack} If $R$ is a symmetric Poisson-CGL extension of length 
$N$ and $i, j \in [1,N]$, $m, n \in \Znn$ are such that
$$
i \leq j \leq s^n(j) \leq s^m(i) \leq N,
$$
then 
$$
\{ y_{[i,s^m(i)]} , y_{[j, s^n(j)]} \} = 
\Om_\lab(e_{[i,s^m(i)]}, e_{[j,s^n(j)]})
y_{[j, s^n(j)]} y_{[i,s^m(i)]}.
$$
\eco

\begin{proof}
Apply \thref{sym-Pprime} to the P-CGL extension 
presentation of $R$ associated to 
$$
\tau':= [j, \ldots, s^m(i),j-1, \ldots, 1, s^m(i) +1, \ldots, N] \in \Xi_N 
$$
and then use \eqref{qkj}, \eqref{bracketykyj}. For this permutation, \thref{sym-Pprime} 
implies that $y_{\tau', s^n(j)-j+1}$ and $y_{\tau', s^m(i)-i+1}$ are 
equal to $y_{[j, s^n(j)]}$ and $y_{[i,s^m(i)]}$, 
respectively.
\end{proof}

\subsection{An identity for Poisson-normal elements}
\label{identPnormal}
For $1 \le j \le l \le N$, set 
\begin{equation}
\label{genO}
O_-^j(l) = \max \{ m \in \Zset_{\geq 0} \mid p^m(l) \geq j \}.
\end{equation}
The following fact follows directly from Theorems \ref{tmainPprime} and \ref{tsym-Pprime}. 

\bco{PprimeRjk} Let $R$ be a symmetric Poisson-CGL extension 
of length $N$ and $1 \leq j \leq k \leq N$. 
Then
$$
\Big\{ y_{[p^{O_-^j(i)}(i), i]} \Bigm| i \in [j,k] , \; s(i)>k \Big\}
$$
is a list of the homogeneous Poisson-prime elements of $R_{[j,k]}$ up to scalar multiples.
\eco

\bco{u-elem} Assume that $R$ is a symmetric Poisson-CGL extension of length $N$,
and $i \in [1,N]$ and $m \in \Zset_{>0}$ are such that $s^m(i) \in [1,N]$.
Then
\begin{equation}
\label{uuu}
u_{[i,s^m(i)]} := 
y_{[i, s^{m-1}(i)]} y_{[s(i), s^m(i)]} - y_{[s(i), s^{m-1}(i)]} y_{[i,s^m(i)]} 
\end{equation}
is a nonzero homogeneous Poisson-normal element of $R_{[i+1,s^m(i)-1]}$ 
which is not a multiple of $y_{[s(i), s^{m-1}(i)]}$ if $m \geq 2$. It Poisson-normalizes the elements of $R_{[i+1,s^m(i)-1]}$
in exactly the same way as $y_{[s(i), s^{m-1}(i)]} y_{[i,s^m(i)]}$ does.  Moreover,
\begin{equation}  \label{uformula}
u_{[i,s^m(i)]} = \psi \prod_{k\in P} y^{m_k}_{[p^{O_-^{i+1}(k)}(k),k]}
\end{equation}
where $\psi \in \kx$,
\begin{equation}  \label{Pset}
P = P_{[i,s^m(i)]} := \{ k\in [i,s^m(i)] \setminus \{i,s(i),\dots,s^m(i)\} \mid s(k) > s^m(i) \},
\end{equation}
and the integers $m_k$ are those from Theorems {\rm\ref{tswap2}, \ref{tswap3}}. Consequently, the leading term of $u_{[i,s^m(i)]}$ 
has the form 
$$
\psi x^{n_{i+1} e_{i+1} + \cdots + n_{s^m(i)-1} e_{s^m(i)-1}}
$$
for some $n_{i+1}, \ldots, n_{s^m(i)-1} \in \Znn$ 
such that $n_{s(i)}= \ldots = n_{s^{m-1}(i)} =0$
and $n_j = n_l$ for all $i+1 \leq j \leq l \leq s^m(i)-1$ with 
$\eta(j) = \eta(l)$.
\eco

\begin{proof}
Consider the P-CGL extension presentations \eqref{taupOre} of $R$ associated to the elements 
\begin{align*}
\tau_{i, s^m(i)-1} &= [i+1, \ldots, s^{m}(i)-1, i, s^m(i), s^m(i)+1, \ldots, N, i-1, \ldots, 1]   \\
\tau_{i, s^m(i)} &= [i+1, \ldots, s^{m}(i)-1, s^m(i), i, s^m(i)+1, \ldots, N, i-1, \ldots, 1] 
\end{align*}
of $\Xi_N$. Applied to these presentations, Theorems \ref{tswap2}, \ref{tswap3}, and \ref{tsym-Pprime} yield all but the last statement. That statement follows from \eqref{uformula}, due to the choice of leading terms for the elements $y_{[p^{O_-^{i+1}(k)}(k),k]}$.
\end{proof}

Since $y_\varnothing=1$ 
and $y_{[i,i]}=x_i$, the case $m=1$ of the corollary states that
\begin{equation}
\label{ucc'}
u_{[i,s(i)]}:= x_i x_{s(i)} - y_{[i,s(i)]}= c_{[i,s(i)-1]}= c'_{[i+1, s(i)]}
\end{equation}
is a nonzero homogeneous Poisson-normal element of $R_{[i+1,s(i)-1]}$, cf.~\thref{y-int}(d).
We set 
$$
u_{[i,i]} :=1.
$$

\bex{OMmn6}
In the case of $R = \OMmn$ (\exref{OMmn}), the result of \exref{OMmn5} shows how to express each $u_{[i, s^l(i)]}$ as a difference of products of solid minors. In particular, if $i = (r-1)n+c$ with $s(i) \in [1,N]$, then $r<m$, $c<n$, and
\begin{equation}  \label{OMmn-uisi}
\begin{aligned}
u_{[i,s(i)]} &:= x_i x_{s(i)} - y_{[i,s(i)]}= t_{rc} t_{r+1, c+1} - \Delta_{[r, r+1], [c, c+1]}  \\
 &= t_{r, c+1} t_{r+1, c} = x_{(r-1)n + c+1} x_{rn + c} .
\end{aligned}
\end{equation}
More generally, if $i = (r-1)n+c$ and $l \in \Znn$ with $s^l(i) \ne +\infty$, one can calculate that
$$
u_{[i, s^l(i)]} = \De_{ [r, r+l-1], [c+1, c+l] } \De_{ [r+1, r+l], [c, c+l-1] } .
$$
\eex

The next result, which is the main result in this subsection, 
describes relationships among the Poisson-normal elements in 
\coref{u-elem}. It will play a key role in normalizing 
the scalars $\kappa$ from \thref{swap2}. Recall \S\ref{PprimePCGL}
for the definition of leading terms $\lt(-)$.

\bth{u-prod} Let $R$ be a symmetric Poisson-CGL extension of length $N$.
For all $i \in [1,N]$ and $m \in \Zset_{>0}$ 
such that $s^{m+1}(i) \in [1,N]$,
\begin{equation}
\label{u-ident}
\lt( u_{[s(i),s^m(i)]} u_{[i,s^{m+1}(i)]})  = 
\lt( u_{[i,s^m(i)]} u_{[s(i),s^{m+1}(i)]} ).
\end{equation}
\eth

\begin{proof} Consider congruences modulo the ideal
$$
W_m := y_{[s(i),s^m(i)]} R_{[i, s^{m+1}(i)]}
$$
of $R_{[i, s^{m+1}(i)]}$. By the definition of the elements $u_{[*,*]}$, 
\begin{align*}
u_{[s(i),s^m(i)]} &\equiv y_{[s(i),s^{m-1}(i)]} y_{[s^2(i), s^m(i)]}  &u_{[i,s^{m+1}(i)]} &\equiv y_{[i, s^m(i)]} y_{[s(i), s^{m+1}(i)]}  \\
u_{[i,s^m(i)]} &\equiv - y_{[s(i),s^{m-1}(i)]} y_{[i,s^m(i)]} &u_{[s(i),s^{m+1}(i)]} &\equiv - y_{[s^2(i), s^m(i)]} y_{[s(i), s^{m+1}(i)]} 
\end{align*}
modulo $W_m$, whence
$$
u_{[s(i),s^m(i)]} u_{[i,s^{m+1}(i)]} \equiv u_{[i,s^m(i)]} u_{[s(i),s^{m+1}(i)]} \quad
\mod W_m.
$$
Consequently,
\begin{equation}
\label{3terms}
u_{[s(i), s^m(i)]} u_{[i,s^{m+1}(i)]}  = u_{[i,s^m(i)]} u_{[s(i),s^{m+1}(i)]}
+ y_{[s(i), s^m(i)]} r 
\end{equation}
for some $r \in R_{[i, s^{m+1}(i)]}$. Since we are done if $r=0$, we may assume that $r \ne 0$.

Note that the leading term of $y_{[s(i), s^m(i)]} r$ has the form $\xi x^f$ where $\xi \in \kx$ and $f = \sum_{l=i}^{s^{m+1}(i)} n_l e_l$ with $n_{s(i)} \ge 1$. On the other hand, the last statement of \coref{u-elem} implies that the leading terms of the other two products in \eqref{3terms} do not contain positive powers of $x_{s(i)}$. Hence, neither of these leading terms is a scalar multiple of $x^f$. The validity of \eqref{u-ident} follows.
\end{proof}

\sectionnew{Chains of mutations in symmetric Poisson-CGL extensions}
\label{mut-sym}

For a given Poisson-CGL extension $R$ one has the freedom of 
rescaling the generators $x_j$ by elements of the base field 
$\KK$. The prime elements $y_j$ from \thref{mainPprime}  obviously depend (again up to rescaling) 
on the choice of $x_j$. 
In this section we prove that for each symmetric P-CGL extension $R$ 
its generators can be rescaled in such a way that all
scalars $\kappa$ in \thref{swap2} become equal to one. 
This implies a mutation theorem, proved in Section \ref{main}, 
associated 
to the elements of $\Xi_N$ via the sequences of prime elements 
from \thref{sym-Pprime}.

Throughout this section we will assume that $R$ is a symmetric P-CGL extension of length $N$.

\subsection{The leading coefficients of the elements $u_{[i,s^m(i)]}$}
\label{lcoeffu}
Recall \eqref{uuu}.
For $i \in [1,N]$ and $m \in \Zset_{\ge0}$ such that $s^m(i) \in [1,N]$, 
let
$$
\pi_{[i,s^m(i)]} \in \kx \quad \mbox{and} 
\quad
f_{[i, s^m(i)]} \in \sum_{j =i+1}^{s^m(i)-1} \Zset_{\geq 0}\, e_j
\subset \Zset^N
$$
be given by
\begin{equation}
\label{pi-f}
\lt(u_{[i,s^m(i)]})= \pi_{[i,s^m(i)]} x^{f_{[i,s^m(i)]}}.
\end{equation}
Note that $\pi_{[i,i]} = 1$ and $f_{[i,i]} = 0$, because $u_{[i,i]} = 1$.

From \eqref{u-ident}, we obtain
\begin{equation}
\label{pifrel}
\begin{aligned}
\pi_{[s(i), s^m(i)]} \pi_{[i, s^{m+1}(i)]} &= \pi_{[i, s^m(i)]} \pi_{[s(i), s^{m+1}(i)]}  \\
f_{[s(i), s^m(i)]} + f_{[i, s^{m+1}(i)]} &= f_{[i, s^m(i)]} + f_{[s(i), s^{m+1}(i)]} 
\end{aligned}
\end{equation}
for $i \in [1,N]$ and $m \in \Zset_{>0}$ with $s^{m+1}(i) \in [1,N]$. A quick induction then shows that
\begin{equation}
\label{pifrel2}
\begin{aligned}
\pi_{[i,s^m(i)]} &= \pi_{[i, s(i)]} \cdots \pi_{[s^{m-1}(i), s^m(i)]}  \\
f_{[i,s^m(i)]} &= f_{[i, s(i)]} + \cdots + f_{[s^{m-1}(i), s^m(i)]}
\end{aligned}
\end{equation}
for $i \in [1,N]$ and $m \in \Zset_{>0}$ with $s^m(i) \in [1,N]$. Consequently, we have

\ble{pi=1}
If $\pi_{[i,s(i)]} = 1$ for all $i\in [1,N]$ such that $s(i) \ne +\infty$, then $\pi_{[i,s^m(i)]} = 1$ for all 
$i \in [1,N]$ and $m \in \Zset_{\ge0}$ such that $s^m(i) \in [1,N]$.
\qed\ele

\subsection{Rescaling of the generators of a symmetric Poisson-CGL extension}
\label{rescalexj}
For $\gab:=(\ga_1, \ldots, \ga_N)  \in (\kx)^N$ and $f:=(n_1, \ldots, n_N) \in \Zset^N$ set
$$
\gab^f := \ga_1^{n_1} \ldots \ga_N^{n_N} \in \kx.
$$

Given $\gab \in (\kx)^N$,
one can rescale the generators $x_j$ of $R$,
\begin{equation}
\label{rescale}
x_j \mt \ga_j x_j, \; \; \forall j \in [1,N],
\end{equation}
by which we mean that one can use $\ga_1 x_1, \ldots, \ga_N x_N$ 
as a new sequence of generators
 and express $R$ as a P-CGL extension of the form
$$
R = \KK[\ga_1 x_1]_p [\ga_2 x_2; \sig_2, \ga_2 \delta_2]_p \cdots [\ga_N x_N; \sig_N, \ga_N \delta_N]_p .
$$
This change of generators obviously does not effect the $\HH$-action 
and the matrix $\lab$, but one obtains a new set of elements $y_k$, $y_{[i,s^m(i)]}$,
$u_{[i, s^m(i)]}$ by applying Theorems \ref{tmainPprime}, \ref{ty-int} and \coref{u-elem} for the new set 
of generators. (Note that this is not the same as substituting \eqref{rescale} 
in the formulas for $y_k$ and $y_{[i,s^m(i)]}$; those other kind of transformed elements 
may not be even prime because \eqref{rescale} does not determine an algebra automorphism.)
The uniqueness part of \thref{y-int} implies that 
the effect of \eqref{rescale} on the elements $y_{[i, s^m(i)]}$ 
is that they are rescaled by the rule
$$
y_{[i,s^m(i)]} \mt \gab^{e_{[i, s^m(i)]}} y_{[i,s^m(i)]}=(\ga_i \ga_{s(i)} \ldots \ga_{s^m(i)}) y_{[i,s^m(i)]}
$$
for all $i \in [1,N]$, $m \in \Zset_{\geq 0}$ such that $s^m(i) \in [1,N]$.
Hence, the effect of \eqref{rescale} on the elements $u_{[i,s^m(i)]}$ is that they are rescaled by
$$
u_{[i,s^m(i)]} \mt 
\gab^{e_{[i, s^m(i)]} + e_{[s(i), s^{m-1}(i)]}} u_{[i,s^m(i)]}= 
(\ga_i \ga_{s(i)}^2 \ldots \ga_{s^{m-1}(i)}^2 \ga_{s^m(i)}) u_{[i, s^m(i)]} .
$$
It follows from \eqref{pi-f} that
the effect of \eqref{rescale} on the scalars $\pi_{[i,s^m(i)]}$ is that 
they are rescaled by 
\begin{multline*}
\pi_{[i,s^m(i)]} \mt 
\gab^{e_{[i, s^m(i)]} + e_{[s(i), s^{m-1}(i)]} - f_{[i,s^m(i)]}} \pi_{[i,s^m(i)]}=  \\
\ga_i \ga_{s^m(i)} \gab^{2 e_{[s(i), s^{m-1}(i)]} - f_{[i,s^m(i)]}} \pi_{[i, s^m(i)]}.
\end{multline*}
(Note that the rescaling \eqref{rescale} has no effect on the integer vector $f_{[i,s^m(i)]}$.)
This implies at once the following fact.

\bpr{resc} Let $R$ be a symmetric CGL extension of length $N$ and rank $\rk(R)$. 
Then there exist $N$-tuples $\gab \in (\kx)^N$ such that after the rescaling 
\eqref{rescale} we have 
\begin{equation}
\label{pi-cond}
\pi_{[i,s(i)]} =1, \; \;   \forall i \in [1,N] \; \; \mbox{such that} \; \; 
s(i) \ne +\infty.
\end{equation}
The set of those $N$-tuples is parametrized by $(\kx)^{\rk(R)}$ and 
the coordinates $\ga_1, \ldots, \ga_N$ of all such $N$-tuples 
are recursively determined by
$$
\ga_i \; \; \mbox{is arbitrary if} \; \; p(i) = - \infty 
$$
and 
\begin{equation}
\label{ga-recursive}
\ga_i = \ga_{p(i)}^{-1} \gab^{f_{[p(i),i]}} \pi_{[p(i), i]}^{-1}, \; \; \text{if} \; \; p(i) \ne -\infty,
\end{equation}
where on the right hand side the $\pi$-scalars are the ones 
for the original generators $x_1, \ldots, x_N$ of $R$.
\epr

Note that the product of the first two terms of the right hand side of \eqref{ga-recursive} is a product of 
powers of $\ga_{p(i)}, \ldots, \ga_{i-1}$ since 
$f_{[p(i), i]} \in \Zset_{\geq 0} e_{p(i)+1} + \cdots + \Zset_{\geq 0} e_{i-1}$.  

\bex{OMmn7}
In the case of $R = \OMmn$, \eqref{OMmn-uisi} shows that the elements $u_{[i, s(i)]}$ all have leading coefficient $1$, and so no rescaling is needed in this case.
\eex

\subsection{Mutations in symmetric Poisson-CGL extensions}
\label{mutaPCGL}
Continue to let $R$ be a symmetric P-CGL extension of length $N$. Let $1\le j \le l \le N$, and recall the definition \eqref{genO} of the function $O_-^j$. 
By \coref{PprimeRjk},
$$
\Bigl\{ y_{[p^{O_-^j(k)}(k), k]} \Bigm| k \in [j,l], \; s(k) >l \Bigr\}
$$
is a list of the the homogeneous prime elements of $R_{[j,l]}$ 
up to scalar multiples.

Let $i \in [1,N]$ and $m \in \Zset_{> 0}$ be such that $s^m(i) \in [1,N]$, and consider the P-CGL extension presentation 
of $R_{[i,s^m(i)]}$ associated to the following order of adjoining 
the generators:
\begin{equation}
\label{gene-int}
x_{i+1}, \ldots, x_{s^m(i)-1}, x_i, x_{s^m(i)}.
\end{equation}
The intermediate subalgebras for this presentation of $R_{[i,s^m(i)]}$ are 
$$
R_{[i+1, i+1]}, \ldots, R_{[i+1, s^m(i)-1]}, R_{[i,s^m(i)-1]}, R_{[i, s^m(i)]}.
$$
We identify $\Zset^{s^m(i)-i+1} \cong \Zset e_i + \cdots + \Zset e_{s^m(i)} \subseteq \Zset^N$, 
and we define a homomorphism 
$$
Y_{[i,s^m(i)]} : \Zset^{s^m(i)-i+1} \to \bigl( \Fract(R_{[i,s^m(i)]}) \bigr)^*
$$
by setting
\begin{equation}
\label{Yismi1}
Y_{[i,s^m(i)]}(e_k) =
\begin{cases}
y_{[p^{O_-^{i+1}(k)}(k),k]}, & \mbox{if} \; \; k \in [i+1, s^m(i)-1]
\\
y_{[i,s^{m-1}(i)]}, & \mbox{if} \; \; k = i
\\
y_{[i, s^m(i)]}, & \mbox{if} \; \; k = s^m(i) .
\end{cases}
\end{equation}

Recall the definition \eqref{pi-f} of the vectors $f_{[i,s^m(i)]}$. It follows 
from \coref{u-elem} that there exists 
a unique vector 
\begin{equation}  \label{gismi}
g_{[i, s^m(i)]} = \sum \{ m_k e_k \mid k \in P_{[i,s^m(i)]} \} \in \Zset^{s^m(i)-i+1} 
\end{equation}
such that 
$$
f_{[i,s^m(i)]} = \sum \Bigl\{ m_k e_{[p^{O_-^{i+1}(k)}(k), k]} \Bigm| 
k \in P_{[i,s^m(i)]} \Bigr\},
$$
where the integers $m_k$ are those from Theorems \ref{tswap2}, \ref{tswap3}.

Set $t := s^m(i)-i+1$. As above, we identify $\Zset^t$ with $\Zset e_i + \cdots + \Zset e_{s^m(i)}$, that is, with the sublattice of $\Zset^N$ 
with basis $e_j$, $j \in [i,s^m(i)]$. 
Define $\sig \in M_t(\Zset)$ by
\begin{equation}
\label{ekprime}
\sig(e_k) :=
\begin{cases}
e_{[p^{O_-^{i+1}(k)}(k),k]}, & \mbox{if} \; \; k \in [i+1, s^m(i)-1]
\\
e_{[i,s^{m-1}(i)]}, & \mbox{if} \; \; k = i
\\
e_{[i, s^m(i)]}, & \mbox{if} \; \; k = s^m(i).
\end{cases}  
\end{equation}
 for all $k \in [i,s^m(i)]$. It is easy to check that 
$$
\sig \; \; {\mbox{is invertible and}} \; \; \sig(\Zset_{\geq 0}^t) \subseteq \Zset_{\geq 0}^t.
$$
This choice of $\sig$ ensures that
\begin{equation}
\label{gfsig}      
\sig(g_{[i, s^m(i)]}) = f_{[i, s^m(i)]} .
\end{equation}                  

The definition of $Y_{[i, s^m(i)]}$ also implies that 
\begin{equation}
\label{ltY}
\lt \left( Y_{[i, s^m(i)]}(e_k) \right) = x^{\sig(e_k)}, \; \;
\forall k \in [i, s^m(i)].
\end{equation}

\ble{leadterm} For all $g \in \Zset_{\geq 0}^t$,
$$
\lt \left( Y_{[i, s^m(i)]}(g) \right) = x^{\sig(g)} .
$$
\ele

\begin{proof} This is immediate from \eqref{ltY}, since $Y_{[i, s^m(i)]}$ and $\sigma$ are homomorphisms.
\end{proof}

\bth{PCGLmuta} Let $R$ be a symmetric Poisson-CGL extension of length $N$.
Assume that
the generators $x_1, \ldots, x_N$ of $R$ are normalized so that 
the condition \eqref{pi-cond} is satisfied. Then for all 
$i \in [1,N]$ and $m \in \Zset_{>0}$ such that $s^m(i) \in [1,N]$,
\begin{equation}
\label{ismi.mut}
y_{[s(i), s^m(i)]} = \begin{cases}
Y_{[i,s(i)]}(e_{s(i)} - e_i) + Y_{[i,s(i)]}(g_{[i,s(i)]} - e_i),  &\text{if} \; \; m=1\\
Y_{[i,s^m(i)]} (e_{s^{m-1}(i)} + e_{s^m(i)} - e_i) + Y_{[i,s^m(i)]} ( g_{[i,s^m(i)]} - e_i),  &\text{if} \; \; m>1. \end{cases}
\end{equation}
\eth

\begin{proof} \leref{leadterm} and Eq. \eqref{gfsig} imply 
$$
\lt \left( Y_{[i, s^m(i)]}(g_{[i, s^m(i)]} ) \right) = x^{f_{[i, s^m(i)]}}.
$$
By Eq. \eqref{uformula}, $u_{[i, s^m(i)]}$ is a scalar multiple of $Y_{[i, s^m(i)]}(g_{[i, s^m(i)]})$, and so we find that
\begin{equation}
\label{qu-scalar}
u_{[i,s^m(i)]} = Y_{[i,s^m(i)]}(g_{[i,s^m(i)]}),
\end{equation}
taking account of \eqref{pi-cond} and \leref{pi=1}.

In case $m>1$, we observe using \eqref{qu-scalar} that
\begin{align*}
y_{[s(i),s^m(i)]} &- Y_{[i, s^m(i)]}(e_{s^{m-1}(i)}+ e_{s^m(i)}- e_i) \\
 &=  y_{[s(i),s^m(i)]} - y^{-1}_{[i,s^{m-1}(i)]} y_{[s(i), s^{m-1}(i)]} y_{[i, s^m(i)]}  \\
 &= y^{-1}_{[i,s^{m-1}(i)]} u_{[i, s^m(i)]} = Y_{[i, s^m(i)]}( g_{[i, s^m(i)]} - e_i ) .
\end{align*}
A similar argument shows that
$$
y_{[s(i),s(i)]} - Y_{[i,s(i)]}(e_{s(i)} - e_i) = x_i^{-1} u_{[i,s(i)]} = Y_{[i,s(i)]}( g_{[i,s(i)]} - e_i ) 
$$
 in the case $m = 1$.
\end{proof}

\bre{appl} Because of \thref{y-int}, all instances when \thref{swap3} is applicable 
to a symmetric P-CGL extension with respect to the P-CGL extension presentations 
associated to the elements of the set $\Xi_N$ are covered by \thref{PCGLmuta}. 
We refer the reader to \prref{cluster-tau-ind} for details.
\ere

\sectionnew{Division properties of mutations between Poisson-CGL extension presentations}
\label{Integr}

We prove in this section that every symmetric Poisson-CGL extension is equal to the intersection of appropriate localizations with respect to sets of $y$-elements.
This plays a key role in Section \ref{main}
where we prove that every symmetric P-CGL extension is a cluster algebra which 
equals the corresponding upper cluster algebra. In the process, we introduce auxiliary permutations $\tau_\bu$, for $\tau \in \Xi_N$, which are used to re-order our sets of cluster variables.

\subsection{One-step mutations}
\label{5.1}
In this subsection we describe the intersection of the localizations of a 
P-CGL extension $R$ by the two sets of  $y$-elements  
in the setting of \thref{swapkk+1}(b).

Assume now that $R$ is a P-CGL extension of length $N$ as in \eqref{itpOre}.
Let $y_1, \ldots, y_N \in R$ be the sequence of 
elements from \thref{mainPprime}. For all $I \subseteq [1,N]$, the set
\begin{equation}
\label{Sloc}
E_I:=
\kx \biggl\{ \prod_{j \in I} y_j^{m_j} \biggm| m_j \in \Zset_{\geq 0} \biggr\}
\end{equation}
is a multiplicative set in $R$. Set $E:= E_{[1,N]}$.
We will say that $y_1^{m_1} \ldots y_N^{m_N}$ is a
\emph{minimal denominator} of a nonzero element $v \in R[E^{-1}]$ if 
$$
v = y_1^{-m_1} \ldots y_N^{-m_N} s
$$
for some $s \in R$ such that $y_j \nmid s$ for all $j \in [1,N]$ 
with $m_j > 0$.

For the rest of this subsection, we assume the setting of \thref{swapkk+1}(b). 
Then for all $I' \subseteq [1,N]$, the set
\begin{equation}
\label{E'set}
E'_{I'}:= \kx \biggl\{ \prod_{j \in I'} (y'_j)^{m_j} \biggm| m_j \in \Zset_{\geq 0} \biggr\}
\end{equation}
is another multiplicative set in $R$. Recall from \thref{swapkk+1}(b) that $y'_j = y_j$ for $j \neq k$. Hence, $E'_{I'} = E_{I'}$ if $k \notin I'$. Denote $E':= E'_{[1,N]}$.   

\ble{nondivide}
For any distinct $j,l \in [1,N]$, we have $y_j \nmid y_l$ and $y_l \nmid y_j$. Moreover, $y_j \nmid y'_k$ and $y'_k \nmid y_j$ for all $j \in [1,N]$.
\ele

\begin{proof} If $1 \le j < l \le N$, then $y_l \nmid y_j$ because $y_l$ and $y_j$ are elements of $R_l$ with $x_l$-degrees $1$ and $0$ respectively. Since $y_l$ and $y_j$ are prime elements of $R$, it follows that $y_j \nmid y_l$ as well.

On replacing the $y_j$ by the $y'_j$, we have $y'_j \nmid y'_l$ and $y'_l \nmid y'_j$ for all $j \ne l$, which yields $y_j \nmid y'_k$ and $y'_k \nmid y_j$ for all $j \ne k$. Since $y_k$ and $y'_k$ are elements of $R_{k+1}$ with $x_{k+1}$-degrees $0$ and $1$ respectively, we conclude as above that $y'_k \nmid y_k$ and $y_k \nmid y'_k$.
\end{proof}

\ble{EIEI'} For all $I, I' \subseteq [1,N]$,
\begin{equation}
\label{EcapE}
E_I \cap E'_{I'} = 
E_{(I \cap I') \backslash \{ k \}} = E'_{(I \cap I') \backslash \{ k \} }.
\end{equation}
\ele

\begin{proof} Let $e \in E_I \cap E'_{I'}$, and write
$$
e = \be \prod_{j=1}^N y_j^{m_j} = \be' \prod_{j=1}^N (y'_j)^{m'_j}
$$
for some $\be, \be' \in \kx$ and $m_j, m'_j \in \Znn$ such that $m_j=0$ for $j \notin I$ and $m'_j=0$ for $j \notin I'$. Since the $y_j$ and $y'_j$ are prime elements of $R$, with $y_j = y'_j$ for all $j \ne k$, it follows from \leref{nondivide} that
$$
m_k = m'_k = 0 \quad \text{and} \quad m_j = m'_j, \; \; \forall j \in [1,N], \; j \ne k,
$$
which then implies $e \in E_{(I \cap I') \backslash \{ k \}} = E'_{(I \cap I') \backslash \{ k \} }$.
\end{proof}

The next theorem  implies that each nonzero element of $R[E^{-1}]$ has a unique minimal 
denominator. It was inspired by \cite[Theorem 4.1]{GLSh2}.

\bth{division} Assume the setting of Theorem {\rm\ref{tswapkk+1}(b)}. Let 
$$
y_1^{m_1} \ldots y_N^{m_N} \quad
\mbox{and} \quad 
y_1^{m'_1} \ldots y_{k-1}^{m'_{k-1}} (y'_k)^{m'_k} y_{k+1}^{m'_{k+1}} \ldots y_N^{m'_N}
$$
be two minimal denominators of a nonzero element $v \in R[E^{-1}] \cap R[(E')^{-1}]$
{\rm(}with respect to the two different localizations\/{\rm)}. Then
$$
m_k = m'_k =0 \quad \mbox{and} 
\quad
m_j= m'_j, \; \; \forall j \in [1,N], \; j \neq k.
$$ 
In particular, 
$$
R[(E_I)^{-1}] \cap R[(E'_{I'})^{-1}] = R[ (E_I \cap E'_{I'})^{-1}] =
R[(E_{(I\cap I')\backslash\{k \}})^{-1}]
$$
for all $I,I' \subseteq [1,N]$.
\eth

\begin{proof} Write
$$
v = y_1^{-m_1} \ldots y_N^{-m_N} s = 
y_1^{-m'_1} \ldots y_{k-1}^{-m'_{k-1}} (y'_k)^{-m'_k} y_{k+1}^{-m'_{k+1}} \ldots y_N^{- m'_N}
s'   
$$
for some nonzero elements $s, s' \in R$. Then
\begin{equation}
\label{ys'y's}
y_1^{m_1} \ldots y_N^{m_N} s' = y_1^{m'_1} \ldots y_{k-1}^{m'_{k-1}} (y'_k)^{m'_k} y_{k+1}^{m'_{k+1}} \ldots y_N^{m'_N} s.
\end{equation}
If $j \in [1,N] \backslash \{k\}$ and $m_j > 0$, then $y_j \nmid s$, so it follows from \leref{nondivide} and the primeness of $y_j$ that $m_j \le m'_j$. The reverse inequality holds by symmetry, whence $m_j = m'_j$.

Equation \eqref{ys'y's} now reduces to $y_k^{m_k} s' = (y'_k)^{m'_k} s$. Since $y_k \nmid s$ if $m_k > 0$, and $y_k \nmid y'_k$ (by \leref{nondivide}), we must have $m_k = 0$, and similarly $m'_k = 0$.

The final statement of the theorem follows.
\end{proof}

\subsection{Permutations $\tau_\bu$}  \label{taubull}
If $R$ is a symmetric P-CGL extension, then for each $\tau \in \Xi_N$ we have a P-CGL extension presentation \eqref{taupOre} for $R$ and corresponding $y$-elements by applying \thref{mainPprime} to this presentation. We label these $y$-elements $y_{\tau,1}, \dots, y_{\tau,N}$; they will become a set of cluster and frozen variables in our main theorem. In order to connect these various sets of variables via mutations, we need to permute the order of each set of variables. This is done by applying a composition of $\tau$ 
and a permutation $\tau_\bu \in \prod_{a \in \Zset} S_{\eta^{-1}(a)}$,
where $\eta \colon [1,N] \to \Zset$ is a function satisfying the 
conditions of \thref{mainPprime} for the original P-CGL extension presentation of $R$.
Note that all terms in the above product are trivial except for the 
terms coming from the range of $\eta$.
For $a$ in the range of $\eta$, denote for brevity $|a|:= |\eta^{-1}(a)|$. Write the set $\eta^{-1}(a)$ in the form
\begin{equation}
\label{eta-1}
\eta^{-1}(a)= \{ \tau(k_1), \ldots, \tau(k_{|a|}) \mid k_1 < \cdots  < k_{|a|} \}
\end{equation}
and order its elements in ascending order. Define $\tau_\bu \in S_N$ by setting
$\tau_\bu(\tau(k_i))$ to be equal to the $i$-th element in the list 
(for all choices of $a$ and $i$). In other words,
\begin{enumerate}
\item[$(\tau_\bu)$] For $l \in [1,N]$, let $\eta^{-1} \eta(l) = \{ t_1 < \cdots < t_m\}$ and
$$i := \big| \tau^{-1} \eta^{-1} \eta(l) \cap [1, \tau^{-1}(l)] \big|.$$
Then $\tau_\bu(l) := t_i$.
\end{enumerate}
Thus, for each level set $L$ of $\eta$, the permutation $\tau_\bu \tau$ gives an order-preserving bijection of $\tau^{-1}(L)$ onto $L$.

The point of applying $\tau$ is
to match the indexing of the $y$-elements with that of the 
$x$-elements in \eqref{taupOre}. (Note that the 
order of the $x$-elements in \eqref{taupOre} is 
$x_{\tau(1)}, \ldots, x_{\tau(N)}$.) The application 
of $\tau_\bu$ then rearranges the $\eta$-preimages 
$\tau(k_1), \ldots, \tau(k_{|a|})$ from \eqref{eta-1}
in ascending order. This is needed because in the setting
of \thref{swapkk+1}(b) the element $y_k$ (not $y_{k+1}$) 
gets mutated. Clearly, $\tau_\bu$ preserves the level sets of 
$\eta$.

Recall that $P(N) = \{ k \in [1,N] \mid s(k) = + \infty \}$
parametrizes the set of homogeneous Poisson-prime elements of $R$,
i.e.,  
\begin{equation}
\label{RPr}
\{y_k \mid k \in P(N) \} \; \; \mbox{is a list of the homogeneous Poisson-prime
elements of $R$}
\end{equation}
up to associates. Define
\begin{equation}
\label{ex-prim}
\ex := [1,N] \setminus P(N) = \{ k \in [1,N] \mid s(k) \neq + \infty \}.
\end{equation}
Since $|P(N)| = \rk(R)$,
the cardinality of this set is $|\ex| = N - \rk(R)$. 

Some properties of the permutations $\tau_\bu \tau$ are given in the following lemma \cite[Lemma 8.16]{GYbig}.

\ble{taubutau} {\rm (a)} For any $\tau \in \Xi_N$, the permutation $(\tau_\bu \tau)^{-1}$ maps $\ex$ bijectively onto the set $\{ k \in [1,N] \mid s_\tau(k) \ne +\infty \}$.

{\rm (b)} Suppose $\tau, \tau' \in \Xi_N$ and $\tau' = \tau (k,k+1)$ for some $k \in [1,N-1]$. If $\eta \tau(k) = \eta \tau(k+1)$, then $\tau'_\bu \tau' = \tau_\bu \tau$.

{\rm (c)} Suppose $\tau, \tau' \in \Xi_N$ and $\tau' = \tau (k,k+1)$ for some $k \in [1,N-1]$. If $\eta \tau(k) \ne \eta \tau(k+1)$, then $\tau'_\bu \tau' = \tau_\bu \tau(k,k+1)$.

{\rm (d)} Let $j \in \ex$. Then there exist $\tau, \tau' \in \Ga_N$ such that $\tau' = \tau (k,k+1)$ for some $k \in [1,N-1]$ with $\eta \tau(k) = \eta \tau(k+1)$ and $\tau_\bu \tau(k) = j$.
\ele

\subsection{Intersections of localizations}  \label{intersectloc}
Throughout this subsection, we will assume that $R$ is a symmetric P-CGL extension of rank $N$ 
as in \deref{symmPCGL}, and we fix a function $\eta : [1,N] \to \Zset$ satisfying the conditions of \thref{mainPprime}. For each $\tau \in \Xi_N$, there is a P-CGL presentation
$$
R = \KK [x_{\tau(1)}] [x_{\tau(2)}; \sig''_{\tau(2)}, \de''_{\tau(2)}] 
\cdots [x_{\tau(N)}; \sig''_{\tau(N)}, \de''_{\tau(N)}]
$$
as in \eqref{taupOre}. Let $p_\tau$ and $s_\tau$ denote the predecessor and successor functions for the level sets of $\eta_\tau := \eta \tau$, which by \coref{steps}(b) can be chosen as the $\eta$-function for the presentation \eqref{taupOre}. Let $y_{\tau,1}, \dots, y_{\tau, N}$ be the corresponding sequence of homogeneous Poisson-prime elements of $R$ from \thref{mainPprime}, and denote
\begin{align*}
\AA_\tau &:= \text{the} \; \KK\text{-subalgebra of} \; R \; \text{generated by} \; \{ y_{\tau,k} \mid k \in [1,N] \}  \\
\TT_\tau &:= \text{the} \; \KK\text{-subalgebra of} \; \Fract(R) \; \text{generated by} \; \{ y_{\tau,k}^{\pm1} \mid k \in [1,N] \}  \\
E_\tau &:= \text{the multiplicative subset of} \; \AA_\tau \; \text{generated by}  \\
 & \qquad\qquad\qquad \qquad\qquad\qquad \kx \sqcup \{ y_{\tau,k} \mid k \in [1,N], \; s_\tau(k) \ne +\infty \}.
\end{align*}
By \prref{PtorusR}, $\TT_\tau$ is a Poisson torus with corresponding Poisson affine space algebra $\AA_\tau$.

Recall the subset $\Ga_N$ of $\Xi_N$ and its total ordering from \eqref{sequence}.

\bth{PCGLalmostcluster}
Let $R$ be a symmetric Poisson-CGL extension of length $N$.

{\rm (a)} $\AA_\tau \subseteq R \subseteq \AA_\tau[E_\tau^{-1}] \subseteq \TT_\tau \subseteq \Fract(R)$, for all $\tau \in \Xi_N$.

{\rm (b)} $R$ is generated as a $\KK$-algebra by $\{ y_{\tau,k} \mid \tau \in \Ga_N, \; k \in [1,N] \}$.

{\rm (c)} $R[E_\tau^{-1}] = \AA_\tau[E_\tau^{-1}]$, for all $\tau \in \Xi_N$.

{\rm (d)} $R = \bigcap_{\tau \in \Ga_N} R[E_\tau^{-1}] = \bigcap_{\tau \in \Ga_N} \AA_\tau[E_\tau^{-1}]$.

{\rm (e)} Let $\inv$ be any subset of $\{ k \in [1,N] \mid s(k) = +\infty \}$. Then
$$
R[y_k^{-1} \mid k \in \inv ] = \bigcap_{\tau \in \Ga_N} R[E_\tau^{-1}] [y_k^{-1} \mid k \in \inv ] .
$$
\eth

\begin{proof} (a) The first, third, and fourth inclusions are clear. For the second, it suffices to show that $x_{\tau(k)} \in \AA_\tau[E_\tau^{-1}]$ for all $k \in [1,N]$, which we do by induction on $k$. The case $k = 1$ is immediate from the fact that $x_{\tau(1)} = y_{\tau,1}$.

Now let $k \in [2,N]$. If $p_\tau(k) = -\infty$, then $x_{\tau(k)} = y_{\tau,k} \in \AA_\tau[E_\tau^{-1}]$. If $p_\tau(k) = l \ne -\infty$, then $l < k$ and $y_{\tau,k} = y_{\tau,l} x_{\tau(k)} - c_{\tau,k}$ for some element $c_{\tau,k}$ in the $\KK$-subalgebra of $R$ generated by $x_{\tau(1)}, \dots, x_{\tau(k-1)}$. By induction, $c_{\tau, k} \in \AA_\tau[E_\tau^{-1}]$. Further, $y_{\tau, l} \in E_\tau$ because $s_\tau(l) = k \ne +\infty$, and thus
$$
x_{\tau(k)} = y_{\tau,l}^{-1} \bigl( y_{\tau, k} + c_{\tau, k} \bigr) \in \AA_\tau[E_\tau^{-1}]
$$
in this case also.

(b) For each $j \in [1,N]$, we have $\tau_{j,j} \in \Ga_N$ with $\tau_{j,j}(1) = j$, and so $y_{\tau_{j,j},1} = x_j$. Thus, in fact, $R$ is generated by $\{ y_{\tau, 1} \mid \tau \in \Ga_N \}$.

(c) This is clear from part (a).

(d) The second equality follows from part (c), and one inclusion of the first equality is obvious.

Let $v \in \bigcap_{\tau \in \Ga_N} R[E_\tau^{-1}]$ be a nonzero 
element. For each $\tau \in \Ga_N$, let 
$$
\prod_{l \in [1,N]} y_{\tau, l}^{m_{\tau, l}} 
$$
be a minimal denominator of $v$ with respect to the localization $R[E_\tau^{-1}]$, where all $m_{\tau, l} \in \Znn$ and $m_{\tau, l} = 0$ when $s_\tau(l) = +\infty$. We first verify the following

{\bf Claim}. Let $\tau, \tau' \in \Ga_N$ such that $\tau' = \tau (k, k+1)$ for some $k \in [1, N-1]$.
\begin{enumerate}
\item $m_{\tau, (\tau_\bu \tau)^{-1}(j)} = m_{\tau', (\tau'_\bu \tau')^{-1}(j)}$, for all $j \in [1,N]$.
\item If $\eta \tau(k)  = \eta \tau(k+1)$, then $m_{\tau, k} = 0$.
\end{enumerate}

If $\eta \tau(k) = \eta \tau(k+1)$, then (1) and (2) follow from \leref{taubutau}(b), \thref{division}, and \thref{swapkk+1}(b). If $\eta \tau(k) \ne \eta \tau(k+1)$, we obtain from \thref{swapkk+1}(a) that $y_{\tau', k} = y_{\tau, k+1}$ and $y_{\tau', k+1} = y_{\tau, k}$, while $y_{\tau', l} = y_{\tau, l}$ for all $l \ne k, k+1$. As a result, $E_{\tau'} = E_\tau$ and we see that $m_{\tau', k} = m_{\tau, k+1}$ and $m_{\tau', k+1} = m_{\tau, k}$, while $m_{\tau', l} = m_{\tau, l}$ for all $l \ne k, k+1$. In this case, (1) follows from \leref{taubutau}(c).

Since all the permutations in $\Ga_N$ appear in the chain \eqref{sequence}, part (1) of the claim implies that 
\begin{equation}
\label{matchm}
m_{\sig, (\sig_\bu \sig)^{-1}(j)} = m_{\tau, (\tau_\bu \tau)^{-1}(j)}, \; \; \forall \sig, \tau \in \Ga_N, \; j \in [1,N].
\end{equation}
For any $j \in \ex$, \leref{taubutau}(d) shows that there exist $\tau, \tau' \in \Ga_N$ such that $\tau' = \tau (k,k+1)$ for some $k \in [1,N-1]$ with $\eta \tau(k) = \eta \tau(k+1)$ and $\tau_\bu \tau(k) = j$. Part (2) of the claim above then implies that $m_{\tau, (\tau_\bu \tau)^{-1}}(j) = 0$. From \eqref{matchm}, we thus get
$$
m_{\sig, (\sig_\bu \sig)^{-1}(j)} = 0, \; \; \forall \sig \in \Ga_N, \; j \in \ex.
$$
In particular, $m_{\id, j} = 0$ for all $j \in \ex$, whence $m_{\id, j} = 0$ for all $j \in [1,N]$ and therefore $v \in R$, which completes the proof of part (d).

(e) Let $v \in \bigcap_{\tau \in \Ga_N} R[E_\tau^{-1}] [y_k^{-1} \mid k \in \inv ]$. For each $\tau \in \Ga_N$, we can write $v$ as a fraction with numerator from $R[E_\tau^{-1}]$ and denominator from the multiplicative set $Y$ generated by $\{ y_k \mid k \in \inv \}$. Hence, choosing a common denominator, we obtain $y \in Y$ such that $vy \in R[E_\tau^{-1}]$ for all $\tau \in \Ga_N$. From part (d), we conclude that $vy \in R$, and thus $v \in R[y_k^{-1} \mid k \in \inv ]$, as required.
\end{proof}

\sectionnew{Symmetric Poisson-CGL extensions and cluster algebras}
\label{main}
In this section we prove that every symmetric Poisson-CGL extension
possesses a cluster algebra structure under a mild additional
assumption on the scalars $\la^*_k$. The cluster algebra structure is constructed
in an explicit fashion. 
We furthermore prove that each of these 
cluster algebras equals the corresponding upper cluster algebra.

Concerning notation: We continue to use $x_1,\dots,x_N$ to denote polynomial variables in P-CGL extensions, and our cluster and frozen variables will be denoted $y_\iota$ for various indices $\iota$.

\subsection{Statement of the main result}
\label{statement}
Fix a symmetric P-CGL extension $R$ of length $N$ and rank $\rk(R)$ 
as in \S \ref{symmP-CGL}. Recall the map $Y : \Zset^N \rightarrow \Fract(R)$ from \eqref{YY'def}, and note that $Y$ can be viewed as a group homomorphism from $\Zset^N$ to $\Fract(R)^*$. Set
\begin{equation}
\label{ybtil}
\ybtil := (Y(e_1), \dots, Y(e_N)) = (y_1, \dots, y_N) \in \Fract(R)^N,
\end{equation}
and recall from \prref{PtorusR} that $\ybtil$ is a sequence of algebraically independent elements generating $\Fract(R)$ over $\KK$. Recall also the skew-symmetric matrices $\lab, \qb \in M_N(\KK)$ from \deref{PCGL} and \eqref{qkj}, and the corresponding skew-symmetric bicharacters $\Om_\lab$, $\Om_\qb$ from \deref{PCGL} and \eqref{Ombfq}.

Consider an arbitrary element $\tau \in \Xi_N \subset S_N$, recall \eqref{tau}, and the associated auxiliary permutation $\tau_\bu$ from \S\ref{taubull}. By \prref{tauPCGL}, we have the P-CGL extension presentation
\begin{equation}
\label{tau-pres}
R = \KK [x_{\tau(1)}]_p [x_{\tau(2)}; \sig''_{\tau(2)}, \de''_{\tau(2)}]_p 
\cdots [x_{\tau(N)}; \sig''_{\tau(N)}, \de''_{\tau(N)}]_p
\end{equation}
with $\sig''_{\tau(k)} = ( h''_{\tau(k)} \cdot)$ where 
$h''_{\tau(k)} \in \HH$, $\forall k \in [1,N]$ and  
\begin{align*}
&h''_{\tau(k)} = h_{\tau(k)}, \; \; \de''_{\tau(k)} = \de_{\tau(k)}, 
\quad \mbox{if} \; \;  \tau(k) = \max \, \tau( [1,k-1]) +1
\\  
&h''_{\tau(k)} = h^\sy_{\tau(k)}, \; \; \de''_{\tau(k)} = \de^*_{\tau(k)},
\quad \mbox{if} \; \; \tau(k) = \min \, \tau( [1,k-1]) -1.
\end{align*}
The $\lab$-matrix of the presentation \eqref{tau-pres} is $\lab_\tau:=\tau^{-1} \lab \tau$
(where as before we use the canonical embedding $S_N \hra GL_N(\Zset)$ 
via permutation matrices). In other words,
the entries of $\lab_\tau$ are given by $(\lab_\tau)_{lj}:= \la_{\tau(l), \tau(j)}$.
Let $y_{\tau,1}, \ldots, y_{\tau,N}$ 
be the sequence of homogeneous Poisson-prime elements from \thref{mainPprime} for the presentation \eqref{tau-pres}.

Let $\qb_\tau \in M_N(\KK)$ be the skew-symmetric matrix with entries 
$$
(\qb_\tau)_{kj} := \sum_{l=0}^{O_-^\tau(k)} \sum_{m=0}^{O_-^\tau(j)} (\lab_\tau)_{p_\tau^l(k), p^m_\tau(j)} , \; \; \forall k,j \in [1,N],
$$
where $O_-^\tau$ denotes the order function corresponding to $p_\tau$,
and let $\Om_{\qb_\tau}$ be the skew-symmetric bicharacter on $\Zset^N$ obtained from $\qb_\tau$ as in \eqref{Om}. Then $\Om_{\qb_\tau}$ is the skew-symmetric bicharacter corresponding to $\Om_\qb$ (recall \eqref{Ombfq}) for the presentation \eqref{tau-pres}. Set
$$
\rb_\tau := (\tau_\bu \tau) \qb_\tau (\tau_\bu \tau)^{-1},
$$
and let $\Om_{\rb_\tau}$ be the skew-symmetric bicharacter on $\Zset^N$ obtained from $\rb_\tau$ as in \eqref{Om}.

Finally, define  group homomorphisms $Y_\tau, \Ytil_\tau : \Zset^N \rightarrow \Fract(R)^*$ such that
\begin{equation}
\label{Ytil}
Y_\tau(e_k) := y_{\tau,k} \; \; \text{and} \; \; \Ytil_\tau(e_k) := y_{\tau, (\tau_\bu \tau)^{-1} (k)} , \; \; \forall k \in [1,N],
\end{equation}
and set
\begin{equation}
\label{ybtiltau}
\ybtil_\tau := (\Ytil_\tau(e_1), \dots, \Ytil_\tau(e_N)) = ( y_{\tau, (\tau_\bu \tau)^{-1} (1)}, \dots, y_{\tau, (\tau_\bu \tau)^{-1} (N)} ) \in \Fract(R)^N.
\end{equation}
It follows from \thref{PCGLalmostcluster}(a) that $\ybtil_\tau$ is a sequence of algebraically independent elements generating $\Fract(R)$ over $\KK$. Since $\id_\bu = \id$, we have $\Ytil_{\id} = Y$ and $\ybtil_{\id} = \ybtil$, as well as $\rb_\id = \qb_\id = \qb$.

For a subset $X \subseteq \ex$, by an $N \times X$ matrix we will mean 
a matrix of size $N \times |X|$ whose columns are indexed by the set
$X$. The set of such matrices with integral entries will be denoted 
by $M_{N \times X} (\Zset)$.

The next theorem contains the main result of the paper. 

\bth{cluster} Let $R$ be a symmetric Poisson-CGL extension of length $N$ and rank $\rk(R)$
as in Definition {\rm\ref{dsymmPCGL}}. Define $\ex \subset [1,N]$ by \eqref{ex-prim}. Assume also that there exist
positive integers $d_i$, $i \in {\mathrm{range}}(\eta)$ such that 
\begin{equation}
\label{d-prop}
d_{\eta(j)} \la_l^\sy = d_{\eta(l)} \la_j^\sy, \quad
\forall j, l \in \ex,
\end{equation}
recall the equality \eqref{la-eq}. Let the sequence of 
generators $x_1, \ldots, x_N$ of $R$ be normalized {\rm(}rescaled\/{\rm)}
so that \eqref{pi-cond} is satisfied {\rm(}recall Proposition {\rm\ref{presc})}. 

Then the following hold:

{\rm(a)} For all $\tau \in \Xi_N$ {\rm(}see \eqref{tau}{\rm)} and $l \in \ex$, 
there exists a unique vector $b_\tau^l \in \Zset^N$ such that $\chi_{\Ytil_\tau(b_\tau^l)}=0$ and 
\begin{equation}
\label{linear-eq}
\Om_{\rb_\tau} ( b^l_\tau, e_j) = 0, \; \; \forall j \in [1,N], \; j \neq l
\quad \mbox{and} \quad
\Om_{\rb_\tau} (b^l_\tau, e_l) = \la_l^*.
\end{equation}
Denote by $\Btil_\tau \in M_{N \times \ex }(\Zset)$ the matrix with columns 
$b^l_\tau$, $l \in \ex$. Let $\Btil:= \Btil_\id$.
 
{\rm(b)} For all $\tau \in \Xi_N$, the pair $(\ybtil_\tau, \Btil_\tau)$ is a 
seed for $\Fract(R)$ and $(\rb_\tau, \Btil_\tau)$ is a compatible pair. The principal part of $\Btil_\tau$ 
is skew-symmetrizable via the integers $d_{\eta(k)}$, $k \in \ex$.

{\rm(c)} All seeds as in {\rm(b)} are mutation-equivalent to each other.
More precisely, they are linked by the following one-step mutations.
Let $\tau, \tau' \in \Xi_N$ be such that
$$
\tau' = ( \tau(k), \tau(k+1)) \tau = \tau (k, k+1)
$$
for some $k \in [1,N-1]$.
If $\eta(\tau(k)) \neq \eta (\tau(k+1))$, then $(\ybtil_{\tau'}, \Btil_{\tau'}) = (\ybtil_\tau, \Btil_\tau)$.
If $\eta(\tau(k)) = \eta (\tau(k+1))$, 
then $(\ybtil_{\tau'}, \Btil_{\tau'}) = \mu_{k_\bu}(\ybtil_\tau, \Btil_\tau)$, where $k_\bu = \tau_\bu \tau(k)$.

{\rm(d)} We have the following equality between the P-CGL extension $R$ and the cluster and upper cluster algebras associated to $(\ybtil, \Btil)$:
$$
R = \AA(\ybtil, \Btil)_\KK = \UU(\ybtil, \Btil, \varnothing)_\KK.
$$
In particular, $\AA(\ybtil, \Btil)_\KK$ is affine and noetherian, and more precisely $\AA(\ybtil, \Btil)_\KK$ is 
generated by the cluster and frozen variables in the seeds parametrized by the subset 
$\Ga_N$ of $\Xi_N$, recall \eqref{sequence}.

{\rm (e)} The cluster algebra structure on $R = \AA(\ybtil, \Btil)_\KK$ is compatible with the given Poisson structure on $R$. Namely, the $N$-tuple of variables in any seed of $\Fract(R)$ mutation-equivalent to $(\xbtil, \Btil)$ is log-canonical. In particular, for each $\tau \in \Xi_N$, the variables $\Ytil_\tau(e_i)$ in the seed $(\ybtil_\tau, \Btil_\tau)$ satisfy
$$
\{ \Ytil_\tau(e_l), \Ytil_\tau(e_j) \} = \Om_{\rb_\tau} (e_l, e_j) \Ytil_\tau(e_l) \Ytil_\tau(e_j), \; \; \forall l,j \in [1,N].
$$

{\rm(f)} Let $\inv$ be any subset of the set $P(N)$ of frozen indices, cf. \eqref{RPr}.
Then 
$$
R[y_k^{-1} \mid k \in \inv] = \AA(\ybtil, \Btil, \inv)_\KK = \UU(\ybtil, \Btil, \inv)_\KK.
$$
Moreover, this cluster algebra structure is compatible with the given Poisson structure on $R[y_k^{-1} \mid k \in \inv]$.
\eth 

\thref{cluster} is proved in \S \ref{recursive}--\ref{complete.main}. The strategy 
of the proof is summarized in \S \ref{overviewmain}.
In \S \ref{clustervar} 
we derive an explicit formula for the cluster variables 
of the seeds that appear in the statement of \thref{cluster}.

\bex{OMmn8}
Let $R = \OMmn$ with the Poisson-CGL extension presentation as in \exref{OMmn}. Define $\eta$ for $R$ as in \exref{OMmn2} and $\ex$ as in \eqref{ex-prim}. Then
$$
\ex = \{ (r-1)n+c \mid r \in [1,m-1], \; c \in [1,n-1] \}.
$$
As shown in \exref{OMmn4}, the P-CGL extension presentation of $R$ is symmetric, and all $\la_j^* = 2$. Hence, \eqref{d-prop} is trivially satisfied with $d_i := 1$ for all $i \in \text{range}(\eta)$. Finally, \exref{OMmn7} shows that condition \eqref{pi-cond} holds.

Therefore all hypotheses of \thref{cluster} are satisfied, and the theorem yields a cluster algebra structure on $R$ compatible with the given Poisson structure. The initial  mutation matrix $\Btil = (b_{ij})$ for $R$ is easily computed:
$$
b_{(r-1)n +c, \, (r'-1)n +c'} = \begin{cases}
\pm1, &\text{if} \; \; r=r', \; c'=c\pm1  \\
 &\text{or} \; \; c=c', \; r'=r\pm1  \\
 &\text{or} \; \; r=r'\pm1, \; c=c'\pm1,  \\
0, &\text{otherwise}
\end{cases} 
\; \; \forall r,r' \in [1,m], \; c,c' \in [1,n] .
$$

From \thref{cluster}(f), we also get Poisson-compatible cluster algebra structures on coordinate rings of various open subvarieties of $\Mmn$ defined by non-vanishing of certain minors. For instance, assume $m \le n$, choose
$$
\inv = \{ rn \mid r \in [1,m] \} \sqcup \{  (m-1)n+c \mid c \in [1,m] \},
$$
and let $R' := R[y_k^{-1} \mid k \in \inv ]$. From \exref{OMmn2}, 
\begin{align*}
y_{rn} &= \De_{[1,r], [n+1-r,n]} \quad (r \in [1,m]),  &y_{(m-1)n+c} = \De_{[m+1-c,m], [1,c]} \quad (c \in [1,m]).
\end{align*}
It follows that $R'$ is the coordinate ring of the generalized double Bruhat cell 
$$
B_m^+ w_1 B_n^+ \cap B^-_m w_2 B^-_n
$$
in $\Mmn$, where $B^\pm_k$ are the standard Borel subgroups of $GL_k(\KK)$ and $w_1$, $w_2$ are the partial permutation matrices
\begin{align*}
w_1 &:= \begin{bmatrix} w^m_\ci &0_{m,n-m} \end{bmatrix} , &w_2 &:= \begin{bmatrix} 0_{m,n-m} &w^m_\ci \end{bmatrix}
\end{align*}
with $w^m_\ci$ the (matrix of the) longest element of $S_m$. To see that
\begin{multline*}
B_m^+ w_1 B_n^+ \cap B^-_m w_2 B^-_n = \{ p\in \Mmn \mid \De_{[1,r], [n+1-r,n]} (p) \ne 0, \; \forall r \in [1,m] \;\;  \\
 \text{and} \;\; \De_{[m+1-c,m], [1,c]} (p) \ne 0, \; \forall c \in [1,m] \},
\end{multline*}
one can use \cite[Proposition 4.1]{BGY}.
\eex

\bre{d-pr} The property \eqref{d-prop} requires that for any particular $j_0 \in \ex$, all the scalars $\la^*_l$ for $l \in \ex$ are positive integer multiples of $\la^*_{j_0}/d_{\eta(j_0)}$. Conversely, if there is some $\la^*_0 \in \kx$ such that $\la^*_l \in \Zset_{> 0}\, \la^*_0$ for all $l \in \ex$, then \eqref{d-prop} is satisfied by choosing $d_{\eta(l)} := \la^*_l/\la^*_0$ for $l \in \ex$, which we can do because $\la_k^* = \la_l^*$ for all $k,l \in \ex$
with $\eta(k) = \eta(l)$ (\prref{la-equal}). This holds, of course, if all $\la^*_l \in \Qset_{> 0}$, which is the case for many P-CGL extensions, such as semiclassical limits of standard uniparameter quantum algebras.

Examples in which non-integer scalars appear among the $\la^*_k$ include semiclassical limits of multiparameter quantum symplectic and euclidean spaces (e.g., see \cite[\S\S2.4, 2.5]{GLaun}). For instance, there are symmetric P-CGL extensions (from \cite[\S2.4]{GLaun}, rewritten in symmetric form) in which $N = 2n$ and $\ex = [1,n-1]$ for an integer $n \ge 2$, and $\la^*_1, \dots, \la^*_{n-1}$ are arbitrary nonzero scalars. When $n \ge 3$ and $\la^*_2/\la^*_1 \notin \Qset_{> 0}$, property \eqref{d-prop} fails.
\ere

\bre{rescaleBtil} For applications, it is useful to determine the exchange matrix $\wt{B}$ \emph{before} the generators $x_k$ have been rescaled to satisfy \eqref{pi-cond}. This is possible because the rescaling does not change $\wt{B}$, as we next note.

Assume $R$ is a symmetric P-CGL extension satisfying all the hypotheses of \thref{cluster} \emph{except} \eqref{pi-cond}. Let the elements $y_k$ and the map $Y = \Ytil_{\id}$ be as above. Now suppose that we rescale the $x_k$ according to \prref{resc}, say with new generators $x'_1,\dots,x'_N$, to make \eqref{pi-cond} hold. Build new $y$-elements from \thref{mainPprime} with respect to these generators, and call them $y'_1, \dots, y'_N$. Each $y'_k$ is a scalar multiple of $y_k$. Since the scalars $\la_{kj}$ and $q_{kj}$ do not change, neither do the bicharacters $\Om_{\lab}$ and $\Om_\qb$. Define $Y'$ analogously to $Y$, using the $y'_k$ in place of the $y_k$. For any vector $b \in \Zset^N$, the element $Y'(b)$ is a scalar multiple of $Y(b)$, and thus $\chi_{Y'(b)} = \chi_{Y(b)}$. Finally, the same elements $h_k, h^*_k \in \h$ which enter into the symmetric P-CGL conditions for the original generators are used with respect to the new generators, which means that the scalars $\la_k$ and $\la^*_k$ do not change under the rescaling. 

Thus, the conditions in \thref{cluster}(a) which uniquely determine the columns of $\wt{B} = \wt{B}_{\id}$ are the same before and after rescaling. Therefore $\wt{B}$ does not change under the rescaling.
\ere

\subsection{Cluster and frozen variables}
\label{clustervar}
The next result gives an explicit formula for the cluster and frozen variables that appear in \thref{cluster}.

\bpr{clustervariables}
Assume the setting of Theorem {\rm\ref{tcluster}}.

{\rm (a)} For $\tau \in \Xi_N$ and $k \in [1,N]$, 
$$
\Ytil_\tau(e_k) = \begin{cases}
y_{[ p^{O_-(k)}( \tau_\bu^{-1}(k) ), \, \tau_\bu^{-1}(k) ]}  &(\text{if} \;\; \tau_\bu^{-1}(k) \ge \tau(1))  \\
y_{[ \tau_\bu^{-1}(k), \ s^{O_-(k)}( \tau_\bu^{-1}(k) ) ]}  &(\text{if} \;\; \tau_\bu^{-1}(k) \le \tau(1)),
\end{cases}
$$
where the predecessor and successor functions are computed
with respect to the original P-CGL extension presentation 
\eqref{itpOre} of $R$.

{\rm (b)} The cluster and frozen variables appearing in the tuples $\ybtil_\tau$ for $\tau$ in $\Xi_N$ or in $\Gamma_N$ are exactly the homogeneous Poisson-prime elements $y_{[i,j]}$. More precisely,
\begin{align*}
\{ \Ytil_\tau(e_k) \mid \tau \in \Xi_N, \; k \in [1,N] \} &= \{ \Ytil_\tau(e_k) \mid \tau \in \Ga_N, \; k \in [1,N] \} \\
&= \{ y_{[i,j]} \mid 1 \le i \le j \le N, \; \eta(i) = \eta(j) \}.
\end{align*}
\epr

\bre{gen} 
\thref{cluster}(d) and \prref{clustervariables} imply that 
the cluster algebra $R = \AA(\ybtil, \wt{B})_\KK$ 
coincides with the subalgebra of $\Fract(R)$ generated by the cluster and frozen 
variables from the (finite) set of seeds $\{ (\ybtil_\tau, \Btil_\tau) \mid \tau \in \Ga_N\}$.
\ere

\begin{proof}[Proof of Proposition {\rm\ref{pclustervariables}}]
(a) Recall that $\Ytil_\tau(e_k) = y_{\tau, (\tau_\bu \tau)^{-1}(k)}$. Set $m := O_-(k)$ and $m_+ := O_+(k)$, and let $L := \eta^{-1}(\eta(k))$ be the level set of $\eta$ containing $k$. Then
$$
L = \{ p^m(k) < p^{m-1}(k) < \cdots < k < s(k) < \cdots < s^{m_+}(k) \}.
$$
Now set $k' := (\tau_\bu \tau)^{-1}(k)$; then $\tau(k') = \tau_\bu^{-1}(k)$. Since $(\tau_\bu \tau)^{-1}$ restricts to an order-preserving bijection of $L$ onto the level set $\tau^{-1}(L)$ of $\eta_\tau$, we see that
$$
\tau^{-1}(L) = \{ p_\tau^m(k') < \cdots < k' < s_\tau(k') < \cdots < s_\tau^{m_+}(k') \}.
$$
Consequently,
\begin{equation}
\label{m+1}
m+1 = \big| \{ j' \in [1,k'] \mid \eta_\tau(j') = \eta_\tau(k') \} \big|
= \big| \{ j \in \tau([1,k']) \mid \eta(j) = \eta(\tau(k')) \} \big|.
\end{equation}

If $\tau(k') \ge \tau(1)$, then $\tau([1,k']) = [\tau(i'), \tau(k')]$ for some $i' \in [1,k']$, and it follows from \eqref{m+1} that $m = \max \{ n \in \Znn \mid p^n(\tau(k')) \in \tau([1,k']) \}$. \thref{sym-Pprime} then shows that
$$
\Ytil_\tau(e_k) = y_{\tau, k'} = y_{[ p^m(\tau(k')), \tau(k')]} = y_{[ p^{O_-(k)}( \tau_\bu^{-1}(k) ), \, \tau_\bu^{-1}(k) ]}
$$
in this case. The case when $\tau(k') \le \tau(1)$ follows similarly.

(b) The first two of the displayed sets are contained in the third by part (a), and the second is contained in the first a priori. It remains to show that the third is contained in the second.

Thus, let $1 \le i \le j \le N$ with $\eta(i) = \eta(j)$. Then $i = p^m(j)$ where
$$
m = \max \{ n \in \Znn \mid p^n(j) \in [i,j] \}.
$$
Set $k' := j-i+1$, and choose $\tau \in \Ga_N$ as follows:
$$
\tau = \begin{cases}
\tau_{1,1} = \id  &(\text{if} \; i = 1)  \\
\tau_{i-1,j} = [i,\dots,j, i-1, j+1, \dots, N, i-1, \dots, 1]  &(\text{if} \; i > 1).
\end{cases}
$$
Then $\tau(1) = i \le j = \tau(k')$ and $\tau([1,k']) = [i,j]$. Consequently,
$$
m = \max \{ n \in \Znn \mid p^n(\tau(k')) \in \tau([1,k']) \},
$$
and \thref{sym-Pprime} shows that
$$
y_{\tau,k'} = y_{[ p^m(\tau(k')), \tau(k')]} = y_{[ i,j ]}.
$$
Therefore $y_{[i,j]} = y_{\tau, (\tau_\bu \tau)^{-1}(k)} = \Ytil_\tau(e_k)$ where $k := \tau_\bu \tau(k')$.
\end{proof}

\bex{OMmn9}
Let $R = \OMmn$ have the Poisson-compatible cluster algebra structure from \thref{cluster} as in \exref{OMmn8}. \prref{clustervariables} and \exref{OMmn5} show that the cluster and frozen variables appearing in \thref{cluster} are exactly the solid minors within $[1,m] \times [1,n]$. For if $l \in \Znn$ and $r,r+l \in [1,m]$, $c,c+l \in [1,n]$, then 
$$
\Delta_{ [r,r+l], [c,c+l]} = y_{[i, s^l(i)]}, \; \; \text{where} \; \; i= (r-1)n+c
$$
(recall \exref{OMmn5}), and the $y_{[i, s^l(i)]}$ are exactly the cluster and frozen variables $\Ytil_\tau(e_k)$ for $\tau \in \Xi_N$ and $k \in [1,N]$.
\eex

\subsection{Auxiliary results}
\label{auxil}
In this subsection we establish two results that will
be needed for the proof of \thref{cluster}. 
The first one uses Theorems \ref{tPCGLmuta} and \ref{tsym-Pprime} to construct 
mutations between pairs of the tuples $\ybtil_\tau$ 
for $\tau \in \Xi_N$. The corresponding mutations of seeds 
(\thref{cluster}(c)) are constructed in \S \ref{recursive}--\ref{allBtil}.

For $g = \sum_j g_j e_j \in \Zset^N$ set
$$
\supp(g) := \{ j \in [1,N] \mid g_j \neq 0 \}.
$$

\bpr{cluster-tau-ind} 
Let $R$ be a symmetric Poisson-CGL extension of length $N$. Assume that the generators 
of $R$ are rescaled so that the condition \eqref{pi-cond} is satisfied.

Let $\tau, \tau' \in \Xi_N$ be such that
$$
\tau' = ( \tau(k), \tau(k+1)) \tau = \tau (k, k+1)
$$
for some $k \in [1,N-1]$ such that $\tau(k) < \tau(k +1)$.

{\rm(a)} If $\eta(\tau(k)) \neq \eta (\tau(k+1))$, then $\rb_{\tau'} = \rb_\tau$ and $\Ytil_{\tau'} = \Ytil_\tau$.

{\rm(b)} Let $\eta(\tau(k)) = \eta (\tau(k+1))$, and set $k_\bu: = \tau_\bu \tau(k)$. Then 
$k_\bu = \tau'_\bu \tau'(k)$ and 
\begin{equation}
\label{Mtautau'}
\Ytil_{\tau'}(e_j)= 
\begin{cases}
\Ytil_\tau(e_j), &\mbox{if} \; \; j \neq k_\bu
\\
\Ytil_\tau(e_{p(k_\bu)} + e_{s(k_\bu)} - e_{k_\bu}) + 
\Ytil_\tau(g - e_{k_\bu}), &\mbox{if} \; \; j =k_\bu
\end{cases}
\end{equation}
for some $g \in \Znn^N$ 
such that $\supp(g) \cap \eta^{-1} \eta(k_\bu) = \varnothing$ and 
$|\supp (g) \cap \eta^{-1}(a)| \leq 1$ for all $a \in \Zset$. Furthermore, 
the vector $e_{p(k_\bu)} + e_{s(k_\bu)} -g \in \Zset^N$ satisfies the identities
\begin{align}
\label{Omtautau'1}
&\Om_{\rb_\tau}( e_{p(k_\bu)} + e_{s(k_\bu)} -g , e_j) = 0, \; \; \forall j \ne k_\bu,
\\
\label{Omtautau'2}
&\Om_{\rb_\tau}( e_{p(k_\bu)} + e_{s(k_\bu)} -g , e_{k_\bu})
= \la^\sy_{k_\bu},
\end{align}
and
\begin{equation}
\label{chiMub}
\chi_{\Ytil_\tau(e_{p(k_\bu)} + e_{s(k_\bu)} -g)}= 0. 
\end{equation}
\epr

Note that the condition $\tau(k) < \tau(k+1)$ is not essential 
since, if $\tau(k) > \tau(k+1)$ and all other conditions are satisfied, 
then one can interchange the roles of $\tau$ and $\tau'$.

\begin{proof} (a) That $\Ytil_{\tau'} = \Ytil_\tau$ follows from \thref{swapkk+1}(a) applied to the P-CGL extension 
presentation \eqref{tau-pres} of $R$ associated to $\tau$, using that $\tau'_\bu \tau' = \tau_\bu \tau (k,k+1)$ (\leref{taubutau}(c)). Now
$$
\lab_{\tau'} = (\tau')^{-1} \lab \tau' = (k,k+1) \tau^{-1} \lab \tau (k,k+1) = (k,k+1) \lab_\tau (k,k+1),
$$
from which it follows that $\qb_{\tau'} = (k,k+1) \qb_\tau (k,k+1)$. Using $\tau'_\bu \tau' = \tau_\bu \tau (k,k+1)$ again, we conclude that
$$
\rb_{\tau'} = (\tau'_\bu \tau') \qb_{\tau'} (\tau'_\bu \tau')^{-1} = (\tau'_\bu \tau') (k,k+1) \qb_\tau (k,k+1) (\tau'_\bu \tau')^{-1} = (\tau_\bu \tau) \qb_{\tau} (\tau_\bu \tau)^{-1} = \rb_\tau.
$$

(b) In this case, $\tau'_\bu \tau' = \tau_\bu \tau$ by \leref{taubutau}(b). Since $\tau([1,j])$ is an interval for all $j \leq k+1$ and $\tau(k) < \tau(k+1)$, 
we have $\tau'([1,k+1]) = \tau([1,k+1]) = [\tau(i), \tau(k+1)]$ and 
$\tau([1,k]) = [\tau(i), \tau(k+1) -1]$ for some $i \in [1,k]$.
On the other hand, the set
$$
\tau'([1,k]) = \tau'([1,k+1]) \setminus \{ \tau'(k+1) \} 
= \tau([1,k+1]) \setminus \{ \tau(k) \}
$$
must be 
also an interval, so $\tau(k)  = \tau(i)$ and $i=k$. Therefore,
\begin{multline}
\label{tau-ident}
\tau'([1,k+1]) = \tau([1,k+1]) = [\tau(k), \tau(k+1)], \\
\tau([1,k]) = [\tau(k), \tau(k+1) -1] \quad
\mbox{and} \quad
\tau'([1,k]) = [\tau(k)+1, \tau(k+1)].
\end{multline}
This implies that 
$$
\tau(k+1) = s^m(\tau(k)) \; \; \mbox{for some} \; \;  m \in \Zset_{>0}
$$ 
and 
\begin{equation}
\label{tau's}
\eta^{-1} (\eta \tau(k)) \cap \tau([1,k+1]) 
= \{\tau(k), s(\tau(k)), \ldots, s^m(\tau(k)) = \tau(k+1) \}.
\end{equation}
From the assumption $\eta(\tau(k)) = \eta (\tau(k+1))$ and the fact that $\tau(j) = \tau'(j)$ for $j \neq k, k+1$
we infer
$$
\tau^{-1} \eta^{-1} (\eta \tau(k)) = 
(\tau')^{-1} \eta^{-1} (\eta \tau(k)).
$$
By the definition of the permutations $\tau_\bu$ and $\tau'_\bu$,
$$
k_\bu = \tau_\bu \tau(k) = \tau'_\bu \tau' (k) = s^{m-1}(\tau(k)).  
$$

From \thref{swap3}(b) we have $y_{\tau',l} = y_{\tau,l}$ for $l \neq k$. 
Eq.~\eqref{Mtautau'} for $j \neq k_\bu$ follows from this.

The identities in \eqref{tau-ident} and \thref{sym-Pprime} give
\begin{align}
\label{2oly}
y_{\tau,k} &= y_{[\tau(k), s^{m-1}(\tau (k))]},  &y_{\tau, k+1} &= y_{[\tau(k), s^m(\tau(k))]},  &y_{\tau', k} &= y_{[s(\tau(k)),s^m(\tau(k))]}.
\end{align}
Recall the definition of the maps $Y_{[i,s^m(i)]}$ from \eqref{Yismi1}. We will construct an isomorphism $\dot w$ from $\Zset^{k+1}$ to $\Zset e_{\tau(k)} +\cdots+ \Zset e_{s^m(\tau(k))}$ such that $Y_{[\tau(k), s^m(\tau(k))]} \dot w$ and $Y_\tau$ agree on suitable $e_j$.

Recall the definition of the set $P_{[\tau(k), s^m(\tau(k))]}$ from \eqref{Pset}, set
\begin{align*}
A &:= \eta(P_{[\tau(k), s^m(\tau(k))]}) = \eta( \tau([1,k+1] )) \setminus \{ \eta(\tau(k)) \}  \\
Q &:= \{ j \in [1,k-1] \mid \eta(\tau(j)) \ne \eta(\tau(l)), \; \; \forall l \in [j+1,k+1] \},
\end{align*}
and note that $\eta \tau$ restricts to a bijection of $Q$ onto $A$. Thus,
\begin{equation}
\label{|QAP|}
|Q| = |A| = | P_{[\tau(k), s^m(\tau(k))]} |.
\end{equation}
The definition of $Q$ also ensures that
\begin{equation}
\label{etatauQ}
\{ t \in [\tau(k), s^m(\tau(k))] \mid \eta(t) = \eta(\tau(j)) \} \subseteq \tau([1,j]), \; \; \forall j \in Q.
\end{equation}

If $j \in Q$ and $\tau(j) \ge \tau(1)$, then $\tau([1,j]) = [\tau(i_j), \tau(j)]$ for some $i_j \in [1,j]$, and we observe that $\tau(j) \in P_{[\tau(k), s^m(\tau(k))]}$. Moreover, the integer $m_j$ corresponding to $j$ in \thref{sym-Pprime} equals $O_-^{\tau(k)+1}(\tau(j))$, and hence we obtain
\begin{equation}
\label{tauj1}
y_{\tau,j} = y_{[p^{m_j}(\tau(j)), \tau(j)]} = Y_{[\tau(k), s^m(\tau(k))]}(e_{\tau(j)}), \; \; \forall j \in Q \; \; \text{with} \; \; \tau(j) \ge \tau(1).
\end{equation}
On the other hand, if $\tau(j) < \tau(1)$, then $\tau([1,j]) = [\tau(j), \tau(i_j)]$ for some $i_j \in [1,j-1]$. Let $m_j$ denote the integer corresponding to $j$ in \thref{sym-Pprime}, and observe that $s^{m_j}(\tau(j)) = \tau(j^-) \in P_{[\tau(k), s^m(\tau(k))]}$ for some $j^- \in [1,j]$. Moreover, $m_j = O_-^{\tau(k)+1}(\tau(j^-))$, and so
\begin{equation}
\label{tauj2}
y_{\tau,j} = y_{[p^{m_j}(\tau(j^-), \tau(j^-)]} = Y_{[\tau(k), s^m(\tau(k))]}(e_{\tau(j^-)}), \; \; \forall j \in Q \; \; \text{with} \; \; \tau(j) < \tau(1).
\end{equation}
Note also that since $\eta \tau(j^-) = \eta \tau(j)$ and $\eta \tau|_Q$ is injective, $\tau(j^-) \ne \tau(i^-)$ for all $i \in Q \setminus \{j\}$ with $\tau(i) < \tau(1)$, while $\tau(j^-) \ne \tau(i)$ for all $i \in Q \setminus \{j\}$ with $\tau(i) \ge \tau(1)$.

In case $m > 1$, we set
\begin{align*}
t &:= \max ( \tau^{-1} \{ s (\tau(k)), \ldots, s^{m-1}(\tau(k)) \})  \\
 &\; = \max \{ j \in [1,k-1] \mid \eta(\tau(j)) = \eta(\tau(k)) \} = p_\tau(k)
\end{align*}
(see \eqref{tau's}). Then either $\tau(t) = s^{m-1}(\tau(k))$ or $\tau(t) = s(\tau(k))$, and \thref{sym-Pprime} yields
\begin{equation}
\label{ytaul}
y_{\tau,t} = y_{[s(\tau(k)), s^{m-1}(\tau(k))]} = Y_{[\tau(k), s^m(\tau(k))]}( e_{s^{m-1}(\tau(k))} ), \; \; \text{if} \; \; m > 1.
\end{equation}

Now choose a bijection $w : [1,k+1] \rightarrow [\tau(k), s^m(\tau(k))]$ such that
\begin{align*}
w(k+1) &= s^m(\tau(k))  \\
w(k) &= \tau(k)  \\
w(t) &= s^{m-1}(\tau(k)), \; \; \text{if} \; \; m > 1  \\
w(j) &= \tau(j), \; \; \forall j \in Q \; \; \text{with} \; \; \tau(j) \ge \tau(1)  \\
w(j) &= \tau(j^-), \; \; \forall j \in Q \; \; \text{with} \; \; \tau(j) < \tau(1),
\end{align*}
and let $\dot w$ denote the isomorphism fron $\Zset^{k+1}$ onto $\sum_{i=\tau(k)}^{s^m(\tau(k))} \Zset e_i$ such that $\dot w(e_j) = e_{w(j)}$ for $j \in [1,k+1]$.
In particular, $\eta w(j) = \eta \tau(j)$ for $j \in Q$, so $\eta w|_Q$ is injective. By construction, $w(Q) \subseteq P_{[\tau(k), s^m(\tau(k))]}$, and so we conclude from \eqref{|QAP|} that
\begin{equation}
\label{wQP}
w|_Q : Q \longrightarrow P_{[\tau(k), s^m(\tau(k))]} \quad\text{is a bijection.}
\end{equation}
Combining Eqs. \eqref{2oly} and \eqref{tauj1}--\eqref{ytaul} with the definition of $w$, we see that
\begin{equation}
\label{MhatMw}
Y_\tau(e_j) = y_{\tau,j} = Y_{[\tau(k), s^m(\tau(k))]}\dot w(e_j), \; \; \begin{cases}
\forall j \in Q \cup \{t,k,k+1\} &(m > 1) \\
\forall j \in Q \cup \{k,k+1\} &(m = 1).
\end{cases}
\end{equation}

The next step is to apply \thref{PCGLmuta}. We do the case $m > 1$ and leave the case $m = 1$ to the reader. (In the latter case, $p(k_\bu) = -\infty$ and $e_{p(k_\bu)} = 0$.) Observe that \eqref{MhatMw} implies that
\begin{equation}
\label{Mhatf}
Y_\tau(f) = Y_{[\tau(k), s^m(\tau(k))]}\dot w(f), \; \; \forall f \in \Zset^N \; \; \text{with} \; \; \supp(f) \subseteq Q \cup \{t,k,k+1\}.
\end{equation}
Thus, taking account of \eqref{2oly}, \eqref{wQP}, and \eqref{MhatMw}, \thref{PCGLmuta} implies that
$$
y_{\tau',k} = 
Y_\tau(e_t + e_{k+1} - e_k) + Y_\tau(g'-e_k),
$$
where $g' := \dot w^{-1}( g_{[ \tau(k), s^m(\tau(k)) ]} )$, recall \eqref{gismi}. Since 
$$
\supp( g_{[ \tau(k), s^m(\tau(k)) ]} ) \subseteq P_{[ \tau(k), s^m(\tau(k)) ]} =w(Q),
$$
we have $\supp(g') \subseteq Q$.  By the definition of $\tau_\bu$ and Eq. \eqref{tau's},
$$
\tau_\bu \tau(t) = \tau_\bu \tau( p_\tau(k) ) = p(k_\bu) \; \; 
\mbox{and} \; \; \tau_\bu \tau(k+1) = s(k_\bu).
$$
Therefore,
$$
\Ytil_{\tau'}(e_{k_\bu}) = Y_{\tau'}(e_k)
= y_{\tau', k}
= \Ytil_\tau(e_{p(k_\bu)} + e_{s(k_\bu)} - e_{k_\bu}) + \Ytil_\tau(\tau_\bu \tau(g')-e_{k_\bu}),
$$
which implies the validity of \eqref{Mtautau'} for $j = k_\bu$, where $g := \tau_\bu \tau(g')$.

Finally, the identities \eqref{Omtautau'1}--\eqref{Omtautau'2} follow from
\thref{Omqval} (applied to the P-CGL presentation \eqref{tau-pres}), the definition of $\rb_\tau$, and the fact that 
$R_{\tau,k+1} = R_{[\tau(k),s^m(\tau(k))]}$, 
see \eqref{tau-ident}. We note that $\eta(k_\bu) = \eta( \tau(k) )$ and $s(\tau(k)) \neq + \infty$, 
$s(k_\bu) \neq + \infty$,
which follow from the definition of $\tau_\bu$ and Eq. \eqref{tau's}. Because of 
this and \prref{la-equal}, $\la^*_{\tau'(k+1)} = \la_{\tau(k)}^\sy = \la_{k_\bu}^\sy$. The identity \eqref{chiMub} follows 
from the fact that $y_{\tau, j}$ and $y_{\tau',j}$ are 
$\HH$-eigenvectors for all $j \in [1,N]$ and Eq. \eqref{Mtautau'}.
\end{proof}

Our second auxiliary result relies on a strong rationality property analogous to \cite[Theorem II.6.4]{BrGo}. Given a ring $A$ equipped with an action of a group $\HH$ by automorphisms, we write $A^\HH$ for the fixed subring $\{ a \in A \mid h\cdot a = a \; \text{for all} \; h \in \HH \}$.

\bth{strongrat}
Let $R$ be a Poisson-CGL extension of length $N$, and $P$ an $\HH$-Poisson-prime ideal of $R$. Then
\begin{equation}
\label{HPrat}
Z_p( \Fract R/P )^\HH = \KK.
\end{equation}
\eth

\begin{proof} We induct on $N$, the case $N = 0$ (when $R = \KK$) being trivial. Now let $N > 0$, and assume the theorem holds for $R_{N-1}$.

The contraction $Q := P \cap R_{N-1}$ is an $\HH$-Poisson-prime ideal of $R_{N-1}$ and is stable under $\sig_N$ and $\de_N$. Hence, $R/QR$ is a Poisson polynomial ring of the form
$$
R/QR = (R_{N-1}/Q)[ x_N; \sig_N, \de_N]_p ,
$$
where we have omitted overbars on cosets and induced maps. Let $S$ be the localization of $R_{N-1}/Q$ with respect to the set of all $\HH$-eigenvectors in $R_{N-1}/Q$, and note that $S$ has no nonzero proper $\HH$-stable ideals. The induced action of $\HH$ on $R_{N-1}/Q$ extends to a Poisson action on $S$, and since $\Fract R_{N-1}/Q = \Fract S$, our induction hypothesis implies that
$$
Z_p(\Fract S)^\HH = \KK.
$$

Our choice of $S$ ensures that the $\HH$-action on $S$ is rational. There are unique extensions of $\sig_N$ and $\de_N$ to a Poisson derivation $\sig_N$ on $S$ and a Poisson $\sig_N$-derivation $\de_N$ on $S$. The differential of the $\HH$-action gives an action of $\Lie \HH$ on $S$, the action of $(h_N\cdot)$ on $S$ extends that on $R_{N-1}/Q$, and $(h_N \cdot) = \sig_N$ on $S$. Finally, $R/QR$ extends to a Poisson polynomial ring
$$
\Rhat := S[ x_N; \sig_N, \de_N]_p ,
$$
and the $\HH$-action on $R/QR$ extends to a rational Poisson action on $\Rhat$. The ideal $P/QR$ induces an $\HH$-Poisson-prime ideal $\Phat$ of $\Rhat$ and there is an $\HH$-equivariant Poisson isomorphism $\Fract \Rhat/\Phat \cong \Fract R/P$, so it suffices to prove that
\begin{equation}
\label{Pratgoal}
Z_p( \Fract \Rhat/\Phat )^\HH = \KK.
\end{equation}

There are two situations to consider, depending on whether $\Rhat$ has any proper nonzero $\HH$-stable Poisson ideals. Assume first that it does not; in particular, $\Phat = 0$.

Consider a nonzero element $u \in Z_p(\Fract \Rhat)^\HH$, set $I := \{ r \in \Rhat \mid ru \in \Rhat \}$, and observe that $I$ is a nonzero $\HH$-stable Poisson ideal of $\Rhat$. Our current assumption implies $I = \Rhat$, whence $u \in \Rhat$. Similarly, $u^{-1} \in \Rhat$, forcing $u \in S$. Thus, $u \in Z_p(\Fract S)^\HH = \KK$. 

Now assume that $\Rhat$ does have proper nonzero $\HH$-stable Poisson ideals. As in the proof of \cite[Proposition 1.2]{GLaun}, $\de_N$ is an inner Poisson $\sig_N$-derivation of $S$, implemented by a homogeneous element $d \in S$ with $\chi_d = \chi_{x_N}$, and the only $\HH$-Poisson-prime ideals of $\Rhat$ are $0$ and $z\Rhat$, where $z := x_N - d$. Moreover, $\Rhat = S[z; \sig_N]_p$. There is an $\HH$-equivariant Poisson isomorphism $\Rhat/z\Rhat \cong S$, and so $Z_p(\Fract \Rhat/z\Rhat)^\HH = \KK$. This establishes \eqref{Pratgoal} when $\Phat = z\Rhat$.

The case $\Phat = 0$ remains. Consider a nonzero element  $u \in Z_p(\Fract \Rhat)^\HH$. As above, $I := \{ r \in \Rhat \mid ru \in \Rhat \}$ is a nonzero $\HH$-stable Poisson ideal of $\Rhat$. If $I \ne \Rhat$, any prime ideal minimal over $I$ is an $\HH$-Poisson-prime ideal (\leref{torusPprime}) and so must equal $z\Rhat$, in which case $\sqrt{I} = z\Rhat$. Whether or not $I = \Rhat$, it follows that $z^n \in I$ for some $n \ge 0$. Hence, $u = v z^{-n}$ for some homogeneous element $v \in \Rhat$.  Note that $\{z,u\} = 0$ implies $\{z,v\} = 0$. Write $v = v_0+ v_1z+ \cdots+ v_tz^t$ for some $v_i \in S$. Now
\begin{align*}
0 &= \{z,v\} = \{z,v_0\} + \{z,v_1\} z + \cdots + \{z,v_t\} z^t  \\
&= \sig_N(v_0)z + \sig_N(v_1) z^2 + \cdots + \sig_N(v_t) z^{t+1},
\end{align*}
whence $h_N \cdot v_i = \sig_N(v_i) = 0$ for $i \in [0,t]$. Then, since $v= uz^n$ and $h_N \cdot u = 0$ by \eqref{h.0}, we get
$$
n \la_N v = h_N \cdot v = \la_N v_1 z + \cdots + t \la_N v_t z^t.
$$
Since $\la_N \ne 0$, it follows that $v_i = 0$ for all $i \ne n$, and so $u = v_n \in S$. Therefore $u \in Z_p(\Fract S)^\HH = \KK$, yielding \eqref{Pratgoal} in the case $\Phat = 0$.
\end{proof}

The next lemma proves uniqueness of integral vectors satisfying bilinear identities 
of the form \eqref{Omtautau'1}--\eqref{Omtautau'2}.

For $\tau \in \Xi_N$, applying \eqref{brackYfYg} to the P-CGL extension \eqref{tau-pres} yields
$$
\{ Y_\tau(f), Y_\tau(g) \} = \Om_{\qb_\tau}(f,g) Y_\tau(f) Y_\tau(g), \; \; \forall f,g \in \Zset^N,
$$
and consequently
\begin{equation}
\label{brackYtiltau}
\{ \Ytil_\tau(f), \Ytil_\tau(g) \} = \Om_{\rb_\tau}(f,g) \Ytil_\tau(f) \Ytil_\tau(g), \; \; \forall f,g \in \Zset^N.
\end{equation}

\ble{unique} Assume that $R$ is a symmetric P-CGL extension of length $N$. For any $\tau \in \Xi_N$, 
$\theta \in \xh$, and $\xi_1, \ldots, \xi_N \in \KK$, there exists 
at most one vector $b \in \Zset^N$ such that $\chi_{\Ytil_\tau(b)}= \theta$ and 
$\Om_{\rb_\tau}(b, e_j) = \xi_j$ for all $j \in [1,N]$.
\ele

\begin{proof}
Let $b_1, b_2 \in \Zset^N$ be such that $\chi_{\Ytil_\tau(b_1)} = \chi_{\Ytil_\tau(b_2)} = \theta$ and 
$\Om_{\rb_\tau}(b_1, e_j) = \Om_{\rb_\tau}(b_2, e_j)$  $= \xi_j$ for all $j \in [1,N]$.
Taking account of \eqref{brackYtiltau}, we find that
$$
\{ \Ytil_\tau(b_1 - b_2), \Ytil_\tau(e_j) \} = 0, \; \; \forall j \in [1,N].
$$
This implies that $\Ytil_\tau(b_1 - b_2)$ belongs to the Poisson center of $\Fract R$, because
by \thref{PCGLalmostcluster}(a), $\Ytil_\tau(e_1), \dots, \Ytil_\tau(e_N)$ generate the field $\Fract R$ over $\KK$.
Furthermore, 
$$
\chi_{\Ytil_\tau(b_1 - b_2)} =0,
$$
so that $\Ytil_\tau(b_1 - b_2)$ is fixed by $\HH$. By \thref{strongrat}, $\Ytil_\tau(b_1 - b_2) \in \KK$. This is only possible if $b_1=b_2$, because $\Ytil_\tau(e_1), \dots, \Ytil_\tau(e_N)$ are algebraically independent over $\KK$.
\end{proof}

\subsection{An overview of the proof of \thref{cluster}}
\label{overviewmain}
In this subsection we give a summary of the strategy of 
our proof of \thref{cluster}.

In \S \ref{statement} we constructed $N$-tuples $\ybtil_\tau \in \Fract(R)^N$ 
associated to the elements 
of the set $\Xi_N$. In order to extend them to seeds of $\Fract(R)$, 
one needs to construct a compatible matrix $\wt{B}_\tau \in M_{N \times \ex}(\Zset)$
for each of them.
This will be first done for the subset $\Ga_N$ of $\Xi_N$ 
in an iterative fashion with respect to the linear ordering \eqref{sequence}.
If $\tau$ and $\tau'$ are two consecutive 
elements of $\Ga_N$ in that linear ordering, then $\tau' = \tau(k,k+1)$ 
for some $k \in [1,N-1]$ such that $\tau(k) < \tau(k+1)$. 
If $\eta(\tau(k)) \neq \eta(\tau(k+1))$ then $\ybtil_{\tau'} = \ybtil_\tau$ 
by \prref{cluster-tau-ind}(a)
and nothing happens at that step. If $\eta(\tau(k)) = \eta(\tau(k+1))$,  
then we use \prref{cluster-tau-ind}(b) to construct $b_\tau^{k_\bu}$ and $b_{\tau'}^{k_\bu}$ 
where $k_\bu := (\tau_\bu \tau)(k)$. Up to $\pm$ sign these vectors 
are equal to $e_{p(k_\bu)} + e_{s(k_\bu)} -g$, where $g \in \Znn^N$ 
is the vector from \prref{cluster-tau-ind}(b). Then we use ``reverse'' mutation to construct $b_\sig^{k_\bu}$ for 
$\sig \in \Ga_N$, $\sig \prec \tau$ in the linear order \eqref{sequence}. 
Effectively this amounts to starting with a cluster algebra in which all 
variables are frozen and then recursively adding more exchangeable variables.

There are two things that can go wrong with this.
Firstly, the reverse mutations from 
different stages might not be synchronized. Secondly, there are many pairs 
of consecutive elements $\tau, \tau'$ for which $k_\bu$ is the same. 
So we need to prove that $b_\sig^{k_\bu}$ is not {\em{overdetermined}}.
We use strong rationality of P-CGL extensions to handle both via \leref{unique}.
This part of the proof (of parts (a) and (b) of \thref{cluster}) 
is carried out in \S \ref{recursive}.

Once $\wt{B}_\id$ is (fully) constructed then the $\wt{B}_\tau$ are constructed 
inductively by applying \prref{cluster-tau-ind} and using the sequences of elements 
of $\Xi_N$ from \coref{steps}(a). At each step \leref{unique} is applied 
to match columns of mutation matrices. This proves parts (a) and (b) 
of \thref{cluster}. Part (c) of the theorem is obtained in a somewhat 
similar manner from
\prref{cluster-tau-ind}. This is done in \S \ref{allBtil}. Once parts (a)--(c) have been proved, we will obtain the algebra equalities in parts (d) and (f) of \thref{cluster} from \thref{PCGLalmostcluster} and the Poisson-compatibility statements in parts (e) and (f) from \thref{P-cluster-main} and \prref{brack-mut}. This is done in \S\ref{complete.main}.

\subsection{Recursive construction of seeds for $\tau \in \Ga_N$}
\label{recursive}
Recall the linear ordering \eqref{sequence} on $\Ga_N \subset \Xi_N$. 
We start by constructing a chain of subsets $\ex_\tau \subseteq \ex$
indexed by the elements of $\Ga_N$ such that 
$$
\ex_\id = \varnothing, \quad \ex_{w_\ci} = \ex \quad
\mbox{and} \quad
\ex_\sig \subseteq \ex_\tau, \; \; \forall \sig, \tau \in \Ga_N, \;
\sig \prec \tau.
$$
This is constructed inductively by starting with $\ex_\id = \varnothing$.
If $\tau \prec \tau'$ are two consecutive elements in the linear ordering, 
then for some $1 \leq i < j \leq N$,
$$
\tau = \tau_{i,j-1}  \; \; 
\mbox{and} \; \; 
\tau'=\tau_{i,j}
$$
(recall \eqref{tauij} and the equalities in \eqref{sequence}).
Assuming that $\ex_\tau$ has been constructed, we define
$$
\ex_{\tau'} := \begin{cases}
\ex_\tau \cup \{ p(j) \},  &\mbox{if} \; p(i) = -\infty, \; \eta(i) = \eta(j)  \\
\ex_\tau,  &\mbox{otherwise}.
\end{cases}
$$
It is clear that this process ends with $\ex_{w_\ci} = \ex$.

The following lemma provides an inductive 
procedure for establishing parts (a) and (b) of \thref{cluster} for $\tau \in \Ga_N$. In the proof, we allow seeds whose exchange matrices are empty.

\ble{ind-ab} Assume that $R$ is a symmetric P-CGL extension of length $N$ 
satisfying \eqref{d-prop}. Assume also that the generators 
of $R$ are rescaled so that the condition \eqref{pi-cond} is satisfied.

Let $\tau \in \Ga_N$. For all $\sig \in \Ga_N$ with $\sig \preceq \tau$,
there exists a unique matrix $\wt{B}_{\sig,\tau}$ in $M_{N \times \ex_\tau}(\Zset)$ whose  columns $b_{\sig, \tau}^l \in \Zset^N$ satisfy 
\begin{equation}
\label{sigtau}
\chi_{\Ytil_\sig(b_{\sig,\tau}^l)}= 0, \quad
\Om_{\rb_\sig} ( b^l_{\sig,\tau}, e_j) = 0, \; \; \forall j \in [1,N], \; j \neq l
\quad \mbox{and} \quad
\Om_{\rb_\sig} (b^l_{\sig,\tau}, e_l) = \la_l^*
\end{equation}
for all $l \in \ex_\tau$.
The matrix $\wt{B}_{\sig, \tau}$ has full rank, and its principal part is skew-symmetrizable 
via the integers $d_{\eta(k)}$, $k \in \ex_\tau$. Moreover, $(\rb_\sig, \Btil_{\sig,\tau})$ is a compatible pair, for all $\sig \preceq \tau$.
\ele  

\begin{proof} The uniqueness statement follows at once from \leref{unique}. 
If a matrix $\wt{B}_{\sig, \tau}$ with the properties \eqref{sigtau} exists, then, in view of \eqref{compscalars}, $(\rb_\sig, \wt{B}_{\sig,\tau})$ is a compatible pair.
The scalars $\be_l^{\sig,\tau} := (\Btil^T_{\sig,\tau} \rb_\sig)_{ll}$ satisfy
$$
\be_l^{\sig,\tau} = \sum_{i=1}^N (\Btil_{\sig,\tau})_{il} (\rb_\sig)_{il} = \sum_{i=1}^N (\Btil_{\sig,\tau})_{il} \Om_{\rb_\sig}(e_i,e_l) = \Om_{\rb_\sig}( b^l_{\sig,\tau}, e_l) = \la^*_l, \; \; \forall l \in \ex_\tau .
$$
Consequently, the principal part of $\Btil_{\sig,\tau}$ is skew-symmetrizable 
via the integers $d_{\eta(k)}$ by \leref{alternative} and the condition \eqref{d-prop}. Moreover, \prref{fullrank} implies that $\Btil_{\sig,\tau}$ has full rank. Thus, $(\ybtil_\sig, \Btil_{\sig,\tau})$ is a seed in $\FF$.

What remains to be proved is the existence statement in the 
lemma. It trivially holds for $\tau =\id$ since 
$\ex_\id = \varnothing$.

Let $\tau \prec \tau'$ be two consecutive elements of $\Ga_N$ 
in the linear ordering \eqref{sequence}. Assuming that the 
existence statement in the lemma holds for $\tau$, we will show 
that it holds for $\tau'$. The lemma will then follow by induction.

As noted above, for some $1 \leq i < j \leq N$ we have
$\tau = \tau_{i,j-1}$ and $\tau'=\tau_{i,j}$. In particular, 
$\tau' = (ij) \tau = \tau(j-i, j-i+1)$ and $\tau(j-i) = i < j = \tau(j-i+1)$,
so \prref{cluster-tau-ind} is applicable to the pair $(\tau, \tau')$, with $k := j-i$. Note that $\tau(k) = i$ and $\tau(k+1) = j$.

If $\eta(i) \neq \eta(j)$, then $\ex_{\tau'} = \ex_\tau$, 
and $\Ytil_{\tau'} = \Ytil_\tau$ and $\rb_{\tau'} = \rb_\tau$ by \prref{cluster-tau-ind}(a).
So, $\Om_{\rb_{\tau'}} = \Om_{\rb_\tau}$.
These identities imply that the following matrices have the properties 
\eqref{sigtau} for the element $\tau' \in \Ga_N$: $\wt{B}_{\sig,\tau'} := \wt{B}_{\sig, \tau}$ 
for $\sig \preceq \tau$ and $\wt{B}_{\tau',\tau'} := \wt{B}_{\tau, \tau}$.

Next, we consider the case $\eta(i) = \eta(j)$. This  implies that 
$j = s^m(i)$ for some $m \in \Zset_{>0}$. This fact and the definition of $\tau_\bu$
give that the element $\tau_\bu \tau(j-i)$ equals the $m$-th element 
of $\eta^{-1}(\eta(i))$ when the elements in the preimage are ordered from least to greatest.
Therefore this element is explicitly given by
\begin{equation}
\label{j-i}
\tau_\bu \tau(j-i) = s^{m-1} p^{O_-(i)}(i).
\end{equation}
Now set
$$
k_\bu:= \tau_\bu \tau(j-i)
$$
as in \prref{cluster-tau-ind}(b).
There are two subcases: (1) $p(i) \ne -\infty$ and (2) $p(i) = -\infty$.

{\bf Subcase (1)}. In this situation $\ex_{\tau'} = \ex_{\tau}$, so we do not need to generate an ``extra column'' for each matrix.
Set $\wt{B}_{\sig, \tau'}:= \wt{B}_{\sig, \tau}$ for $\sig \in \Ga_N$, 
$\sig \preceq \tau$. Eq. \eqref{sigtau} for the pairs $(\sig, \tau')$ 
with $\sig \preceq \tau$ follows from the equality $\ex_{\tau'} = \ex_\tau$.

Next we deal with the pair $(\sig,\tau') = (\tau', \tau')$.
Applying the inductive assumption \eqref{sigtau} for $\wt{B}_{\tau, \tau}$ and \prref{cluster-tau-ind}(b) shows that the vector $g \in \Znn^N$
has the properties 
$$
\chi_{\Ytil_\tau(b_{\tau,\tau}^{k_\bu})} = 0
= \chi_{\Ytil_\tau(e_{p(k_\bu)}+e_{s(k_\bu)} -g)}  
$$
and 
$$
\Om_{\rb_\tau}(b_{\tau,\tau}^{k_\bu}, e_t) = 
\Om_{\rb_\tau}( e_{p(k_\bu)}+e_{s(k_\bu)} -g, e_t), \; \; \forall t \in [1,N].
$$
\leref{unique} implies that $e_{p(k_\bu)}+e_{s(k_\bu)} -g = b_{\tau,\tau}^{k_\bu}$.
It follows from this and Eq. \eqref{Mtautau'} that 
\begin{equation}
\label{Ytiltau'}
\begin{aligned}
\Ytil_{\tau'}(e_j) &= \Ytil_\tau(e_j), \; \; \forall j \ne k_\bu ,  \\
\Ytil_{\tau'}(e_{k_\bu}) &= \Ytil_\tau( - e_{k_\bu} + [b_{\tau,\tau}^{k_\bu}]_+) + \Ytil_\tau( - e_{k_\bu} - [b_{\tau,\tau}^{k_\bu}]_-),
\end{aligned}
\end{equation}
whence $\ybtil_{\tau'} = \mu_{k_\bu}(\ybtil_\tau)$.

We set $\Btil_{\tau',\tau'} := \mu_{k_\bu}(\Btil_{\tau,\tau})$, so that $(\ybtil_{\tau'}, \Btil_{\tau',\tau'}) = \mu_{k_\bu}(\ybtil_\tau, \Btil_{\tau,\tau})$. Since the columns of $\Btil_{\tau,\tau}$ satisfy \eqref{sigtau} and the entries $\Ytil_\tau(e_j)$ of $\ybtil_\tau$ are $\HH$-eigenvectors,  \leref{equivar-mut}
implies that the columns of $\wt{B}_{\tau', \tau'}$ have the property $\chi_{\Ytil_{\tau'}(b^l_{\tau', \tau'})}= 0$
for all $l \in \ex_\tau = \ex_{\tau'}$. 

According to \eqref{brackYtiltau},
\begin{align}
\{ \Ytil_\tau(e_l), \Ytil_\tau(e_j) \} &= \Om_{\rb_\tau}(e_l, e_j) \Ytil_\tau(e_l) \Ytil_\tau(e_j), \; \; \forall l,j \in [1,N], \label{Ytt}  \\
\{ \Ytil_{\tau'}(e_l), \Ytil_{\tau'}(e_j) \} &= \Om_{\rb_{\tau'}}(e_l, e_j) \Ytil_{\tau'}(e_l) \Ytil_{\tau'}(e_j), \; \; \forall l,j \in [1,N].  \label{Yt't'}
\end{align}
On the other hand, since the $\Ytil_\tau(e_i)$ and $\Ytil_{\tau'}(e_i)$ are the entries of $\ybtil_\tau$ and $\ybtil_{\tau'} = \mu_{k_\bu}(\ybtil_\tau)$, it follows from \eqref{Ytt} and \prref{brack-mut} that
\begin{equation}
\label{Yt't'2}
\{ \Ytil_{\tau'}(e_l), \Ytil_{\tau'}(e_j) \} = \Om_{\mu_{k_\bu}(\rb_\tau)}(e_l, e_j) \Ytil_{\tau'}(e_l) \Ytil_{\tau'}(e_j), \; \; \forall l,j \in [1,N]. 
\end{equation}
Comparing \eqref{Yt't'} and \eqref{Yt't'2}, we see that $\Om_{\rb_{\tau'}} = \Om_{\mu_{k_\bu}(\rb_\tau)}$, and therefore
\begin{equation}
\label{rtau'-mukbu}
\rb_{\tau'} = \mu_{k_\bu}(\rb_\tau).
\end{equation}

Now \prref{pair-mut}(b) implies that $\Btil^T_{\tau',\tau'} \rb_{\tau'} = \Btil^T_{\tau,\tau} \rb_\tau$. Consequently,
$$
\Om_{\rb_{\tau'}}( b^l_{\tau',\tau'}, e_j) = e^T_l \Btil^T_{\tau',\tau'} \rb_{\tau'} e_j = e^T_l \Btil^T_{\tau,\tau} \rb_\tau e_j = \Om_{\rb_\tau}( b^l_{\tau,\tau}, e_j) = \de_{jl} \la^*_l
$$
for all $l \in \ex_\tau$ and $j \in [1,N]$. Thus, the columns of $\Btil_{\tau',\tau'}$ satisfy the conditions \eqref{sigtau}, which completes the proof of the inductive step of the lemma in this subcase.

{\bf Subcase (2)}. In this case, $k_\bu = s^{m-1}(i) = p(j)$ and $\ex_{\tau'} = \ex_\tau \sqcup \{k_\bu\}$. Define 
the matrix $\wt{B}_{\tau, \tau'}$ by 
$$
b^l_{\tau, \tau'} = 
\begin{cases}
b^l_{\tau, \tau} \,, & \mbox{if} \; \; l \neq k_\bu
\\
e_{p(k_\bu)} + e_{s(k_\bu)} - g, & \mbox{if} \; \; 
l = k_\bu,
\end{cases}
$$
where $g \in \Znn^N$ is the vector from \prref{cluster-tau-ind}(b).
Applying the assumption \eqref{sigtau} for $\wt{B}_{\tau, \tau}$ and \prref{cluster-tau-ind}(b), we obtain 
that the matrix $\wt{B}_{\tau, \tau'}$ has the properties \eqref{sigtau}. 
We set $\wt{B}_{\tau', \tau'}:= \mu_{k_\bu} (\wt{B}_{\tau, \tau'})$.

As in subcase (1), using \leref{equivar-mut} and Propositions \ref{ppair-mut}(b) and \ref{pbrack-mut},  
one derives that $\wt{B}_{\tau', \tau'}$ satisfies the properties \eqref{sigtau}.

We are left with constructing $\wt{B}_{\sig, \tau'} \in M_{N \times \ex_{\tau'}}(\Zset)$ 
for $\sig \in \Ga_N$, $\sig \prec \tau$.
We do this by a downward induction on the linear ordering \eqref{sequence} in a fashion that 
is similar to the proof of the lemma in the subcase (1). Assume that $\sig \prec \sig'$ is a pair of consecutive 
elements of $\Ga_N$ such that $\sig' \preceq \tau$. As in the beginning of the subsection, 
we have that for some $1 \leq i\spcheck < j\spcheck \leq N$,
$$
\sig = \tau_{i\spcheck,j\spcheck-1}  \; \; 
\mbox{and} \; \; 
\sig'=\tau_{i\spcheck,j\spcheck},
$$
so
$$
\tau' = (i\spcheck j\spcheck) \tau = \tau(j\spcheck-i\spcheck, j\spcheck-i\spcheck+1).
$$
Assume that there exists a matrix $\wt{B}_{\sig', \tau'} \in M_{N \times \ex_{\tau'}} (\Zset)$ 
that satisfies \eqref{sigtau}. We define the matrix $\wt{B}_{\sig, \tau'} \in M_{N \times \ex_{\tau'}} (\Zset)$
by 
$$
\wt{B}_{\sig, \tau'}:= 
\begin{cases}
\wt{B}_{\sig', \tau'}, & \mbox{if} \; \; \eta(i\spcheck) \neq \eta(j\spcheck)
\\
\mu_{k\spcheck_\bu}(\wt{B}_{\sig, \tau'}), & \mbox{if} \; \; 
\eta(i\spcheck) = \eta(j\spcheck),
\end{cases}
$$
where
$$
k\spcheck_\bu := \sig_\bu \sig(j\spcheck-i\spcheck).
$$
Analogously to the proof of the lemma in the subcase (1), 
using Propositions \ref{pcluster-tau-ind} and \ref{ppair-mut}(b)
and Lemmas \ref{lunique} and \ref{lequivar-mut}, 
one proves that the matrix $\wt{B}_{\sig, \tau'}$ has the properties \eqref{sigtau}.
This completes the proof of the lemma.
\end{proof}

\bre{restr} It follows from \leref{unique} that the matrices $\wt{B}_{\sig, \tau}$ 
in \leref{ind-ab} have the following restriction property:

{\em{For all triples $\sig \preceq \tau \prec \tau'$ of elements of $\Ga_N$,
$$
b^l_{\sig, \tau} = b^l_{\sig, \tau'}, \; \; \forall l \in \ex_\tau.
$$
In other words, $\wt{B}_{\sig, \tau}$ is obtained from $\wt{B}_{\sig, \tau'}$ by removing 
all columns indexed by the set $\ex_{\tau'} \backslash \ex_\tau$.}}

This justifies that \leref{ind-ab} gradually enlarges 
a matrix $\wt{B}_{\sig, \sig} \in M_{N \times \ex_\sig}(\Zset)$ 
to a matrix $\wt{B}_{\sig, w_0} \in M_{N \times \ex}(\Zset)$, 
for all $\sig \in \Ga_N$. In the case of $\sig = \id$, we start 
with an empty matrix ($\ex_\id = \varnothing$) and obtain a matrix $\wt{B}_{\id, w_\ci} \in M_{N \times \ex}(\Zset)$
which will be the needed mutation matrix for the initial cluster variables $\ybtil$.
\ere

\begin{proof}[Proof of Theorem {\rm\ref{tcluster}(a)(b)} for $\tau \in \Ga_N$] Change $\tau$ to $\sig$ in these statements.
These parts of the theorem for the elements of $\Ga_N$ follow from
\leref{ind-ab} applied to $(\sig,\tau) = (\sig,w_\ci)$. For all $\sig \in \Ga_N$ we set 
$\wt{B}_\sig := \wt{B}_{\sig, w_\ci}$ and use that $\ex_{w_\ci} = \ex$. 
\end{proof}

\subsection{Proofs of parts (a), (b) and (c) of \thref{cluster}}
\label{allBtil}
Next we establish \thref{cluster}(a)(b) in full generality. 
This will be done by using the result of \S \ref{recursive} for $\tau=\id$ 
and iteratively applying the following proposition.

\bpr{tau-seed-muta} Let $R$ be a symmetric Poisson-CGL extension of length $N$ 
satisfying \eqref{d-prop}. Assume also that the generators 
of $R$ are rescaled so that the condition \eqref{pi-cond} is satisfied.
Let $\tau, \tau' \in \Xi_N$ be such that
$$
\tau' = ( \tau(k), \tau(k+1)) \tau = \tau (k, k+1)
$$
for some $k \in [1,N-1]$ with $\tau'(k) < \tau'(k +1)$ and
$\eta(\tau(k)) = \eta(\tau(k+1))$. Set $k_\bu := \tau_\bu \tau(k)$.

Assume that there exists an $N \times \ex$ matrix $\wt{B}_\tau$ with integral 
entries whose columns 
$b_\tau^l \in \Zset^N$, $l \in \ex$ satisfy 
\begin{equation}
\label{B-Om}
\chi_{\Ytil_\tau(b_\tau^l)}= 0, \quad
\Om_{\rb_\tau} ( b^l_\tau, e_j) = 0, \; \; \forall j \in [1,N], \; j \neq l
\quad \mbox{and} \quad
\Om_{\rb_\tau} (b^l_\tau, e_l) = \la_l^*
\end{equation}
for all $l \in \ex$. Then its principal part is skew-symmetrizable 
and the columns $b_{\tau'}^j \in \Zset^N$, $j \in \ex$ 
of the matrix $\mu_{k_\bu}(\Btil_\tau)$ satisfy  
\begin{equation}
\label{B-Om'}
\chi_{\Ytil_{\tau'}(b_{\tau'}^l)}= 0, \quad 
\Om_{\rb_{\tau'}} ( b^l_{\tau'}, e_j) = 0, \; \; \forall j \in [1,N], \; j \neq l
\quad \mbox{and} \quad
\Om_{\rb_{\tau'}} (b^l_{\tau'}, e_l) = \la_l^*
\end{equation}
for all $l \in \ex$. Furthermore,
\begin{equation}
\label{rtautau'}
\rb_{\tau'} = \mu_{k_\bu}(\rb_\tau)
\end{equation}
and
\begin{equation}
\label{bbg}
b_\tau^{k_\bu} = -b_{\tau'}^{k_\bu} = e_{p(k_\bu)} + e_{s(k_\bu)} - g,
\end{equation}
where $g \in \Znn^N$
is the vector from Proposition {\rm\ref{pcluster-tau-ind}}. 
\epr

\bre{diff}
The statements of \leref{ind-ab} and \prref{tau-seed-muta} have many similarities 
and their proofs use similar ideas. However, we note that there is a
principal difference between the two results. In the former case we have no 
mutation matrices to start with and we use \prref{cluster-tau-ind}(b) to gradually add
columns. In the latter case we already have a mutation matrix for one 
seed and construct a mutation matrix for another seed. 
\ere

\begin{proof}[Proof of Proposition {\rm\ref{ptau-seed-muta}}]
The fact that the principal part of $\wt{B}_\tau$ is skew-symmetrizable 
follows from \leref{alternative} and the condition \eqref{d-prop}, since $(\rb_\tau, \Btil_\tau)$ is a compatible pair by \eqref{B-Om}. 
The assumptions 
on $\wt{B}_\tau$ and \prref{cluster-tau-ind}(b) imply
$$
\chi_{\Ytil_\tau(b_\tau^{k_\bu})} = \chi_{\Ytil_\tau(e_{p(k_\bu)} + e_{s(k_\bu)} - g)} 
$$
and
$$
\Om_{\rb_\tau} (b^{k_\bu}_\tau, e_j) = \Om_{\rb_\tau}(e_{p(k_\bu)} + e_{s(k_\bu)} - g, e_j), 
\; \forall j \in [1,N].
$$ 

By \leref{unique}, $b^{k_\bu}_\tau = e_{p(k_\bu)} + e_{s(k_\bu)} - g$. The 
mutation formula for $\mu_{k_\bu}(\wt{B}_\tau)$ also gives that $b^{k_\bu}_{\tau'} = - b^{k_\bu}_\tau$,
so we obtain \eqref{bbg}. 
Analogously to the proof of \eqref{rtau'-mukbu}, one establishes \eqref{rtautau'}. 
Finally, all identities in \eqref{B-Om'}
follow from the general mutation facts in \prref{pair-mut}(b) and \leref{equivar-mut}.
\end{proof}

\begin{proof}[Proof of Theorem {\rm\ref{tcluster}(a)(b)} for all $\tau \in \Xi_N$]
Similarly to the proof of \leref{ind-ab}, 
the uniqueness statement in part (a) follows from \leref{unique}.
We will prove the existence statement in part (a) by an inductive argument on $\tau$. 
Once the existence of the matrix $\wt{B}_\tau$ with the stated properties 
is established, the fact that the principal part of $\wt{B}_\tau$ is 
skew-symmetrizable follows from \leref{alternative} and 
the condition \eqref{d-prop}, and the compatibility of the pair $(\rb_\tau, \Btil_\tau)$ follows from \eqref{linear-eq} and \eqref{compscalars}.
Hence, $(\ybtil_\tau, \wt{B}_\tau)$ 
is a seed and this yields part (b) of the theorem
for the given $\tau \in \Xi_N$.

For the existence statement in part (a)
we fix $\tau \in \Xi_N$. By \coref{steps}(a), there exists a sequence
$\tau_0 = \id, \tau_1, \ldots, \tau_n = \tau$ in $\Xi_N$ with the property that
for all $l \in [1,n]$,
$$
\tau_l = ( \tau_{l-1}(k_l), \tau_{l-1} (k_l+1)) \tau_{l-1} = \tau_{l-1}(k_l, k_l+1)
$$
for some $k_l \in [1,N-1]$ such that $\tau_{l-1}(k_l) < \tau_{l-1}(k_l +1)$.
In \S \ref{recursive} we established the validity of \thref{cluster}(a) for the identity element 
of $S_N$. By induction on $l$ we prove the validity of \thref{cluster}(a) for $\tau_l$. 
If $\eta(\tau_{l-1}(k_l)) \neq \eta(\tau_l(k_l))$, then \prref{cluster-tau-ind}(a) implies 
that $\Ytil_{\tau_l} = \Ytil_{\tau_{l-1}}$ and we can choose $\wt{B}_{\tau_l} = \wt{B}_{\tau_{l-1}}$. If $\eta(\tau_{l-1}(k_l)) = \eta(\tau_l(k_l))$, then \prref{tau-seed-muta} proves that 
the validity of \thref{cluster}(a) 
for $\tau_{l-1}$ implies the validity of \thref{cluster}(a) 
for $\tau_l$. In this case $\wt{B}_{\tau_l} := \mu_{(k_l)_\bu} (\wt{B}_{\tau_{l-1}})$, 
where $(k_l)_\bu = (\tau_{k_l})_\bu  \tau_{k_l} (k_l)$.
This completes the proof of \thref{cluster}(a)(b). 
\end{proof}

\begin{proof}[Proof of Theorem {\rm\ref{tcluster}(c)}] The one-step mutation statement 
in part (c) of \thref{cluster} and \coref{steps}(a) imply that all seeds 
associated to the elements of $\Xi_N$ are mutation-equivalent 
to each other.

In the rest we prove the one-step mutation statement in part (c) 
of the theorem.
If $\eta(\tau(k)) \neq \eta(\tau(k+1))$, then the 
statement follows from \prref{cluster-tau-ind}(a).

Now let $\eta(\tau(k)) = \eta(\tau(k+1))$. We have that either $\tau(k) < \tau(k+1)$ or
$\tau'(k) = \tau(k+1) < \tau'(k+1) = \tau(k)$. In the first case we 
apply \prref{tau-seed-muta} to the pair $(\tau, \tau')$ and 
in the second case to the pair $(\tau', \tau)$. The  one-step mutation statement in \thref{cluster}(c) follows from this, the uniqueness statement in part (a) of the theorem,
and the involutivity of mutations of seeds.
\end{proof}

\subsection{Completion of the proof of \thref{cluster}} 
\label{complete.main}

Recall the setting of \S \ref{intersectloc}. For all $\tau \in \Xi_N$, define the multiplicative subsets
$$
E_\tau := \Big\{ \alpha \Ytil_\tau(f) \bigm| \alpha \in \kx, \; f \in \sum_{j \in \ex} \Zset_{\geq 0} e_j \Big\} \subset R.
$$
In view of \leref{taubutau}(a) and the definition of $\Ytil_\tau$, we see that $E_\tau$ is generated (as a multiplicative set) by
$$
\kx \sqcup \{ y_{\tau,k} \mid k \in [1,N], \; s_\tau(k) \ne +\infty \},
$$
matching the definition used in \S \ref{intersectloc}.

\begin{proof}[Proof of Theorem {\rm\ref{tcluster}(d)--(f)}] By \thref{cluster}(c), for all $\tau \in \Xi_N$ the seeds $(\ybtil_\tau, \wt{B}_\tau)$ are mutation-equivalent to each other. 
For each $j \in [1,N]$, we have $\tau_{j,j} \in \Ga_N$ with $\tau_{j,j}(1) = j$, and so $\ol{y}_{\tau_{j,j},1} = \ol{y}_{[j,j]} = x_j$ by \thref{sym-Pprime}. Thus,
each generator $x_j$ of $R$ is a cluster or frozen 
variable for a seed associated to some $\tau \in \Ga_N$. Hence,
\begin{equation}
\label{inc1}
R \subseteq \AA(\ybtil, \wt{B})_\KK.
\end{equation}
The Laurent phenomenon (\thref{phenom}) implies that 
\begin{equation}
\label{inc2}
\AA(\ybtil, \wt{B})_\KK \subseteq \UU (\ybtil, \wt{B}, \varnothing)_\KK.
\end{equation}

For all $\tau \in \Xi_N$, we have $\UU (\ybtil, \wt{B}, \varnothing)_\KK = \UU (\ybtil_\tau, \wt{B}_\tau, \varnothing)_\KK$ in view of part (c) of the theorem and \thref{U=Abar}. Moreover, $\UU (\ybtil_\tau, \wt{B}_\tau, \varnothing)_\KK \subseteq \TT\AA (\ybtil_\tau, \wt{B}_\tau, \varnothing)_\KK \subseteq  R[E_\tau^{-1}]$ because $y_{\tau,i} \in R$ for all $i \in [1,N]$ and $y_{\tau, (\tau_\bu \tau)^{-1} (k)} \in E_\tau$ for all $k \in \ex$. Consequently,
\begin{equation}
\label{inc3}
\UU (\ybtil, \wt{B}, \varnothing)_\KK \subseteq \bigcap_{\tau \in \Ga_N} R[E_\tau^{-1}].
\end{equation} 
Combining the inclusions \eqref{inc1}, \eqref{inc2}, \eqref{inc3} and \thref{PCGLalmostcluster}(d) leads to 
$$
R \subseteq \AA(\ybtil, \wt{B})_\KK \subseteq \UU (\ybtil, \wt{B}, \varnothing)_\KK
\subseteq \bigcap_{\tau \in \Ga_N} R[E_\tau^{-1}] = R,
$$
which establishes all equalities in \thref{cluster}(d).

For part (f) we obtain inclusions
\begin{equation}
\label{chain}
\begin{aligned}
R[y_k^{-1} \mid k \in \inv]  &\subseteq \AA(\ybtil, \wt{B}, \inv)_\KK  \\
&\subseteq \UU (\ybtil, \wt{B}, \inv)_\KK
\subseteq \bigcap_{\tau \in \Ga_N} R[y_k^{-1} \mid k \in \inv][E_\tau^{-1}]
\end{aligned}
\end{equation}
in the same way as above. \thref{PCGLalmostcluster}(e) and \eqref{chain} imply the validity of the first statement in
part (f) of the theorem. 

It remains to verify the Poisson-compatibility statements in parts (e) and (f).

By part (b) of the theorem, $(\rb_\tau, \Btil_\tau)$ is a compatible pair, for each $\tau \in \Xi_N$. In particular, $(\qb, \Btil) = (\rb_\id, \Btil_\id)$ is a compatible pair. Eq.~\eqref{bracketykyj} in \prref{PtorusR} says that
\begin{equation}
\label{brack-Omq}
\{ y_l, y_j\} = \Om_\qb(e_l,e_j) y_l y_j , \; \; \forall l,j \in [1,N].
\end{equation}
Consequently, \thref{P-cluster-main} implies that the cluster algebras $\AA(\ybtil, \Btil, \inv)_\KK$ are Poisson-compatible.

Given any $\tau \in \Xi_N$, we again appeal to \coref{steps}(a) to obtain a sequence
$\tau_0 = \id, \tau_1, \ldots,\allowbreak \tau_n = \tau$ in $\Xi_N$ with the property that
for all $i \in [1,n]$,
$$
\tau_i = ( \tau_{i-1}(k_i), \tau_{i-1} (k_i+1)) \tau_{i-1} = \tau_{i-1}(k_i, k_i+1)
$$
for some $k_i \in [1,N-1]$ such that $\tau_{i-1}(k_i) < \tau_{i-1}(k_i +1)$. We prove by induction on $i \in [1,n]$ that
\begin{equation}
\label{induct-compat}
\{ \Ytil_{\tau_i}(e_l), \Ytil_{\tau_i}(e_j) \} = \Om_{\rb_{\tau_i}}(e_l, e_j) \Ytil_{\tau_i}(e_l) \Ytil_{\tau_i}(e_j), \; \; \forall l,j \in [1,N].
\end{equation}
The case $i=0$ is just \eqref{brack-Omq}. Now let $i>0$. If $\eta(\tau_{i-1}(k_i)) \ne \eta(\tau_i(k_i))$, then $\rb_{\tau_i} = \rb_{\tau_{i-1}}$ and $\Ytil_{\tau_i} = \Ytil_{\tau_{i-1}}$ by \prref{cluster-tau-ind}(a), and so the case $i$ of \eqref{induct-compat} is immediate from the case $i-1$. If $\eta(\tau_{i-1}(k_i)) = \eta(\tau_i(k_i))$, then $\rb_{\tau_i} = \mu_{(k_i)_\bu}(\rb_{\tau_{i-1}})$ by \prref{tau-seed-muta} and $\Ytil_{\tau_i} = \mu_{(k_i)_\bu}(\Ytil_{\tau_{i-1}})$ by statement (c) of \thref{cluster}. In this case, \prref{brack-mut} shows that the case $i-1$ of \eqref{induct-compat} implies the case $i$.
\end{proof}

The chain of embeddings \eqref{chain} and \thref{PCGLalmostcluster}(e) also imply the 
following description of the upper cluster algebra in \thref{cluster}
as a finite intersection of mixed polynomial-Laurent polynomial algebras.

\bco{up-cl} In the setting of Theorem {\rm\ref{tcluster}},
$$
\UU (\ybtil, \wt{B}, \inv)_\KK =  \bigcap_{\tau \in \De} \TT\AA(\xbtil_\tau, \Btil_\tau, \inv)
$$
for every subset $\De$ of $\Ga_N$ which is an interval with respect to the 
linear ordering \eqref{sequence} and has the property that for each $k \in \ex$ 
there exist two consecutive elements $\tau=\tau_{i, j-1} \prec \tau' = \tau_{i,j}$ 
of $\De$ such that $\eta(i) = \eta(j)$ and $k = \tau_\bu \tau(j-i)$, recall \S {\rm\ref{recursive}}.

In particular, this property holds for $\De = \Ga_N$. 
\eco



\end{document}